\documentclass[11pt,fleqn]{article}
\usepackage[top=1.1in, left=1in, right=1in, bottom=1.1in]{geometry}
\usepackage{amsmath,amssymb,amsthm,amsfonts,verbatim}
\usepackage{enumerate}
\usepackage{microtype}
\usepackage{url}
\usepackage{algpseudocode}
\usepackage{tikz}
\usetikzlibrary{arrows}
\usepackage{amsmath,amssymb}
\usepackage[justification=centering]{caption}
\usepackage[justification=centering]{subcaption}
\usepackage[english]{babel}
\usepackage{array}
\usepackage{psfrag}
\usepackage{graphicx}
\usepackage{color}
\usepackage{setspace}
\usepackage{hyperref}
\usepackage{algorithm}
\usepackage{algpseudocode}
\usepackage[utf8]{inputenc}
\usepackage{cleveref}
\usepackage{mathrsfs} 
\usepackage{authblk}
\usepackage{threeparttable}
\usepackage[titletoc,title]{appendix}
\usepackage{xcolor}
\usepackage{xfrac}
\usepackage{caption}
\usepackage[font={small}]{caption}
\usepackage{adjustbox}

\newtheorem{proposition}{Proposition}[section]

\newtheorem{theorem}{Theorem}[section]

\newtheorem{definition}[theorem]{\bf Definition}

\definecolor{red}{rgb}{1,0.2,0.2}

\newtheorem{ARrule}{Accept/Reject Rule}

\def\F{{\mathcal F}}

\def\O{{\mathcal O}}

\def\T{{\mathcal T}}
\def\R{{\mathbb R}}

\def\E{{\mathbb E}}
\def\R{{\mathbb R}}

\def\ds{\displaystyle}
\newcommand{\norm}[1]{\left\lVert#1\right\rVert}

\def\hs{\hspace{2cm}}

\def\e{\epsilon}

\def\s0t{\sup_{t \in [0,T]}}

\def\ds{\displaystyle}

\def\beq{\begin{equation}}
\def\eeq{\end{equation}}
\def\barr{\begin{array}}
\def\earr{\end{array}}

\def\e{{\epsilon}}
\def\F{{\mathcal F}}

\def\O{{\mathcal O}}

\def\E{{\mathbb E}}

\def\ds{\displaystyle}
\def\Z20{{\mathbb{Z}^2_0}}
\def\T2{\mathbb{T}^2}
\def\hs{\hspace{2mm}}

\hypersetup{colorlinks=false,linkbordercolor=red,linkcolor=blue,pdfborderstyle={/S/U/W 1}}

\makeatletter
\newcommand*{\rom}[1]{\expandafter\@slowromancap\romannumeral #1@} 
\makeatother

\hypersetup{colorlinks=true,linkbordercolor=red,linkcolor=blue,pdfborderstyle={/S/U/W 1}}

\def\R{\mathbb{R}}

\title{An efficient jet marcher for computing the quasipotential for 2D SDEs}%
\author[1]{Nicholas Paskal\thanks{npaskal@umd.edu}}
\author[1]{Maria K. Cameron\thanks{mariakc@umd.edu}}
\affil[1]{\small{Department of Mathematics, University of Maryland, College Park, MD 20742, USA}}

\begin{document}
\maketitle

\abstract{
We present a new algorithm, the efficient jet marching method (EJM), for computing the quasipotential and its gradient for two-dimensional SDEs. The quasipotential is a potential-like function for nongradient SDEs that gives asymptotic estimates for the invariant probability measure, expected escape times from basins of attractors, and maximum likelihood escape paths. The quasipotential is a solution to an optimal control problem with an anisotropic cost function which can be solved for numerically via Dijkstra-like label-setting methods. Previous Dijkstra-like quasipotential solvers have displayed in general 1st order accuracy in the mesh spacing. However, by utilizing higher order interpolations of the quasipotential as well as more accurate approximations of the minimum action paths (MAPs), EJM achieves second-order accuracy for the quasipotential and nearly second-order for its gradient. Moreover, by using targeted search neighborhoods for the fastest characteristics following the ideas of  Mirebeau, EJM also enjoys a reduction in computation time. This highly accurate solver enables us to compute the prefactor for the WKB approximation for the invariant probability measure and the Bouchet-Reygner sharp estimate for the expected escape time for the Maier-Stein SDE. Our codes are available on GitHub.
}

%
%
\section{Introduction}
\label{sec:quasi:intro}
Consider the following stochastic differential equations in $\mathbb{R}^n$:
\begin{equation}\label{eq_SDE}
	dX_t = b(X_t)dt + \sqrt{\e} dW_t, \hs X_0 = x_0 \in \mathbb{R}^n,
\end{equation}
where $b: \mathbb{R}^n \to \mathbb{R}^n$, $W_t$ is a standard Brownian motion in $\mathbb{R}^n$ and $\e>0$ is a positive parameter. When $\e>0$ is small, the solution $X_t$ typically follows very closely the trajectories of the deterministic system 
\begin{equation}\label{eq_ODE} 
dY_t = b(Y_t)dt, \hs Y_0 = x_0.
\end{equation}
If, in addition, the vector field $b:\mathbb{R}^n \to \mathbb{R}^n$ admits a globally attracting equilibrium at $\O \in \mathbb{R}^n$, the solution $X_t$ will spend exponentially long time periods in $\e^{-1}$ near the attractor $\O$ before experiencing $\O(1)$ size deviations from $\O$. 

The long-time behavior can be characterised by the invariant probability measure $\mu$ of $\eqref{eq_SDE}$. A key result of the Freidlin-Wentzell theory \cite{freidlin12} is that the measure $\mu$ is approximately \textit{Gibbsian} for small $\e > 0$, in the sense that
\begin{equation}\label{eq_invar}
\frac{d\mu}{dx}(x) \asymp \exp\Big(-\frac{U(x)}{\e} \Big),
\end{equation}
where the symbol $\asymp$ denotes logarithmic equivalence in the $\e \to 0$ limit. The exponent $U:\mathbb{R}^n \to \mathbb{R}$ in $\eqref{eq_invar}$ is the so called \textit{quasipotential}, which is the primary object of interest in this paper. 

\subsection{An overview}
The quasipotential $U(x)$, whose definition is given in Section \ref{sec_background}, can be interpreted as the minimum cost of moving from the attractor $\O$ to the point $x$, as measured by an appropriate action functional. Thus, a logical way to compute $U(x)$ for a given $x$ is to conduct this minimization numerically. Indeed, one can construct a numerical version of the action function via quadrature and then perform a high-dimensional minimization over a suitably rich path space. Since methods of this form invariably return the minimum action paths (MAP) themselves, we refer to them as \textit{path-based} methods. 

In practice, techniques conducting the minimizations of two types of action functionals are used. The Minimum Action Method (MAM) \cite{e04} and Adaptive Minimum Action Method (AMAM) \cite{zhou08} find MAPs by minimizing the Freidlin-Wentzell action, while the Geometric Minimum Action Method (GMAM) \cite{heymann08} and the recent Ritz method \cite{kikuchi20} do so by minimizing the geometric action. By design, these techniques are applicable to both finite and infinite dimensional problems. For example, the GMAM was used to find the MAP between two stationary solutions of a 2D reaction-diffusion partial differential equation in \cite{heymann08} and to find the MAP between two solitary waves of different amplitudes of nonlinear Schr\"{o}dinger equations in \cite{poppe18}. 

An important advantage of path-based methods is that they are computationally cheap
and suitable for any dimension. However, they have a number of key drawbacks. 
First, they allow for computation of the quasipotential along a \textit{single} MAP, 
but do not provide a mechanism for efficient computation of the quasipotential over an entire region of space. 
In particular, they are ill-suited for the task of identifying quasipotential minima. 
Second, path-based techniques work well if the MAP is relatively simple, but
their convergence tends to stall if the MAP exhibits spiraling or other complicated behavior. These issues become even more severe if the MAP has infinite length.
Finally, the MAPs obtained by these methods are \textit{biased} by the initial guess, typically taken to be straight line segments, 
and may only be \textit{local} minimizers of the action that lead to incorrect estimates of the quasipotential. 
While various attempts have been undertaken to address these issues
\cite{tao18,lin19}, none of the proposed solutions have become commonly
used due to complexity and lack of robustness.

\begin{figure}[h]
  \includegraphics[width=0.9\textwidth]{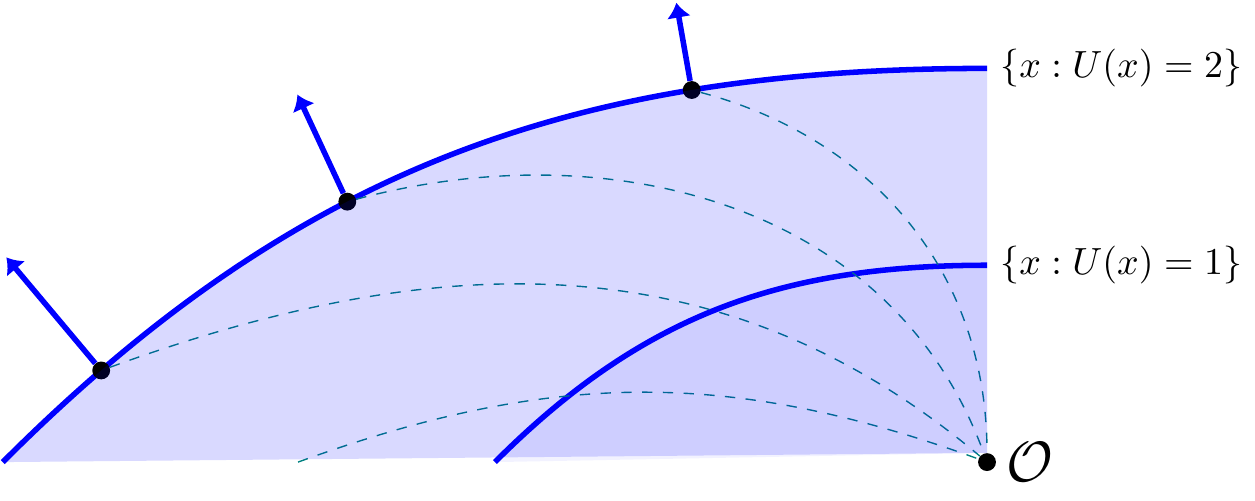}
\caption{The quasipotential problem as a wavefront propagation problem.}
\label{fig:quasi:wave}       
\end{figure}

This article focuses on \emph{Dijkstra-like} \textit{mesh-based} quasipotential solvers, descendants of Dijkstra's algorithm for identifying the shortest path in a network \cite{dijkstra59}, and, more specifically, of Sethian's Fast Marching Method for solving the eikonal equation \cite{sethian96,sethian99}. 
These methods treat the quasipotential $U(x)$ as the first hitting time at location $x$ of a fictitious wavefront propagating outwards from the point attractor $\O$, 
as illustrated in Figure \ref{fig:quasi:wave}. 
In this analogy, the level set $\{x:U(x) = T\}$ represents the state of the wavefront at ``time'' $T$. The MAP passing through a point $x$ can be seen 
as representing the trajectory of the particular particle that was the first to reach $x$.
The advantages of methods of this type are 
that they 
$(i)$ find the quasipotential in a whole region of space and allow us to estimate the invariant probability measure, 
$(ii)$ require no initial guess,
and 
$(iii)$ allow for reconstruction of the MAP passing through any given point $x$.
Their primary disadvantage is that they suffer from the curse of dimensionality and hence are limited to two or three dimensions.

 All current mesh-based quasipotential solvers  \cite{cameron12,cameron18,dahiya18,yang19,lorenz63}
inherit the MAP search algorithm from Sethian's and Vladimirsky's 
\emph{ordered upwind method} (OUM) for solving static Hamilton-Jacobi equations with \emph{bounded} anisotropic speed functions \cite{sethian01,sethian03}. 
Roughly speaking, this involves an exhaustive search over rather large neighborhoods, with significant 
rationalizations in later solvers called \emph{ordered line integral methods} (OLIMs) \cite{cameron18,dahiya18,yang19,lorenz63}.
An important aspect of the problem of computing the quasipotential is that its anisotropic 
speed function is \emph{unbounded}, and the optimal choice of the radius for the search neighborhood remains an important issue; 
only rules of thumb based on detailed studies of particular systems have been proposed so far \cite{cameron18,dahiya18,yang19}.  
Like OUM, the OLIMs are only first-order accurate, though their error constants are by orders of magnitude smaller than in the adjusted OUM \cite{cameron12}.

\subsection{A brief summary of main results}
In this work, we introduce a new solver,  the \emph{efficient jet marcher} (EJM), that exhibits second-order convergence and uses a targeted search algorithm for segments of the MAPs. The aforementioned issue of choosing the radius for the search neighborhood is completely eliminated.
This solver serves as a pilot project whose purpose is to explore new technical solutions for enhancing accuracy and reducing computation time. 
Therefore, we limit ourselves to the simplest setting, two-dimensional SDEs with point attractors, 
and focus on the development of a highly accurate, fast, and robust numerical algorithm.

The \textit{first} key feature is the use of higher-order approximations of (a) the quasipotential between mesh points and (b) local segments of MAPs. Concerning (a), our solver propagates the quasipotential {\it and} its gradient by taking advantage of the geometric relationship between the gradient of the quasipotential, and the direction of the MAP.
This enables us to use Hermite cubic interpolation of $U$ between mesh points rather than linear interpolation as done in \cite{cameron12,cameron18,dahiya18,yang19}. The idea of using Hermite interpolation originally comes from the jet scheme for solving an advection equation \cite{nave10}. Concerning (b), our solver approximates MAP segments using cubic curves, while all previous quasipotential solvers have approximated them with straight line segments. We remark that a similar family of jet marching methods (JMMs) for isotropic eikonal equations was recently introduced in \cite{potter20}. In comparison with the eikonal JMMs, our EJM for computing the quasipotential is significantly more nuanced and complicated as the optimal control problem for the quasipotential is notably harder than the eikonal equation.
The improved approximation scheme provides our method with an $\O(h^2)$ convergence rate for the eikonal and $\sim\O(h^{1.8})$ for its gradient with respect to the mesh spacing $h$, as we demonstrate empirically on a collection of test problems.
To our knowledge, this is the first such quasipotential solver with second-order in $h$ accuracy. 

%
 
The aforementioned improvement in accuracy does come at the cost of an increase in runtime due to the more computation-intensive MAP approximations. However, this increase in runtimes is offset by the algorithm's \textit{second} key feature: implementation of pre-computed anisotropic stencils, inspired by the Anisotropic Stencil Refinement (ASR) algorithm of Mirebeau \cite{mirebeau14}. The anisotropic unbounded speed function in the quasipotential problem requires the searching of very large ``neighborhoods'' of mesh points in the search for MAPs. Mirebeau's ASR algorithm pre-computes smaller, more targeted neighborhoods, which results in a significantly smaller number of total MAP searches at a comparable overall accuracy. Adoption of a modified version of these ideas provided EJM with a significant reduction in computation time. 

{ The EJM is implemented in C++. Our codes are available on GitHub \cite{EJMpaskal}.}

We compare the performance of the EJM algorithm with other mesh-based quasipotential solvers on a variety of different drift fields $b$ and demonstrate a dramatic gain in efficacy.

The highly accurate numerical solution not only for the quasipotential but also for its gradient allows us to solve the transport equation for the WKB prefactor { that contains second derivatives of the quasipotential. This yields} a sharp estimate for the invariant probability measure. We do so for the Maier-Stein SDE \cite{MS1993} at two parameter values: $\beta=3$ where the SDE is nongradient but the quasipotential is smooth, and $\beta=10$ where the gradient of the quasipotential is discontinuous along an interval of the $x$-axis. As predicted in \cite{maier96}, the WKB approximation yields a good estimate for the invariant measure at $\beta=3$ but not at $\beta=10$ where the prefactor blows up. In addition, we use our numerical solutions to evaluate the Bouchet-Reygner formula for the expected escape time from the basin of the left attractor of the Maier-Stein SDE at various values of $\epsilon$. We compare our estimates for the invariant measure and the expected escape time to those obtained by means of the transition path theory \cite{eve2006} and observe a good agreement.

{ As you will see, the high accuracy and efficiency of the EJM come at the cost of its complexity. Pros and cons of the EJM and its possible extensions  to three-dimensional SDEs and SDEs with geometric noise will be discussed in Section \ref{sec:conclusions}. }

The rest of the paper is organized as follows. Background on the Freidlin-Wentzell theory 
and the optimal control problem for the quasipotential is given in Section \ref{sec_background}.
Relevant numerical algorithms that have motivated our solver are described in Section \ref{sec:quasi:dijkstra}.
The new solver, the efficient jet marcher (EJM), is detailed in Section \ref{sec:quasi:hasr}.
A full analysis of the accuracy-speed trade-off is provided in Section \ref{sec:quasi:results}.
The application to the Maier-Stein SDE is found in Section \ref{sec:MS}. 
{ The results are summarized and perspectives are discussed} in Section \ref{sec:conclusions}. 
Some calculations necessary for implementation of the EJM are placed in Appendix A. Some details on the transition path theory are provided in Appendix B.

%
%

\section{Background}
\label{sec_background}
Let $X^\e(t)$ be the solution to SDE $\eqref{eq_SDE}$ in $\R^2$ for a fixed $\e > 0$ where $W_t$ is the standard Brownian motion in $\R^2$. We assume that the drift field $b:\R^2\to\R^2$ is a smooth vector field that admits a stable attracting equilibrium at the point $\O$. Suppose that the set $D \subset \R^2$ is the basin of attraction of $\O$ for the deterministic system $\eqref{eq_ODE}$. When $\e$ is small, the process $X_t$ will remain close to the attractor $\O$ for large periods of time until a sufficiently strong realization of the driving noise $W_t$ drives the $X_t$ farther away from $\O$. We are interested in quantifying these deviations from the deterministic dynamics. 
In this section, we provide a primer on Freidlin-Wentzell theory to define the quasipotential, describe its importance in obtaining asymptotic probability estimates, 
and formulate the optimal control problem for the quasipotential.

\subsection{Freidlin-Wentzell theory}
The $\e \downarrow 0$ dynamics described by the Freidlin-Wentzell theory are driven by a least action principle. Indeed, for any fixed $T> 0$ the asymptotic \textit{likelihood} of the solution $X(t)$ following a given trajectory $ \phi \in C([0,T];\R^2)$ is controlled by the Freidlin-Wentzell action functional $S_T:C([0,T];\R^2 \to [0,+\infty]$ given by
\begin{equation*}
	S_T(\phi) := 
	\begin{cases}
	\ds \frac{1}{2} \int_0^T \norm{\dot{\phi}_t - b(\phi_t)}^2 dt, & \text{if $\phi$ is absolutely continuous},
	\\  +\infty, &\text{otherwise}.
	\end{cases}
\end{equation*}
Formally speaking, we say that the families of solutions $\{X^\e\}_{\e >0}$ of equation $\eqref{eq_SDE}$ satisfy a large deviations principle on $C([0,T];\R^2)$ with action function $S_T$. Effectively, this means that the probability of $X_\e$ lying in a tube around $\phi$ decays at the exponential rate $\exp(-S_T(\phi)/\e)$ as $\e \downarrow 0$.

The \textit{quasipotential} $U:D \to [0,+\infty]$, defined with respect to the stable attractor $\O$ and for $x$ in its domain of attraction, is the cost of the \textit{cheapest} path, as measured by the Freidlin-Wentzell action functional, that starts at $\O$ and terminates at $x$. That is,
\begin{equation}\label{eq_quasi}
U(x) := \inf \Big\{ S_T(\phi): T > 0, \ \ \phi \in C([0,T]; \R^n), \ \ \phi(0) = \O, \ \ \phi(T) =x \Big\}.
\end{equation}

The function $U$ determines much of the behavior of equation $\eqref{eq_SDE}$ in the asymptotic long-time limit for small $\e>0$. As mentioned in the Introduction, the leading order term of the stationary measure $\mu_\e$ of $\eqref{eq_SDE}$ behaves like equation $\eqref{eq_invar}$ (Theorem 4.3 of \cite{freidlin12}). Via this asymptotics on the stationary measure, the quasipotential $U$ can be used to obtain estimates on the mean first exit times from bounded domains. 

Indeed, consider the problem of the first exit of the solution $X_t^\e$ from a bounded domain $\tilde{D}  \subset D$. Provided that $b$ does not grow too steeply, one can show that the solution $X_t^\e$ starting at $\O$ will almost surely exit $\tilde{D}$ eventually due to the noise. In fact, the mean first exit $\tau^\e:= \inf \{t > 0: X_t^\e \notin\tilde{ D}\}$ can be shown to obey the inverse asymptotic growth rate of the stationary measures. Namely,
\begin{equation*}
	\mathbb{E} \tau^\e \asymp \exp \Big(\frac{\inf_{x \in \partial \tilde{D}} U(x)}{\e} \Big).
\end{equation*}
In effect, the average first exit time is entirely controlled by the cost of the cheapest path from $\O$ to the boundary $\partial \tilde{D}$. In fact, if there exists a unique minimizer $x^* \in \partial \tilde{D}$, then as $\e \downarrow 0$ the actual escape trajectory concentrates with overwhelming probability closer and closer to the minimizing path $\phi^*$ from $\O$ to $x^*$.

{\color{black} A similar analysis holds in the case where the drift admits multiple attractors. In such systems, an important phenomenon is the probabilistic switching from one \textit{metastable} attractor to another. For stochastic differential systems, the leading order calculation of mean transition rates can be quantified in the same manner. This is often a topic of interest in application. For instance, in \cite{Nolting16}, they utilize this framework to understand transitions between different stable population configurations in an ecological model.  }

\subsection{The optimal control problem for the quasipotential}
\label{sec:optcontrol}
In actuality, the infimum in definition \ref{eq_quasi} over $T> 0$ is not achieved. Rather, the corresponding minimum action path departs $\O$ with infinitesimally small speed and the minimizer of the action exists only as a $T \to \infty$ limit. Nonetheless, the graph of the path does exist in practice. 
To identify this minimum action path, it is convenient to reparametrize and to work with a modified action functional called the geometric action functional
\begin{equation}\label{eq_action_geometric}
	\tilde{S}_L(\varphi) := 
	\begin{cases}
	\ds \int_0^L \Big[ \norm{\dot{\varphi}_r} \norm{b(\varphi_r)}- \dot{\varphi_r} \cdot b(\varphi_r) \Big] dr, & \text{if $\varphi$ is absolutely continuous},
	\\  +\infty, &\text{otherwise}.
	\end{cases}
\end{equation}
By a simple reparametrization argument (see Section 3 in \cite{cameron12}), 
it can be shown that the definition of the quasipotential $U(x)$ with $S$ replaced by $\tilde{S}$ is equivalent. 
{\color{black} The paths $\phi$ and $\varphi$ relate via $\dot{\phi} = \|b(\varphi)\|\dot{\varphi}$.}
Therefore, an alternative definition which we stick to throughout the remainder is
\begin{equation}
\label{eq_quasi_geometric}
	U(x)  := \inf \Big\{ \tilde{S}_L(\varphi): \  L > 0, \ \varphi \in C([0,L];\R^2), \ \varphi_0 = \O, \ \varphi_L = x \Big\}.
\end{equation}
Note that the integrand of the geometric action can be written as $s(\varphi,\dot{\varphi})\|\dot{\varphi}\|$, where
\begin{equation}\label{eq_slowness}
	s(x,v) = \norm{b(x)} - b(x)\cdot \frac{v}{\norm{v}},\quad {\rm where}  \quad v\equiv \dot{\varphi}.
\end{equation}
The quantity $s(x,v)$ can be interpreted as the anisotropic slowness function at point $x$ in the direction $v$. 
Note that $s(x,v)$ is bounded from above by $2\|b(x)\|$ which makes $U(x)$ Lipschitz-continuous in any bounded region in $\R^2$.
On the other hand, $s(x,v)$ is zero at any equilibrium point and whenever the directions 
of $v$ and  $b(x)$ coincide. As a result, the speed function $f(x,v)\equiv 1/s(x,v)$ is unbounded from above.

Using  definition \eqref{eq_quasi_geometric} of the quasipotential, the infimum in $U$ may be achieved, although this is not guaranteed as infinite length minimum action paths are possible. Moreover, in \cite{cameron12} it is shown \texttt{?}hat the direction of the minimizing path $\varphi \in C([0,L];\R^2)$ passing through $x$ can be related to the gradient of the normal to the level sets of $U$ by the relation
\begin{equation}
\label{eq:gradUdir}
	\nabla U(\varphi) = \norm{b(\varphi)}\dot{\varphi} - b(\varphi).
\end{equation}

The quasipotential allows us to decompose the drift field $b(x)$ into a sum of potential and rotational components:
\begin{equation}
\label{eq:ortdec}
b(x) = -\frac{1}{2}\nabla U(x) + l(x),\quad{\rm where}\quad l(x):=b(x) +\frac{1}{2}\nabla U(x).
\end{equation}
Both, the trajectories of $dx = b(x)dt$ and the MAPs follow the rotational component $l(x)$. However, at the same time, 
the trajectories descend down the potential component $-\tfrac{1}{2}U(x)$, while the MAPs climb up it, as one can see by recasting \eqref{eq:gradUdir} as 
$\dot{\phi} = b(\phi)+\nabla U(\phi) \equiv \tfrac{1}{2}\nabla U(\phi) + l(\phi)$.

We can investigate the invariant measure further by considering again the Fokker-Planck equation with the WKB ansatz 
that the density takes the form of Gibbs density multiplied by a subexponential pre-factor. Namely, we take
\begin{equation}
\label{eq:WKB}
p^\e(x) = C(x) \exp \Big( -\frac{U(x)}{\e}  \Big).
\end{equation}
Justification for the WKB approximation is given in e.g. \cite{talkner87,maier96}. 
Plugging \eqref{eq:WKB} into the Fokker-Planck equation and grouping terms by order in $\e$ yields
the Hamilton-Jacobi-Bellman equation  for the quasipotential
\begin{equation}
 \label{eq_HJB}
b \cdot \nabla U + \frac{1}{2} \norm{\nabla U}^2 = 0,
 \end{equation}
and the transport equation for the leading order term in the exponential prefactor
\begin{equation}
\label{eq:transport}
(\nabla \cdot b +\frac{1}{2}\Delta U) C + (b+\nabla U) \cdot \nabla C = 0.
\end{equation}

In general, the quasipotential is not continuously differentiable. 
At worst, it is  Lipschitz continuous \cite{freidlin12} and the Lebesgue measure of the set of points for which it is not differentiable is $0$. 
Moreover, it can be shown that when $U$ is continuously differentiable, it is a classical solution to \eqref{eq_HJB} endowed with the boundary condition  $U(\O) = 0$, 
and a viscosity solution otherwise. 

Furthermore, even if the classical solution to \eqref{eq_HJB} with the boundary condition  $U(\O) = 0$ exists, it is not unique. 
There are at least two solutions, the quasipotential defined by \eqref{eq_quasi_geometric}, and $U\equiv 0$. 
In general, one can check that there are as many solutions as there are invariant subspaces for the linearized at $\O$ drift field.
Therefore, we solve the optimal control problem \eqref{eq_quasi_geometric}  rather that the boundary-value problem for \eqref{eq_HJB}. 

%
%

\section{Mesh-based label-setting algorithms}
\label{sec:quasi:dijkstra}
In this section, we describe the structure of Dijkstra-like solvers of anisotropic eikonal equations. We also provide a description of Mirebeau's design for the anisotropic (ASR) algorithm \cite{mirebeau14}. 

\subsection{Dijkstra-like eikonal solvers}
Consider the geometric action expressed in the line integral form
\begin{equation}\label{eq:quasi:action_generic}
	\tilde{S}_L(\varphi) = \int_0^L s(\varphi_r,\dot{\varphi}_r) \norm{\dot{\varphi}_r}dr \equiv \int_0^L \frac{\norm{\dot{\varphi}_r}}{f(\varphi_r,\dot{\varphi}_r) }dr.
\end{equation}
Written in this way, $\tilde{S}_L(\varphi)$ gives the cost of trajectory $\varphi$ using the instantaneous anisotropic cost function $s(x,v)$. If we are to proceed with the wavefront analogy introduced in Section \ref{sec:quasi:intro}, we can interpret $f(x,v)$ as the instantaneous speed of particles moving in the $v$ direction at point $x$. At each point along the wavefront, one can envision more particles being spawned and sent out in all directions with speed $f(x,v)$. For this reason, we call the function $f(x,v)$ the \textit{speed} function and $s(x,v)$ the \textit{slowness} function. The quasipotential $U(x)$ will then measure the quickest time for  one of these particles to hit the point $x$, while the MAP passing through $x$ will trace that particle's path back to its origin at $\O$. 

An alternative perspective can be taken by expressing the HJB equation $\eqref{eq_HJB}$ for the quasipotential in the form
\begin{equation}\label{eq:quasi:eikonal}
	 F(x,\hat{n})\norm{\nabla U(x)}= 1, \hs U(\O) = 0,
\end{equation}
where $\ds F(x,\hat{n}) = \frac{1}{-2b(x)\cdot \hat{n}(x)}$ and $\ds \hat{n}(x) = \frac{\nabla U(x)}{\norm{\nabla U(x)}}$ is the outward pointing normal vector to the level set of the quasipotential. Here, the quantity $F(x,\hat{n}) \geq 0$ can be interpreted as the speed at which the wavefront expands outward in the normal direction. 

For a general $F$, equation $\eqref{eq:quasi:eikonal}$ is referred to as an anisotropic eikonal equation or static Hamilton-Jacobi equation. When the front speed $F$ does not depend on direction $v$, we refer to it just as an eikonal equation. In general, the front speed $F$ can be reconstructed from the particle speed function $f$ by taking the weighted dual norm 
\begin{equation}\label{eq:quasi:dual}
F(x,u) = \max_{v \neq 0} \frac{ u \cdot v }{ \norm{v}}f(x,v).
\end{equation}

Note that in the isotropic case, it is clear from $\eqref{eq:quasi:dual}$ that $F(x) = f(x)$. Namely, particle speeds are equal in all directions. In particular, this implies that the MAPs will always be normal to the wavefront. This is not the case for the quasipotential problem, where instead the direction of travel is given by 
\begin{equation}
\label{eq:dirphi}
\dot{\varphi} = \frac{b(\varphi)+\nabla U(\varphi)}{b(\varphi)}
\end{equation}
(see e.g. \cite{cameron12} for details). As we will see, this lack of orthogonality will pose some computational difficulties.

\begin{figure}
  \includegraphics{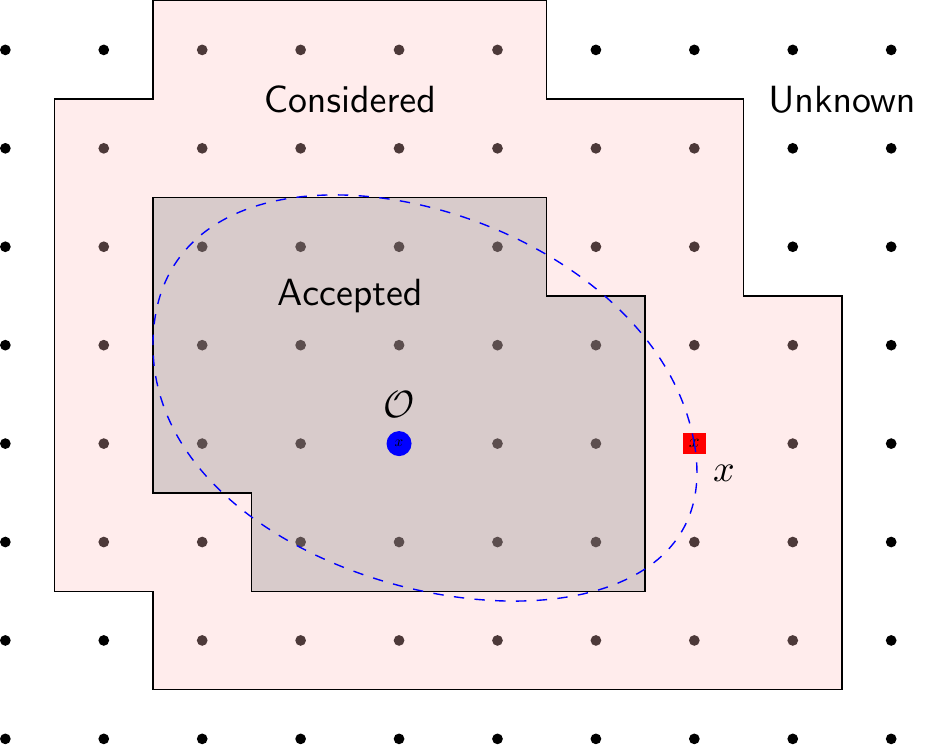}
\caption{Algorithm snapshot as mesh point $x$ is added to the {\sf Accepted} list. The dashed blue line represents a possible implied current state of the wavefront.}
\label{fig:quasi:groups}       
\end{figure}

We now introduce the main ideas of Sethian's fast marching algorithm for computing the viscosity solution $U(x)$ to equation $\eqref{eq:quasi:eikonal}$ in $2$-dimensions. Let $\mathcal{X}$ be a discretization of the domain $D$ into a uniform rectangular mesh with common horizontal and vertical spacing $h > 0$. We assume for simplicity that the attractor $\O$ is a member of the mesh $\mathcal{X}$. Dijkstra-like algorithms rely on the partition of $\mathcal{X}$ into the following three disjoint groups:
\begin{itemize}
	\item {\sf Unknown}: \textit{Mesh points for which no value of $U$ has been computed. The value of $U$ defaults to $+\infty$.}
	\item {\sf Considered}: \textit{Mesh points for which a tentative value of $U$ has been computed.}
	\item {\sf Accepted}: \textit{Mesh points for which a final, immutable value of $U$ has been computed.}
\end{itemize}
Initially, all mesh points will begin in the {\sf Unknown} category. As the algorithm proceeds, points will be gradually moved from {\sf Unknown} to {\sf Considered} as tentative values of $U$ are computed. For each mesh point, one of those tentative values will eventually be finalized, at which time that point will be moved into the {\sf Accepted} category. Figure \ref{fig:quasi:groups} displays a snapshot in time of what this setup may look like, while Algorithm \ref{alg:quasi:fastmarch} provides a template for fast marching Dijkstra-like eikonal solvers.

\begin{algorithm}
\raggedright
    \caption{Fast marching method template for solving eikonal equations}
    \label{alg:quasi:fastmarch}
	\textbf{Initialization} 
	\\
	Start with all mesh points in $\mathcal{X}$ in {\sf Unknown} and set them  to $U = +\infty$.
	\\
	Set $U(\O) = 0$ and add $\O$ to {\sf Considered}.
	\\
	\textbf{Main Body} 	
	\begin{algorithmic}[1]
	\While{$(\mathrm{{\sf Considered}}$ is non-empty$)$}
		\State Set $\ds x := \arg\min_{z \in \mathcal{X}} \Big\{ U(z) : \ z  \in \mathrm{{\sf Considered}} \Big\}$.
		\State Switch $x$ from $\mathrm{{\sf Considered}}$ to $\mathrm{{\sf Accepted}}$ and finalize its current value of $U(x)$.		
		\For {each \textit{neighbor} $y$ of $x$ such that $y \notin \mathrm{{\sf Accepted}}$} 
			\State Compute a value $U_{\mathrm{new}}(y)$ using the values of $U(x)$ and possibly $U(z)$ for other  
			$z$ in ${\sf Accepted}$.
			\State Set $U(y):= \min(U(y),U_\mathrm{new}(y))$. 
		\State Switch $y$ to $\mathrm{{\sf Considered}}$ if it was previously in {\sf Unknown}.
		\EndFor		
	\EndWhile
\end{algorithmic}
\end{algorithm}

Let us remark on some of the key steps of Algorithm \ref{alg:quasi:fastmarch}. First, in line 2, the mesh point $x$ with the smallest tentative value of $U$ among all {\sf Considered} points is selected. This point is effectively the next point to be ``hit'' by the expanding wavefront, and one can envision the level set $\{r \in \R^d:U(r) = U(x)\}$ as actually representing the current state of the front, as illustrated in Figure \ref{fig:quasi:groups}. Since the minimizer will need to be extracted from the {\sf Considered} list at each iteration, the {\sf Considered} list is often given a heap-sort structure in practice so that the argmin can be computed in $\O(\log N)$ operations.

Next, we note that the definition of \textit{neighbors} (line 4) and the process for computing estimates $U_{\mathrm{new}}(y)$ (line 5) are left unspecified in Algorithm \ref{alg:quasi:fastmarch}. These are the main areas in which Dijkstra-like eikonal solvers will differ, and we describe some of the possible choices for these procedures below. 

Finally, we note that after computing a new tentative value $U_{\mathrm{new}}(y)$, this value replaces the previous tentative value only if it is smaller. This particular choice of Accept/Reject rule tries to filter out those estimates of $U(y)$ that come from mesh points $x$ which are not necessarily near the MAP passing through $y$. 

\subsubsection{Computation of $\mathbf{U_{\mathrm{new}}}$}
We discuss techniques for computing the update value $U_\mathrm{new}(y)$ in line 5 of Algorithm \ref{alg:quasi:fastmarch}. Suppose, as in Algorithm \ref{alg:quasi:fastmarch} that $x$ has just been switched to {\sf Accepted} and we are now interested in using $U(x)$ to compute a tentative value of $U(y)$ for a nearby ``neighbor'' $y$. Appropriate definitions of ``neighbor'' are made precise in Section \ref{sec:neib}.

\begin{figure}
  \includegraphics[width=0.5\textwidth]{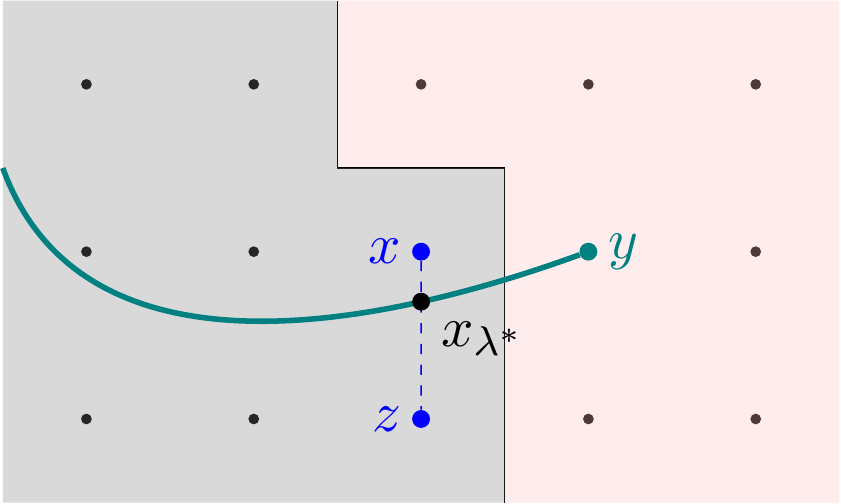}
\caption{The MAP of {\sf Considered} or {\sf Unknown} point $y$ passes between nearby {\sf Accepted} points $x$ and $z$.}
\label{fig:quasi:trajectory}       
\end{figure}

To compute $U(y)$, we seek to identify the MAP that passes through $y$ (Figure \ref{fig:quasi:trajectory}). In general, this MAP will not pass directly through $x$; however, it may pass between $x$ and another neighboring mesh point $z$ that also lies in the {\sf Accepted} group. In this case, the MAP will pass through the point $x_{\lambda^*} = (1-\lambda^*) x + \lambda^* z$ for some $\lambda^* \in [0,1]$. By noting that $U(y)$ is the minimum cost (geometric action in the quasipotential case) over all paths terminating at $y$, we can write
\begin{align}\label{eq:quasi:update_exact}
	U(y)  \nonumber &=  \min_{\lambda \in [0,1]} \Big[  U(x_\lambda) 	
	\\ & + \inf_{\substack{L > 0\\ \varphi \in C([0,L];\R^2)}} \left\{\int_0^L \frac{\norm{\dot{\varphi}(r)}dr}{f(\varphi(r),\dot{\varphi}(r))} : \varphi(0) = x_{\lambda}, \varphi(L) = y \right\} \Big].
\end{align}
The key decisions to be made concern how best to approximate the right hand side of equation $\eqref{eq:quasi:update_exact}$. This requires determining (a) how to interpolate $U(x_\lambda)$, since $x_\lambda$ in general lies between mesh points, (b) what path-space to take the inner minimization over, and (c) what quadrature rule to use to approximate the action integral.  

The simplest answers to these three questions are to (a) linearly interpolate $U(x_\lambda)$ between the finalized values of $U(x)$ and $U(z)$, (b) use a linear path connecting $x_\lambda$ and $y$, so that no inner minimization need occur at all, and (c) use right endpoint quadrature (at $y$) of the action integral. This route is taken, for instance, by the Ordered Upwind Method (OUM) of \cite{sethian03}. With these approximations, $\eqref{eq:quasi:update_exact}$ becomes
\begin{equation}\label{eq:quasi:update_linear}
	U_\mathrm{new}(y) = \min_{\lambda \in [0,1]} \Big[ (1-\lambda)U(x) + \lambda U(z) + \frac{\norm{y-x_\lambda}}{f(y,\frac{y-x_\lambda}{\norm{y-x_\lambda}})} \Big].
\end{equation}
This univariate minimization can then easily be carried out numerically. For the isotropic case where $f(y,v) = f(y)$, the minimizer $\lambda^*$ can be found analytically by solving a quadratic equation. 

The ordered line integral methods (OLIMs) \cite{cameron18,dahiya18} for solving for the quasipotential also use linear interpolation
for $U(x_{\lambda})$ and linear paths, but they employ higher order quadrature rules instead of the right-hand rule. The most efficient choice turned out to be the midpoint rule, which reduced the error constants by two-to-three orders of magnitude compared to right-hand rule quadrature. 

Update formula $\eqref{eq:quasi:update_linear}$ in general leads to an $\O(h)$ convergence rate. It is perfectly valid for the quasipotential problem, however in EJM we will opt for higher order approximations of $\eqref{eq:quasi:update_exact}$ in pursuit of a $\O(h^2)$ convergence rate.

 \subsubsection{Neighborhoods} 
 \label{sec:neib}
 The remaining undiscussed components of Algorithm \ref{alg:quasi:fastmarch} are $(i)$ how to choose the neighbors $y$ to update from $x$ (line 4) and 
 $(ii)$ how to choose points $z$ to pair with a given $x$ and $y$ for the computation of $\eqref{eq:quasi:update_linear}$ (line 5). 
 
The goal of \textit{neighborhood} design is to assure that every mesh point $y$ experiences at least one of the $\eqref{eq:quasi:update_exact}$ updates with an $x$ and $z$ that straddle its MAP as in Figure \ref{fig:quasi:trajectory}. Since we do not a priori know where the MAPs will be coming from, many updates will necessarily have to be performed where $x$ and $z$ do \textit{not} surround the MAP at all. These minimizations should result in boundary ($\lambda=0$ and $\lambda = 1$) minimizers of equation $\eqref{eq:quasi:update_exact}$. Moreover, these solutions will typically result in proposed values of $U_{\mathrm{new}}(y)$ that are \textit{larger} than the true value. Therefore, they should eventually be discarded by the Accept/Reject rule of line 6 if a proper interior solution is found. However, for this to happen, the neighborhoods must be large enough to include a valid $\triangle xzy$ triangle where the MAP passes between $x$ and $z$. We refer to such a triangle $\triangle xzy$ of mesh points, where $U(x) < U(y)$, $U(z) < U(y)$, and $y$'s MAP passes between $x$ and $z$, as a \textit{causal} triangle.

Due to this filtering system, it is in general harmless to increase the size of the neighborhoods, since any extraneous updates will automatically be filtered out by the Accept/Reject rule in the end, {\color{black} though the rule might be sophisticated -- see the end of Section \ref{sec:2ptu}}.
 The most obvious drawback of using very large neighborhoods is simply the additional computation time present in performing the additional minimizations $\eqref{eq:quasi:update_exact}$. The main goal of design is then to create neighborhoods that are just large enough to ensure that at least one causal triangle is checked for each $y$. 
If  a ball-shaped neighborhood is used for a mesh point $y$, the minimal necessary radius of the ball is determined by the anisotropy ratio of the speed function  \cite{sethian03}:
\begin{equation*}
r(y) := h\frac{\sup_{v} f(y,v)}{\inf_{v} f(y,v)}.
\end{equation*}
Note that for the optimal control problem for the quasipotential \eqref{eq_quasi_geometric}, this ratio is infinite.
In the next section, we discuss the construction of specialized star-shaped neighborhoods proposed by Mirebeau \cite{mirebeau14}.


\subsection{Anisotropic stencil refinement}\label{sec:quasi:mesh:mirebeau}
In this section we describe the Anisotropic Stencil Refinement (ASR) algorithm of Mirebeau \cite{mirebeau14} for solving anisotropic eikonal equations. 
The ASR algorithm follows the prescription of Algorithm \ref{alg:quasi:fastmarch} and uses the approximation scheme $\eqref{eq:quasi:update_linear}$, 
but with pre-computed stencils serving as the ``neighborhoods''. 
These stencils are significantly smaller than the neighborhoods used in \cite{sethian03}, so that the ASR algorithm is significantly faster than OUM and only slightly less accurate. 

Suppose we wish to solve the minimum cost problem associate with an action of the form $\eqref{eq:quasi:action_generic}$ where $f$ has a finite anisotropy ratio. 
For each mesh point $y$, we first construct the \textit{stencil} $\mathcal{N}(y)$, representing the candidates for $x$ and $z$ in a potentially causal triangle $\triangle xzy$. The neighborhood referred to in line 4 of Algorithm \ref{alg:quasi:fastmarch} will then be the \textit{reversed stencil}
\begin{equation*}
\mathcal{N}^{-1}(x):= \{ y \in \mathcal{X}: x \in \mathcal{N}(y)\},
\end{equation*}
which represents the points $y$ of which $x$ might be a member of a causal triangle.

Next, we define the function 
\begin{equation*}
\mathcal{F}(x,v):= s(x,v) \norm{v},
\end{equation*}
representing the integrand in the action functional $\eqref{eq:quasi:action_generic}$. We assume that for each $x$, the mapping $\mathcal{F}(x,\cdot):\R^2 \to [0,+\infty)$ is an \textit{asymmetric} norm, that is, it is subadditive and positive definite, but only satisfies positive homogeneity: $\F(x,\lambda v) = \lambda \F(x,v)$ for all $\lambda \geq 0$. 

\begin{definition}
Let $x$ be fixed. The vectors $u$ and $v$ are said to form an acute angle with respect to an asymmetric norm { $\mathcal{F}(x,\cdot):\R^2 \to [0,+\infty)$ (or $\F$-acute angle, for short) provided that
\begin{equation}\label{eq:quasi:mirebeau:acute}
	\mathcal{F}(x, u + \delta v) \geq \mathcal{F}(x,u), \hs \text{and} \hs \mathcal{F}(x,v + \delta u) \geq \mathcal{F}(x,v),
\end{equation} 
for all $\delta \geq 0$.}
\end{definition}
We note that if the asymmetric norm {$\F$} is differentiable, it is easy to see that acuteness conditions $\eqref{eq:quasi:mirebeau:acute}$ are equivalent to
\begin{equation*} {
	u \cdot \nabla_v \mathcal{F}(x,v) \geq 0, \hs \text{and} \hs v \cdot \nabla_v \mathcal{F}(x,u) \geq 0.}
\end{equation*}

The objective will then be to construct a stencil $\mathcal{N}(y)$ for any $y$ that is a collection of directions such that neighboring line segments terminating at $y$ form $\mathcal{F}(y,\cdot)$-acute angles. That is, for each mesh point $y$, we seek to construct a finite (and rotationally ordered) collection of mesh points $\mathcal{N}(y) = \{y_k\}_{k = 1}^{n_y}$, such that $y-y_k$ and $y-y_{k+1}$ form a $\mathcal{F}(y,\cdot)$-acute angle for each $k = 1,...,n_y$. 
Note that when $s(x,v) = s(x)$, the norm $\mathcal{F}(x,\cdot)$ is a multiple of the standard Euclidean norm so that $\mathcal{F}(x,\cdot)$-acute angles become ordinary Euclidean acute angles. Thus, the standard $4$-point diamond neighborhood is an admissible stencil, since neighboring directions are separated by right angles.  

The guarantee that such a stencil creates causal triangles is provided by the following causality property.
\begin{proposition}[Proposition 1.3 of \cite{mirebeau14}: Causality Property]
\label{prop:quasi:causality}
	Let $\mathcal{F}(x,\cdot)$ be an asymmetric norm on $\R^2$. Let $u,v \in \R^2$ be linearly independent and let $d_u,d_v \in \R$. Assume that $u$ and $v$ form an $\F$-acute angle. Define
	\begin{equation}\label{eq:quasi:causality}{
		d_w := \min_{t \in [0,1]}t d_u + (1-t) d_v +\mathcal{F}(x,tu+(1-t)v),}
	\end{equation}
	and assume that this minimum is not attained for $t \in \{0,1\}$. Then $d_u < d_w$ and $d_v < d_w$.

\end{proposition}

We do not provide the proof here, but rather defer it to \cite{mirebeau14}. It can be done rather succinctly using Lagrange multipliers. We do, however, provide a geometric proof in Section \ref{sec:quasi:hasr} of a similar statement in the quasipotential case that sheds more light on why this acuteness condition is relevant. 

Let us remark on Proposition \ref{prop:quasi:causality}. The right-hand side of $\eqref{eq:quasi:causality}$ is precisely the right-hand side of the update formula 
$\eqref{eq:quasi:update_linear}$ for $U_{\mathrm{new}}(y)$, if $u = x$, $v = z$, $d_u = U(x)$ and $d_v = U(z)$. 
Thus Proposition \ref{prop:quasi:causality} can be interpreted as saying that if the $xy$ and $zy$ legs of the $\triangle xzy$ triangle form a 
$\mathcal{F}(y,\cdot)$-acute angle and an interior minimum is found, then necessarily $U(x) < U_{\mathrm{new}}(y)$ and $U(z) < U_{\mathrm{new}}(y)$. 
This means, assuming $U(x) > U(z)$ without loss of generality, that when $x$ is switched to {\sf Accepted}, 
$z$ would already be an {\sf Accepted} point so that the triangle $\triangle xzy$ would yield a valid interior update, as desired.

The final step is then to define a method for creating causal stencils for a given problem. In \cite{mirebeau14}, 
Mirebeau defines stencils via Algorithm \ref{alg:quasi:stencil}, in which the base stencil is a $4$-point diamond that will be refined as necessary, 
until all neighboring angle are $\mathcal{F}(y,\cdot)$-acute. 

\begin{algorithm}[ht]
	\begin{algorithmic}[l]
	\For{$y \in \mathcal{X}$}
		\State Set $L:=[(1,0)]$ and $M:= [(1,0),(0,-1),(-1,0),(0,1)]$.
		\While{$M$ is non-empty}
			\State Set $u$ and $v$ to the last elements of $L$ and $M$, respectively.
			\If{$u$ and $v$ are $\mathcal{F}(y,\cdot)$-acute}
			\State Remove $v$ from $M$ and append $v$ to $L$.
			\Else
			\State Append $u+v$ to $M$.
			\EndIf
		\EndWhile
		\State Set $\mathcal{N}(y) := y- L$.
	\EndFor
	\State \textbf{Output} $\{\mathcal{N}(y)\}_{y \in \mathcal{X}}$.
\end{algorithmic}
\caption{FM-ASR Stencil Design \cite{mirebeau14}}
\label{alg:quasi:stencil}
\end{algorithm}

%
%

\section{Our algorithm: Efficient Jet Marcher}\label{sec:quasi:hasr}
Like the previously discussed solvers, our efficient jet marcher (EJM)  follows the general template of Algorithm \ref{alg:quasi:fastmarch}. Unlike the previously discussed solvers, EJM will instead use a higher order MAP approximation scheme as compared to $\eqref{eq:quasi:update_linear}$, as well as a modified version of Mirebeau's pre-computed stencils. The full structure of the solver is shown in Algorithm \ref{alg:quasi:EJM}. 

We remark briefly on key observations about Algorithm \ref{alg:quasi:EJM}; detailed descriptions are provided in the subsequent sections. The first observation is that the gradient $\nabla U$, is now treated as part of the solution and is included in all of the update computations. This will be necessary to perform the higher order MAP interpolations in the update step of line 7 (Section \ref{sec:quasi:update}).

As in the OLIM methods \cite{cameron18,dahiya18}, we apply a slightly more involved initialization process. This is typically only needed when the MAPs exhibit significant curvature near the origin. The most common technique is to simply initialize $U$ and $\nabla U$ with the exact solution corresponding to a linearized version of $b$ around the attractor $\O$.

Concerning the construction of anisotropic stencils, we follow the ideas of Mirebeau \cite{mirebeau14}, discussed in Section \ref{sec:quasi:mesh:mirebeau}. However, we also use the fact that stencils only depend on the angle of $b(x)$ and not its magnitude, to save the stencils only for a binned collection of possible angles of $b$. This drastically reduces the memory requirement and the runtime of pre-processing phase. Moreover, this allows us to skip the stencils entirely and compute the reversed stencils directly, since these are ultimately the objects used in the body of Algorithm \ref{alg:quasi:EJM}. We also use a slightly different stencil construction algorithm. All of which is described in Section \ref{sec:quasi:stencil}.

The if-statement and subsequent fail-safe method, mentioned in lines 3 and 4, are used to prevent values of $U(x)$ computed from one-point updates (boundary solutions to the update minimization problem $\eqref{eq:quasi:update_exact}$) from becoming finalized values. As we will see, due to the assumption of linear MAPs in their construction,
the stencils are not perfect and will occasionally fail to find causal triangles when the mesh is not sufficiently refined. If these failures are not caught and corrected, the higher order accuracy may not be achieved. This fail-safe is called only very rarely. When called, it searches a much larger area than the stencil until it finds a successful triangle update (interior solution of $\eqref{eq:quasi:update_exact}$).

Finally, the most important difference lies in the structure of the update procedure (line 7) and prescription for computing the right-hand side of $\eqref{eq:quasi:update_exact}$. This procedures is responsible for the $\O(h^2)$ error convergence rate. It is discussed in detail in Section \ref{sec:quasi:update}.

\begin{algorithm}
\raggedright
    \caption{Efficient Jet Marching Algorithm (EJM)}
    \label{alg:quasi:EJM}
	\textbf{Initialization and Pre-processing}
	\\ Start with all mesh points in $\mathcal{X}$ in {\sf Unknown} and set them to $U = +\infty$ and $\nabla U = (\infty,\infty)$.
	\\ 
	Initialize the $U$ and $\nabla U$ values of an 8-point rectangular neighborhood of the attractor $\O$. 
	\\
	Switch each of these points into $\mathrm{{\sf Considered}}$. 
	\\ 
	For each $\theta_k = \frac{2\pi k}{N_{\mathrm{bins}}}$, $k = 0,...,N_{\mathrm{bins}}-1$, build the reversed stencils $\mathcal{N}^{-1}{(\theta_k)}$ (Section \ref{sec:quasi:stencil}).
	\\		
	For each $x$  { in the 8-point neighborhood of the attractor $\O$},
	assign $\mathcal{N}^{-1}(x) := \mathcal{N}^{-1}(\theta_k)$ where $k$ is such that
	\begin{equation*}
		-\angle b \in \Big(\theta_k - \frac{\pi}{N_{\mathrm{bins}}}, \theta_k + \frac{\pi}{N_{\mathrm{bins}}} \Big],
	\end{equation*}
	 { where $\angle b$ denotes the angle between $b$ and the positive direction of the $x$-axis, i.e., the polar angle of $b$.}\\
	\textbf{Main Body} 	
	\begin{algorithmic}[1]
	\While{$(\mathrm{{\sf Considered}}$ is non-empty$)$}
		\State Set $x :=  \ds \arg\min_{z \in \mathcal{X}} \{ U(z) : \ z  \in \mathrm{{\sf Considered}} \}$.
		\If {$x$'s last update is a one-point update} 
			 \State Run the fail-safe on $x $ (Section \ref{sec_update_fail}). 
		
		\EndIf 
		\State Switch $x$ from $\mathrm{{\sf Considered}}$ to $\mathrm{{\sf Accepted}}$.		
		\For {each $y \in \mathcal{N}^{-1}(x)$ such that $y \notin \mathrm{{\sf Accepted}}$} 
			\State Update $U(y)$  and $\nabla U(y)$ from $x$ and possibly other $z \in \mathrm{{\sf Accepted}}$ (Section \ref{sec:quasi:update}).
		\State Switch $y$ to $\mathrm{{\sf Considered}}$ if it was previously in {\sf Unknown}.
		\EndFor		
	\EndWhile
\end{algorithmic}
\end{algorithm}

\subsection{Anisotropic stencils}\label{sec:quasi:stencil}
We return to the anisotropic stencil ideas discussed in Section \ref{sec:quasi:mesh:mirebeau}. For the quasipotential problem, we have
\begin{equation*}
\mathcal{F}(y,v) = \norm{b(y)} \norm{v}- b(y) \cdot v.
\end{equation*}
This is subadditive and positive homogeneous, but not positive definite, since 
\begin{equation*}
\mathcal{F}(y,\lambda b(y)) = 0,
\end{equation*} for any $y$ and any $\lambda > 0$. Nonetheless, the machinery of the ASR method can still be applied, for instance, by considering instead a modified function 
\begin{equation*}
\mathcal{F}^\alpha(y,v) = \norm{b(y)} \norm{v} - \alpha b(y) \cdot v,
\end{equation*}
where $\alpha = 1 -\delta$ for some $\delta >0$ very small. Such an $\mathcal{F}^\alpha(y,\cdot)$ is an asymmetric norm for any $y$ such that $b(y) \neq 0$.

For each $y \in \mathcal{X}$, we seek to construct a stencil $\mathcal{N}(y)$ such that neighboring points $y_k,y_{k+1} \in \mathcal{N}(y)$ form line segments $u=\frac{y-y_k}{\norm{y-y_k}}$ and $v =\frac{y-y_{k+1}}{\norm{y-y_{k+1}}}$ that satisfy the acuteness conditions $\eqref{eq:quasi:mirebeau:acute}$, which, in the quasipotential case, reduces to the conditions
\begin{equation*}
u \cdot v \geq u \cdot (-b(y)), \hs \text{and} \hs u \cdot v \geq v \cdot (-b(y)).
\end{equation*}
It is not intuitively clear where these acuteness conditions come from. To motivate them, we provide a simple geometric proof of the following proposition, which is essentially a corollary of Proposition \ref{prop:quasi:causality} for the specific case of the quasipotential problem.

\begin{proposition}
Suppose the drift field $b$ is smooth, fix $y \in \R^2$ and let $\varphi$ denote the MAP passing through $y$. Suppose that $D$ is a rotationally ordered collection of unit directions in $\R^2$ such that any neighbors $\hat{u},\hat{v} \in D$ satisfy 
\begin{equation}\label{eq:quasi:acute_mine}
	\hat{u} \cdot \hat{v} \geq \max \Big(- \hat{u} \cdot \hat{b},- \hat{v} \cdot \hat{b},0\Big),
\end{equation}
where $\hat{b} = b(y)/\|b(y)\|$. Then there exist neighbors $\hat{u},\hat{v} \in D$ and $h > 0$ small enough that $\varphi$ passes between $A := y + h\hat{u}$ and $B:= y + h\hat{v}$ while { $U(A) < U(y)$ and $U(B) < U(y)$}. 
\end{proposition}
\begin{proof}
Suppose first that $U$ is continuously differentiable at $y$. Let $\hat{\tau}$ be the tangent vector of the level set { $\{z:U(z) = U(y)\}$ at $y$}. Picking scale $h >0$ such that $U$ and $\varphi$ are locally flat, the setup looks like Figure \ref{fig_proofcausality}. By taking the dot product of both sides of equation $\nabla U(y) = \|b(y)\|\dot{\varphi}-b(y)$  (see \eqref{eq:gradUdir}) with $\hat{\tau}$, we immediately have that
	\begin{equation*}
		\Big(\frac{\dot{\varphi}}{\norm{\dot{\varphi}}}\Big) \cdot \hat{\tau} = \hat{b} \cdot \hat{\tau} =:\cos(\theta).
	\end{equation*}	
	
\begin{figure}
\includegraphics[width=0.5\textwidth]{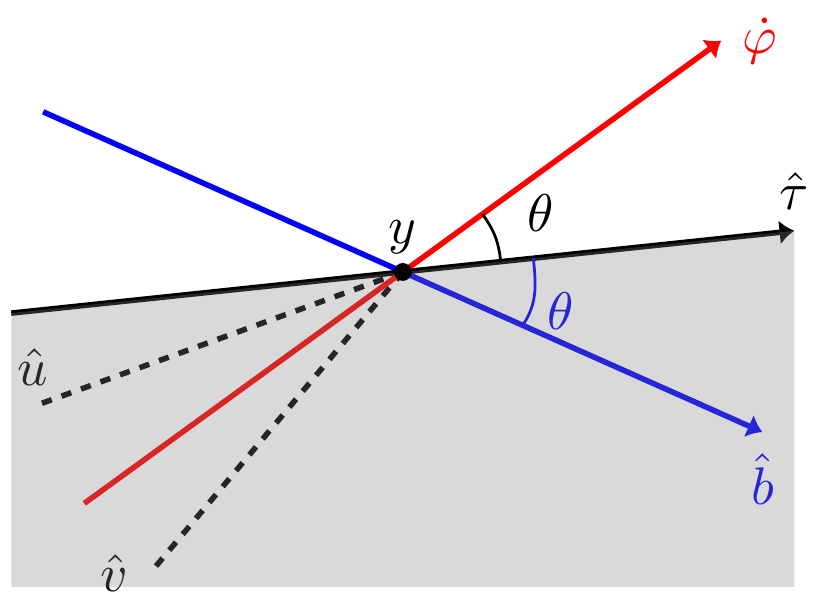}
\caption{If $\eqref{eq:quasi:acute_mine}$ is satisfied, there must exist neighbors $\hat{u}$, $\hat{v} \in D$ that sandwich the incoming MAP $\varphi$ and both lie in the shaded region $\{x:U(x) < U(z)\}$. }
\label{fig_proofcausality}
\end{figure}

	The result then follows by noting that there \textit{must} exist a direction $\hat{u} \in D$ lying between $-\dot{\varphi}$ and $-\hat{\tau}$, as in Figure \ref{fig_proofcausality}. We can see this by realizing that there \textit{cannot} exist neighbors $\hat{u},\hat{v} \in D$ that straddle the wedge created by $-\hat{\tau}$ and $-\dot{\varphi}$. If there did exist such neighbors, $\hat{u}$ between $-\hat{b}$ and $-\hat{\tau}$ and $\hat{v}$ between $\hat{b}$ and $-\dot{\varphi}$, then we would immediately contradict the acuteness requirement, since
	\begin{equation*}
		 -\hat{b} \cdot \hat{u} > -\hat{b} \cdot (-\hat{\tau}) = \cos \theta = -\Big(\frac{\dot{\varphi}}{\norm{\dot{\varphi}}}\Big)  \cdot (- \hat{\tau}) > \hat{u} \cdot \hat{v}.
	\end{equation*}
	Hence nearest neighbors cannot straddle the wedge, and there must be neighbors $\hat{u},\hat{v} \in D$ surrounding $\dot{\varphi}$ that both lie in the half-plane of directions opposite $\nabla U$.
\end{proof}

We now seek to construct stencils such that neighboring directions satisfy acuteness condition $\eqref{eq:quasi:acute_mine}$. 
To do this, we make the following two approximations to reduce computation time and memory requirements.
\begin{enumerate}[(1)]
	\item We take the reversed stencils to be
\begin{equation}\label{eq:quasi:stencil_approx}
	\mathcal{N}^{-1}(y) = -\mathcal{N}(y),
\end{equation}
so that the \textit{reversed} stencils are readily obtained by choosing a set of directions satisfying 
\begin{equation*}
	\hat{u} \cdot \hat{v} \geq \max(\hat{u} \cdot \hat{b}, \hat{v} \cdot \hat{b}, 0)
\end{equation*} 
and the original stencils $\mathcal{N}(y)$ need not be computed at all.
\item We save reversed stencils only for $N_{\mathrm{bin}}$ different uniformly spaced values of the angle of $b(y)$. 
\end{enumerate}

Approximation (1) is taken to speed up the computation time of the pre-processing phase of Algorithm \ref{alg:quasi:EJM}. The process of inverting the stencil can be time consuming for fine meshes, since its complexity is $\mathcal{O}(K(\mathcal{N})h^{-2})$ where $K(\mathcal{N})$ is the average stencil cardinality. Moreover, the approximation $\eqref{eq:quasi:stencil_approx}$ is very close to exact on a scale $h > 0$ such that the MAPs are approximately flat.

Approximation (2) is taken to reduce memory requirements of the algorithm. Here, we take use of the fact that the metric $\mathcal{F}(y,\cdot)$ only depends on the angle of $b(y)$. Therefore, the reversed stencils can be saved for a common set of $\theta = \angle b$ values rather than one for every mesh point $y \in \mathcal{X}$. Angles are binned as mentioned in the pre-processing phase of Algorithm \ref{alg:quasi:EJM}. We do note however, that this shortcut is only possible because we are using a mesh that is translation invariant, so that we need only to store the shifts in the reversed stencil. This would not be possible if instead we were using a non-uniform mesh. 

To construct the reversed stencils, one option is to follow the prescription of Algorithm \ref{alg:quasi:stencil}. However, since $\mathcal{F}(y,\cdot)$ is not positive definite, one should modify the process or else it will never terminate, as it will attempt to refine the reversed stencil in the direction of $-b(y)$ indefinitely. One solution, is to instead use the modified version of $\mathcal{F}^\alpha$, mentioned above, for some $\alpha$ very slightly below $1$. The value of $\alpha$ will control how refined the reversed stencil is in the $-b(y)$ direction: our typical values is $0.9999$. An alternative option is simply to cut off the production algorithm once the leg in the direction of $-b(y)$ reaches some maximum threshold value. 

\begin{figure}
\includegraphics{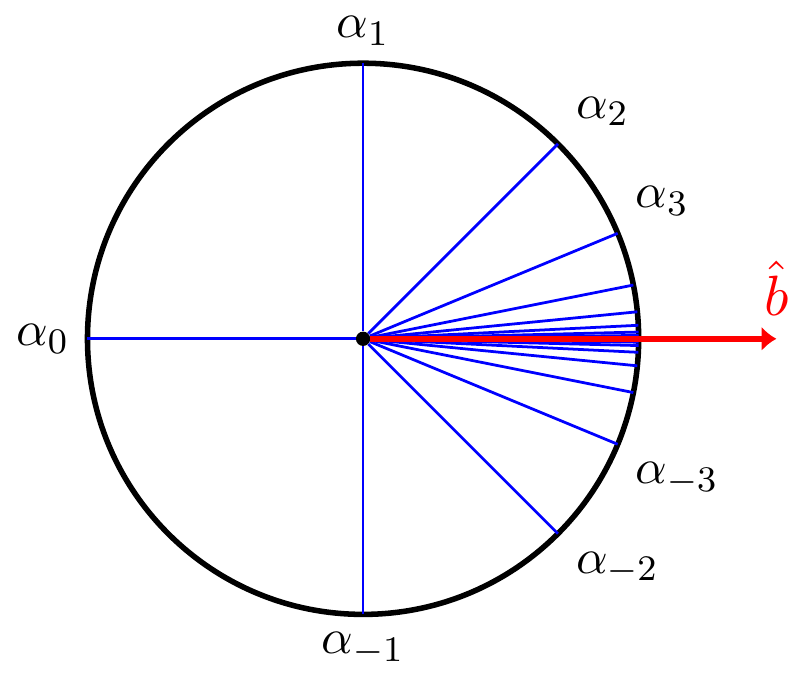}
\caption{The infinite collection of angles $ \{\alpha_k\}_{k=0}^\infty$ with $\alpha_k$ given by $\eqref{eq:quasi:perfect}$ satisfy $\eqref{eq:quasi:action_generic}$ with equality except right at the direction $\hat{b}$.}
\label{fig:quasi:perfect}
\end{figure}

We opted to construct the reversed stencils via another route. 
Instead, we noted that for a given value of $\angle b$, a minimal set of directions satisfying relations $\eqref{eq:quasi:acute_mine}$ 
can be easily written down. In fact, consider the infinite collection of angles $\{\alpha_k\}_{k=-\infty}^\infty$ where
\begin{equation}\label{eq:quasi:perfect}
	\alpha_k = \angle \hat{b} + \mathrm{sgn}(k) \frac{\pi}{2^{|k|}},
	\end{equation} 
which is displayed in Figure \ref{fig:quasi:perfect}. Here, it is clear that $\alpha_{k+1} - \alpha_k = \alpha_k - \angle b$ so that relation $\eqref{eq:quasi:acute_mine}$ holds with equality. Such a simple characterization of an ``ideal'' stencil is certainly not possible for the general case of an anisotropic eikonal equation. Thus, rather than following Algorithm \ref{alg:quasi:stencil} which is designed to work for the general case, we seek to construct stencils by directly discretizing rotated versions of the collection shown in Figure \ref{fig:quasi:perfect}.

\begin{figure}
\includegraphics[width=0.75\textwidth]{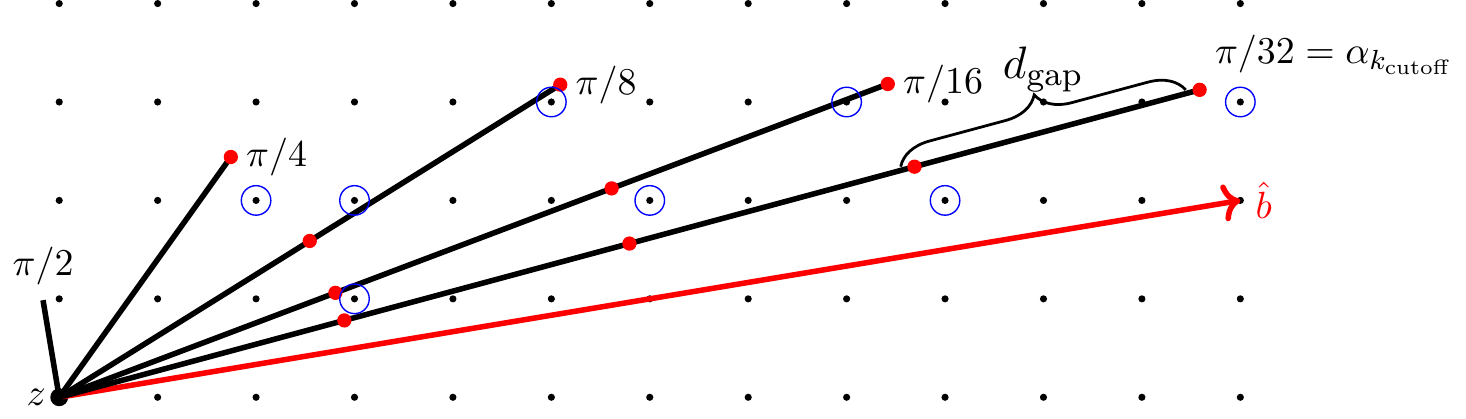}
\caption{The stencil consists of all of the blue circled mesh points, along with the analog procedure done on the other side of $\hat{b}$ and overlayed with a standard $8$-point neighborhood of $z$.}
\label{fig_stencil_mine}
\end{figure}

The construction of our discretization of Figure \ref{fig:quasi:perfect} is given in Algorithm \ref{alg:mystencil}, with the supporting Figure \ref{fig_stencil_mine}. The reversed stencil consists of a discretization of $|k|-1$ equi-spaced points along each of the $\alpha_k$ rays in Figure \ref{fig:quasi:perfect}, up until some maximum $k$ value. Since these are all non-integer points, their nearest neighbor in $\mathbb{Z}^2$ is added instead. In addition to these points, the standard 8-point neighborhood is included. 

\begin{figure}
(a)\includegraphics[width=0.4\textwidth]{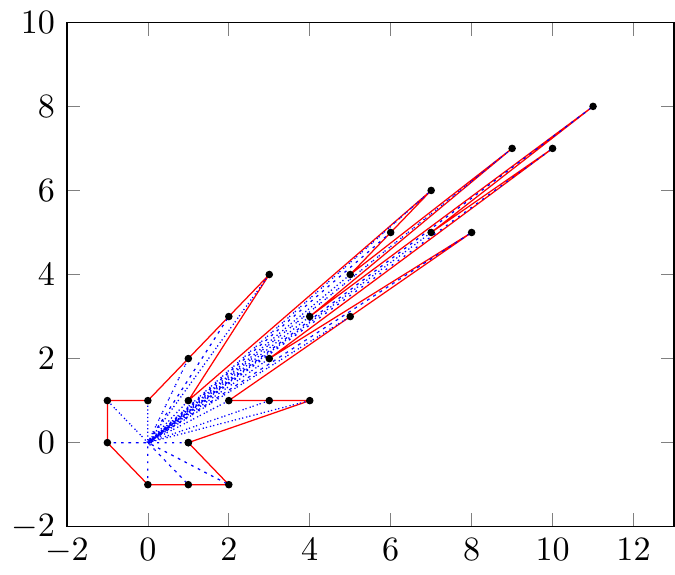}
(b)\includegraphics[width=0.4\textwidth]{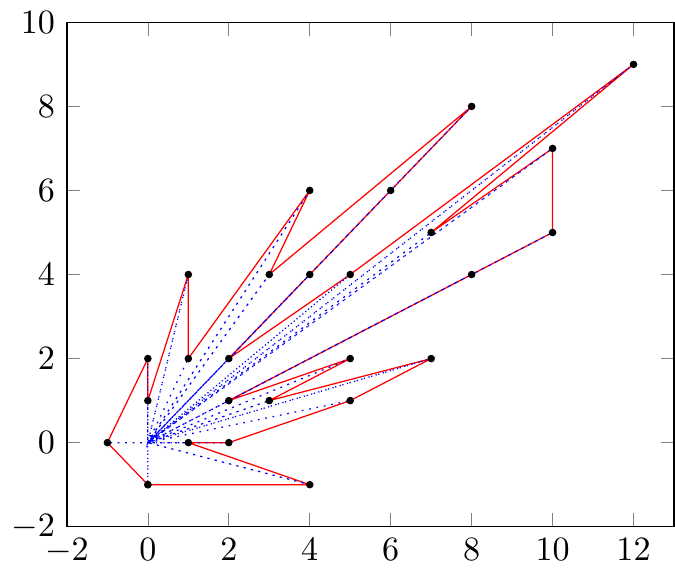}
\caption{Sample Reversed Stencils created using Algorithms \ref{alg:quasi:stencil} (a) and \ref{alg:mystencil} (b), respectively. Both of these stencils have the same cardinality.}
\label{fig:quasi:stencil_example}
\end{figure}

Algorithm \ref{alg:mystencil} is a rather simple method for  discretizing Figure \ref{fig:quasi:perfect} in a way that is both non-hollow and elongated in the direction of $\hat{b}$. Example stencils created via the two methods are displayed in Figure \ref{fig:quasi:stencil_example}. The stencil created via Algorithm \ref{alg:mystencil} is slightly less concentrated in the direction of anisotropy and has denser interior: properties that have shown to perform slightly better in practice by causing fewer calls of the fail-safe procedure.

\begin{algorithm}
\raggedright
\caption{Construction of dense oblong reversed stencils}\label{alg:mystencil}
Let $\theta$ be the binned value of $-\angle b$.

Set $k_{\mathrm{cutoff}} = $ largest $k$ to include in discretization of $\alpha_k$ (Figure \ref{fig:quasi:perfect}). Value $k_{\mathrm{cutoff}} = 4$ is typical.

Set $d_{\mathrm{gap}} = $ spacing between sampled points on each ray (in multiples of $h$). Value $d_{\mathrm{gap}} = 3$ is typical.

\begin{algorithmic}[1]
\State Add $\{(1,0), (1,1),(0,1),(-1,1),(-1,0),(-1,-1), (0,-1),(1,-1)\}$ to $\mathcal{N}^{-1}(\theta)$.

\State Let $M$ be an empty list of points.
\For{$k = 2:k_{\mathrm{cutoff}}$}
	\For{$n = 1:|k|-1$}
	\State Add $n \cdot d_{\mathrm{gap}} \cdot (\cos(\theta + \frac{\pi}{2^k}),\sin(\theta+\frac{\pi}{2^k})$ to $M$
	\EndFor
\EndFor
\For{each point $(x,y) \in M$}
	\State Add the closest integer pair $({\sf round}(x),{\sf round}(y))$ to $\mathcal{N}^{-1}(\theta)$.
\EndFor

\end{algorithmic}

\end{algorithm}

\subsection{Update procedure}
\label{sec:quasi:update}
In this section, we discuss how the update procedure (line 7 of Algorithm \ref{alg:quasi:EJM}) is conducted. 
In EJM, this procedure is considerably more involved than in the OUM where a quadratic equation is solved, and in the OLIMs where a 1D minimization problem is solved. 
We provide a full description of the update process in Algorithm \ref{alg:quasi:update} at the end of this section.

As mentioned in Section \ref{sec:quasi:hasr}, we henceforth distinguish between the process of checking for interior and exterior solutions to $\eqref{eq:quasi:update_exact}$, which we refer to as one-point updates and triangle updates, respectively. 
That is, one-point updates compute values of $U(y)$ and $\nabla U(y)$ from a single mesh point $x$, while triangle updates, 
when successful, compute values of $U(y)$ and $\nabla U(y)$ from two mesh points $x$ and $z$.

\subsubsection{One-point updates} 
As in Algorithm \ref{alg:quasi:EJM}, we suppose that mesh point $x$ has just been switched to {\sf Accepted}, 
and the mesh point $y \in \mathcal{N}^{-1}(x)$ is to be updated from $x$. Let $\varphi$ be the MAP that passes through $y$. 
In view of equation $\eqref{eq:quasi:update_exact}$, it follows that
\begin{equation}\label{eq:quasi:update_onept_exact}
	U(y)  \geq   U(x) 
	+ \inf_{L > 0,\varphi \in C([0,L];\R^d)} \{\int_0^L \frac{\norm{\dot{\varphi}(r)}dr}{f(\varphi(r),\dot{\varphi}(r))} : \varphi(0) = x, \varphi(L) = y \},
\end{equation}
with equality if and only if $x$ lies on the MAP $\varphi$. 
The one-point update will be a numerical approximation of this right-hand side. 
The resulting proposed value of $U(y)$ will in general only be a good estimate of the true value if $x$ lies very near $\varphi$. 
It will be a significant overestimate otherwise.

In approximating the right-hand side of $\eqref{eq:quasi:update_onept_exact}$, we opt to conduct the minimization over the two-parameter family of paths { $\tilde{\varphi}_{\alpha,\beta}:[0,L]\to \R^2$,  $L:=\|x-y\|$,} consisting of cubic curves with fixed endpoints at $x$ and $y$. This family is parametrized by the entry and exit angles $\alpha$ and $\beta$ (Figure \ref{fig:updates}(a)), while each individual path $\tilde{\varphi}_{\alpha,\beta}$ is parametrized by its coordinate in the $y-x$ direction. Complete formulas can be found in Appendix A.

For a given path $\tilde{\varphi}_{\alpha,\beta}$, we approximate the integral in $\eqref{eq:quasi:update_onept_exact}$ by using  
Simpson's quadrature rule. We also adopt the notation $U_{\mathrm{new}}^x(y)$ and $\nabla U_{\mathrm{new}}^x(y)$ to denote the proposed values of $U(y)$ and $\nabla U(y)$, respectively, from the one-point update from mesh point $x$. Hence, we set
{
\begin{align}
 U_{\mathrm{new}}^x(y) := U(x) &+  \min_{\alpha,\beta \in [0,2\pi]} \frac{\norm{y-x}}{6} \Big[  \frac{\norm{\dot{\tilde{\varphi}}_{\alpha,\beta}(0)}}{f(x,\dot{\tilde{\varphi}}_{\alpha,\beta}(0))} \notag
	\\ & + 4 \frac{\norm{\dot{\tilde{\varphi}}_{\alpha,\beta}(\sfrac{L}{2})} }{f(\tilde{\varphi}_{\alpha,\beta}(\sfrac{L}{2}),\dot{\tilde{\varphi}}_{\alpha,\beta}(\sfrac{L}{2}))} 
	 + \frac{ \norm{\dot{\tilde{\varphi}}_{\alpha,\beta}(L)}}{f(y,\dot{\tilde{\varphi}}_{\alpha,\beta}(L))} \Big]. \label{eq:quasi:update_1pt_cubic}
\end{align}
}
We conduct this two-dimensional minimization over $\alpha$ and $\beta$ by using Newton's method. Details are provided in Appendix A.

In view of relation $\eqref{eq:gradUdir}$, the proposed one-point update value of $\nabla U(y)$ is
{
\begin{equation*}
 \nabla U_{\mathrm{new}}^{x}(y) = \norm{b(y)} \dot{\tilde{\varphi}}_{\alpha^*,\beta^*}(L) - b(y),
\end{equation*}
where $\alpha^*$ and $\beta^*$ are the minimizing angles of $\eqref{eq:quasi:update_1pt_cubic}$. 
}
We remark that if the numerical solver fails to find a minimizer of $\eqref{eq:quasi:update_1pt_cubic}$, then the values $\alpha= 0$ and $\beta = 0$ are used to compute $U_{\mathrm{new}}^x(y)$ and $\nabla U_{\mathrm{new}}^x(y)$. The proposed update value then represents an approximation of using a linear MAP and Simpson's rule quadrature. From our empirical observations, such a failure is extremely rare unless there is a significant amount of curvature of the MAP over the one-point update in question. In scenarios where the MAP is very flat, minimizing values of $\alpha$ and $\beta$ are very close to $0$ and the Newton solver has little difficulty finding them.


\subsubsection{Triangle Updates} 
\label{sec:2ptu}
Conversely, suppose now that points $x$ and $z$ are {\sf Accepted} and that we seek to update the values of $U(y)$ and $\nabla U(y)$ using $x$ and $z$. The triangle update will search for interior (in $\lambda$) minimizers of $\eqref{eq:quasi:update_exact}$ corresponding to a situation described by Figure \ref{fig:quasi:trajectory}.

\begin{figure}
(a)\includegraphics[width=0.4\textwidth]{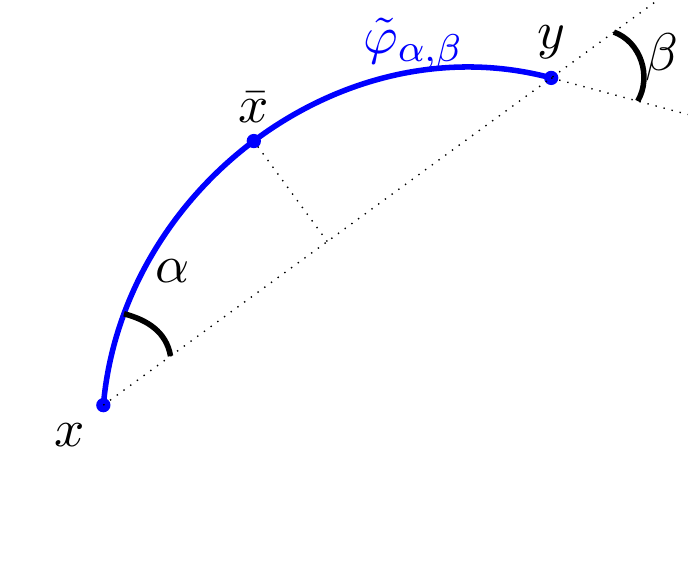}
(b)\includegraphics[width=0.4\textwidth]{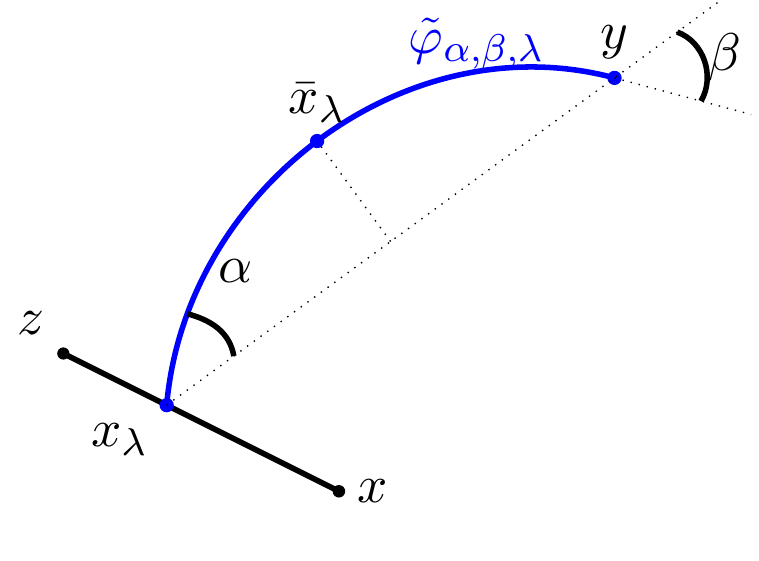}
\caption{Diagrams for the (a) one-point updates and (b) triangle updates. 
The MAP is approximated by minimizing the geometric action over the families $\tilde{\varphi}_{\alpha,\beta}$ and $\tilde{\varphi}_{\alpha,\beta,\lambda}$, respectively.}
\label{fig:updates}
\end{figure}

As in the case of the one-point update, we approximate the MAP with a cubic polynomial terminating at $y$. However, here the starting point $x_{\lambda} = (1-\lambda)x + \lambda z$ is allowed to vary on along the $xz$ line segment as part of the minimization. 

In approximating the right-hand side of $\eqref{eq:quasi:update_exact}$, we opt to conduct the minimization over the three-parameter family of paths { $\tilde{\varphi}_{\alpha,\beta,\lambda} :[0,L] \to \R^2$, $L:=\|y-x_{\lambda}\|$, } consisting of cubic curves with fixed endpoints at $x_\lambda$ and $y$ (Figure \ref{fig:updates}(b)). As before, this family is parametrized by the entry and exit angles $\alpha$ and $\beta$, in addition to the value of $\lambda \in [0,1]$. Each path $\tilde{\varphi}_{\alpha,\beta,\lambda}$ is parametrized by its coordinate in the $y-x_\lambda$ direction. In particular, the coordinate systems vary with the value of parameter $\lambda$.

Notice that in the one-point update, we did not use the value of $\nabla U(x)$ in order to perform the update. We use it here, however, to interpolate the value of $U(x_\lambda)$. In particular, we interpolate $U(x_\lambda)$ with the unique cubic Hermite polynomial $p:[0,1] \to \R$ satisfying boundary conditions
\begin{equation*}
\begin{cases}
p(0) = U(x), \hs p'(0) = \nabla U(x) \cdot (z-x),
\\ p(1) = U(z), \hs p'(1)=\nabla U(z) \cdot (z-x).
\end{cases}
\end{equation*}

Then, similar to the case of the one-point update, we take a Simpson's rule approximation of the integral in $\eqref{eq:quasi:update_exact}$. As before, we adopt the notation $U_{\mathrm{new}}^{x,z}(y)$ and $\nabla U_{\mathrm{new}}^{x,z}(y)$ to denote the proposed values of $U(y)$ and $\nabla U(y)$, respectively, from the triangle update of $y$ from $x$ and $z$. Thus, we set
{
\begin{align}\label{eq:quasi:update_triangle_cubic}
	\nonumber U_{\mathrm{new}}^{x,z}(y)& := \min_{\lambda \in [0,1], \alpha,\beta \in \R} \Big( p(\lambda) + \frac{\norm{y-x_\lambda}}{6} \Big[ \frac{\norm{\dot{\tilde{\varphi}}_{\alpha,\beta,\lambda}(0)}}{f(x_\lambda,\dot{\tilde{\varphi}}_{\alpha,\beta,\lambda}(0))} 
	\\ & + 4 \frac{\norm{\dot{\tilde{\varphi}}_{\alpha,\beta,\lambda}(\sfrac{L}{2})} }{f(\tilde{\varphi}_{\alpha,\beta,\lambda}(\sfrac{L}{2}),\dot{\tilde{\varphi}}_{\alpha,\beta,\lambda}(\sfrac{L}{2}))} + \frac{ \norm{\dot{\tilde{\varphi}}_{\alpha,\beta,\lambda}(L)}}{f(y,\dot{\tilde{\varphi}}_{\alpha,\beta,\lambda}(L))} \Big] \Big).
\end{align}
}
As before, we perform the three-dimensional minimization over $\alpha,\beta$ and $\lambda$ using Newton's method with exact derivatives, see Appendix A. 
Similarly, the proposed update value of the gradient is 
{
\begin{equation*}
 \nabla U_{\mathrm{new}}^{x,z}(y) = \norm{b(y)}\dot{\tilde{\varphi}}_{\alpha^*,\beta^*,
 \lambda^*}(L) - b(y),
\end{equation*}
}
where $\alpha^*$, $\beta^*$, $\lambda^*$ are the minimizing parameter values  of $\eqref{eq:quasi:update_triangle_cubic}$.

We are only interested in interior solutions since exterior solutions are handled by separate one-point updates. 
As such, we terminate the numerical solver once the working value of $\lambda$ leaves the interval $[0,1]$ and return \texttt{failure} of the triangle update. Since the triangle update should only succeed if the MAP passes through the $xz$ line segment, the vast majority of attempted triangle updates should result in failure.

We now state the full update procedure (Algorithm \ref{alg:quasi:update}). 
Justification for lines 2 and 7,  Accept/Reject Rule \ref{con:update}), and the choice of points $z$ for the triangle updates (line 4) are given in Section \ref{sec:quasi:obstacle}.

\begin{algorithm*}
\raggedright
    \caption{Update neighbors $y$ of $x$ (line 7 of Algorithm \ref{alg:quasi:EJM})}\label{alg:quasi:update}
    Let $x$ be newly {\sf Accepted} and $y \in \mathcal{N}^{-1}(x)$ be {\sf Unknown} or {\sf Considered}.
\begin{algorithmic}[1]
	\State Compute one-point update values $U_{\mathrm{new}}^x(y)$ and $\nabla U_{\mathrm{new}}^x(y)$ from $\eqref{eq:quasi:update_1pt_cubic}$.
	\If{$U_{\mathrm{new}}^x(y) < U(y)$}
		\State Set $U(y):= U_{\mathrm{new}}^x(y)$ and $\nabla U(y):=\nabla U_{\mathrm{new}}^x(y)$.
	\EndIf
	\For{$z$ in the 8-point nearest neighborhood of $x$ (Figure \ref{fig:quasi:updates})}
		\If{$z$ is {\sf Accepted}}
			\State Compute triangle update values $U_{\mathrm{new}}^{x,z}(y)$ and $\nabla U_{\mathrm{new}}^{x,z}(y)$ from $\eqref{eq:quasi:update_triangle_cubic}$.	
			\If{triangle update is successful and the Accept/Reject Rule \ref{con:update} is met}
				\State Set $U(y):= U_{\mathrm{new}}^{xz}(y)$ and $\nabla U(y):=\nabla U_{\mathrm{new}}^{xz}(y)$.
				\EndIf
		\EndIf
	\EndFor
    \end{algorithmic}
\end{algorithm*}

\begin{ARrule}
\label{con:update}
Let $U_{\mathrm{new}}^{x,z}(y)$ and $\nabla U_{\mathrm{new}}^{x,z}(y)$ be the values proposed from a 
successful triangle update with minimizing $\lambda$ value $\lambda^*$. The proposed values are {\sf Accepted} if and only if
\begin{enumerate}[(A)]
	\item $U_{\mathrm{new}}^{x,z}(y)$ is smaller than \textit{all} previous proposed one-point update values $U_{\mathrm{new}}^\cdot(y)$ of $y$.
	\item If the current tentative value of $U(y)$ came from a \textit{triangle update} $\triangle x_{\mathrm{old}} z_{\mathrm{old}}y$, then we must have
	\begin{equation*}
		\norm{x_{\lambda^*}-y} < \norm{x_{\lambda^*_{\mathrm{old}}} - y },
	\end{equation*}
	where $x_{\lambda^*_{\mathrm{old}}} = (1-\lambda^*_{\mathrm{old}})x_{\mathrm{old}} + \lambda^*_{\mathrm{old}}z_{\mathrm{old}}$ and $\lambda^*_{\mathrm{old}}$ is the minimizing $\lambda$ from the $x_{\mathrm{old}} z_{\mathrm{old}}y$ update.
\end{enumerate}
\end{ARrule}

\begin{figure}
\includegraphics[width=0.5\textwidth]{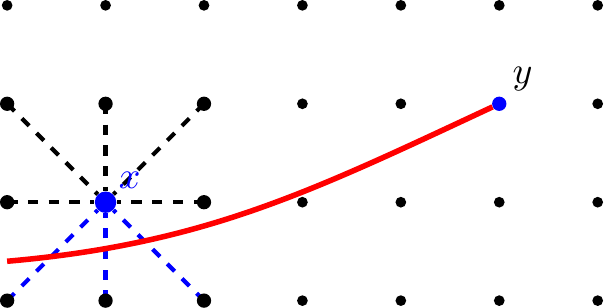}
\caption{Possible choices of $z$ for triangle updates. Dashed lines indicate bases used for triangle updates, while blue dashed lines indicate bases that should result in successful triangle updates.
 A triangle update is performed for each $z$ in the 8-point neighborhood of $x$, as in Algorithm \ref{alg:quasi:update}.
}
\label{fig:quasi:updates}
\end{figure}


\subsection{Practical difficulties}
\label{sec:quasi:obstacle}
There are several key algorithmic difficulties that arise in EJM but are not present in the previous first-order quasipotential solvers  \cite{cameron12,cameron18,dahiya18,yang19,lorenz63}. 
We describe these in detail in this section, in order to justify the decisions made in Algorithm \ref{alg:quasi:update} and the fail-safe method, discussed in the next section. 

\textbf{Large triangles.} First, we note that with cubic interpolation, the error in the approximation of $U(x_\lambda)$ in the triangle updates grows more quickly with respect to the length of the $xz$ line segment. If we suppose that $U(x)$, $U(z)$, $\nabla U(x)$ and $\nabla U(z)$ are equal to the correct values, then the 
error for the Hermite interpolation scales as $\mathcal{O}(\|x-z\|^4)$ whereas error of the linear interpolation grows as $\O(\|x-z\|^2)$. 

As such, it is critical to reduce the size of the $xz$ leg of triangle updates as much as possible. To this purpose, 
we opt to run a triangle update with $x$ and \textit{each} of its $8$-point nearest neighbors as $z$ (see Figure \ref{fig:quasi:updates}).
 This limits the $xz$ leg to a length of either $h$ or $\sqrt{2}h$, compared to a potentially much larger $xz$ length when neighboring elements of the stencil are used as in \cite{mirebeau14}.
 The cost of this improvement in accuracy, however, is the increase of the total number of calls for triangle update by a factor of around 3.

A more difficult but related problem to handle is the reduction of the size of the update length $\|x_\lambda-y\|$ in both the triangle updates and the one-point updates. 
Here, due to Simpson's quadrature rule we expect the local error of a 
$\eqref{eq:quasi:update_triangle_cubic}$ update to scale with update length as $\O(\|x_\lambda-y\|^5)$ and can be of either sign. 
There are cases where it is unavoidable to have a large update length $\|x_\lambda-y\|$ when the angle between the MAP and level set of $U$ is small (see Figure \ref{fig_ellball} for example). In fact, we see very clearly that in those situations the accuracy benefits of the EJM algorithm are significantly reduced (see Section \ref{sec:quasi:results}).

Our solution to minimizing this triangle leg is the introduction of item (B) in Accept/Reject Rule \ref{con:update}. 
Instead of using the rule $U(y):=\min(U(y),U_{\mathrm{new}}(y))$ of Algorithm \ref{alg:quasi:fastmarch}, we accept a new triangle update over a previous triangle update, 
only if the new update length $\|x_\lambda-y\|$ leg is smaller than the previous update length. 
We do not, however, require that the new triangle update proposes a smaller value of $U(y)$ 
then the previous triangle update value because the value with a longer update length may be artificially too small. 

In addition, to reduce the effect of the increased quadrature error of formulas $\eqref{eq:quasi:update_1pt_cubic}$ and 
$\eqref{eq:quasi:update_triangle_cubic}$ for large triangles, we increase the refinement of Simpson's rule quadrature. 
Specifically, once the minimizing path $\tilde{\varphi}_{\alpha^*,\beta^*,\lambda^*}$ of 
$\eqref{eq:quasi:update_triangle_cubic}$ is found, the \textit{actual} proposed update 
value $U_{\mathrm{new}}^{x,z}(y)$ is computed with a Simpson's rule approximation of 
$\eqref{eq:quasi:update_exact}$ along $\tilde{\varphi}_{\alpha^*,\beta^*,\lambda^*}$, 
but with a \textit{larger} number of nodes sampled along $\tilde{\varphi}_{\alpha^*,\beta^*,\lambda^*}$. 
For the number of nodes, we typically use the smallest odd number greater than $1 + \|x_{\lambda^*}-y\|/h$. 
Thus, for update lengths of $h$ and $\sqrt{2}h$, the number of nodes used in the Simpson's rule evaluation is the usual $3$, but for larger triangles, the number of nodes is proportional to $\|x_{\lambda^*}-y\|$. 
This refinement significantly reduces quadrature error when there is a large number of long-distance updates.

\textbf{Missed triangle updates.} Second, the vast reduction in the error from triangle updates $\eqref{eq:quasi:update_triangle_cubic}$ causes the error from one-point updates to be unacceptably large by comparison. The error discrepancy between the two types of updates is still present in earlier quasipotential solvers, but is of far smaller magnitude. In fact, with the EJM method, a small handful of mesh points whose final updates are one-point updates will destroy the expected $\O(h^2)$ error convergence rate. It is therefore essential that the final {\sf Accepted} value of $U(y)$ for each $y$ comes from a triangle update rather than a one-point update. This is precisely the purpose of the fail-safe discussed in Section \ref{sec_update_fail}. 

\begin{figure}
(a)\includegraphics{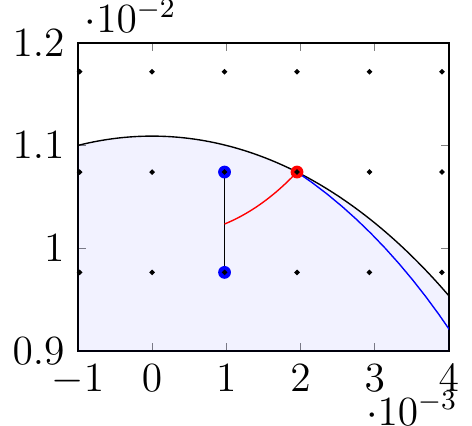}
(b)\includegraphics{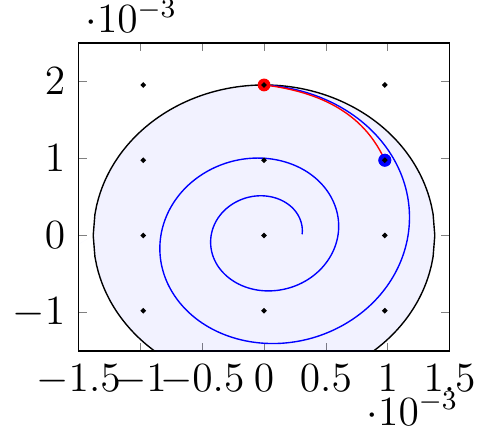}
\caption{Pitfalls. (a):  A situation where a local minimizer (red) { of the triangle update optimization problem $\eqref{eq:quasi:update_triangle_cubic}$} exists which does not correspond to the true MAP (blue).  (b): A situation where no successful 
triangle update occurs for a point $y$ near the attractor $\O$ because the curvature of the MAP (blue) is too great.}
\label{fig_fake}
\end{figure}

\textbf{Undesirable local minima.} Finally, there is a third more subtle problem introduced by adopting item (A) of Accept/Reject rule \ref{con:update}:  
namely, that the right-hand side of $\eqref{eq:quasi:update_triangle_cubic}$ may have local minima that have nothing to do with the desired MAP. An example of such a local minimum that does not correspond to a MAP is shown in Figure \ref{fig_fake}(a). Under the traditional rule $U(y):=\min(U(y),U_{\mathrm{new}}(y))$, the presence of these local minima is not an issue since these ``fake'' paths necessarily correspond to \textit{overestimates} of $U(y)$ and are hence filtered out. This is not the case when we just use item (B) as the accept/reject rule, since we are filtering based on update length.

To mitigate this problem, we introduce item (A) of Accept/Reject Rule \ref{con:update}. With (A), a valid triangle update must propose a value of $U(y)$ smaller than any previously proposed one-point update values of $U(y)$. This filters out any ``fake'' paths that lead to answers worse than the \textit{best} of the one-point updates. This is certainly not a guarantee that \textit{all} of the other possible local minima will be prevented, but rather, a guarantee that the ones that do slip by, will provide update values at least as good as the best one-point update.

\subsection{Fail-safe}
\label{sec_update_fail}

As mentioned above, it is critical that {\sf Accepted} mesh points have a final value of $U$ coming from a triangle update, rather than from a one-point update. To ensure this is the case we run the following fail-safe procedure on any mesh points whose final update is a one-point update. 

\begin{figure}
\includegraphics{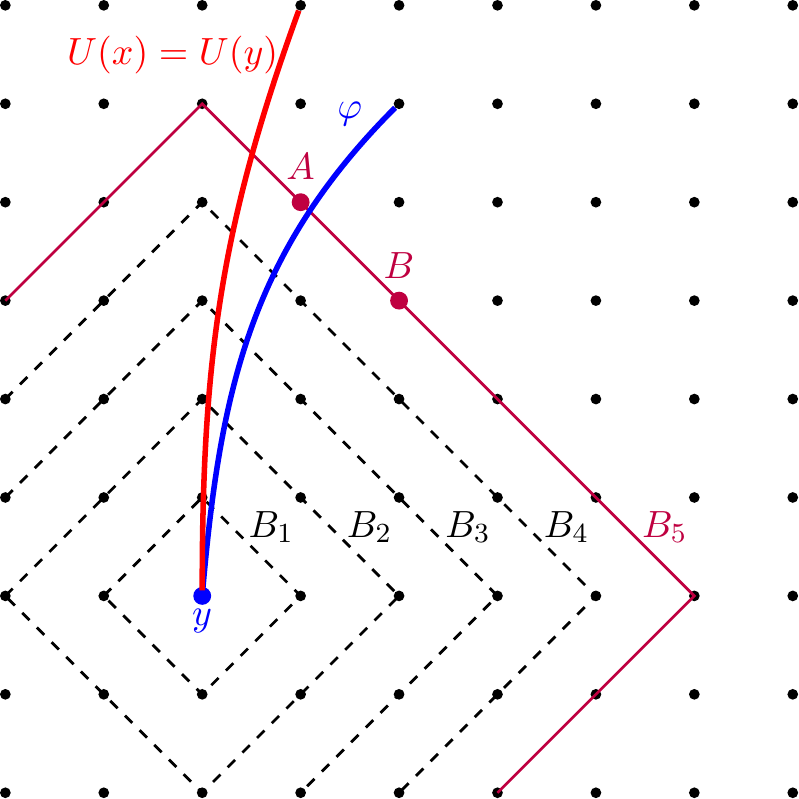}
\caption{The fail-safe will check $B_k$ for each $k$ until a successful two-point update is found. Here the first time a two-point update should be successful is on iteration $k = 5$ with $x = A$ and $z = B$.}
\label{fig_ellball}
\end{figure}

Let $x$ be the {\sf Considered} mesh point with the smallest value of $U$ in {\sf Considered} as in Algorithm \ref{alg:quasi:EJM}. Suppose that $x$'s most recent update is a one-point update. As shown in Figure \ref{fig_ellball}, let $B_k$ be the set of mesh points at distance $kh$ in $\ell^1$ norm away from $y$. We simply search each $B_k$, starting from $k =1$, for a pair of {\sf Accepted} neighbors $A$ and $B$ in $B_k$ that provide a successful triangle update of $x$. If we find such a pair, we use the associated triangle update values of $U(x)$ and $\nabla U(x)$. If the procedure does not succeed through some large $k$ threshold (we typically use $k = 20$), we revert to the prior one-point update value of $U(x)$. 

It is important to note that since the value of $U(x)$ is changed (and typically lowered) \textit{after} $x$ was selected from the {\sf Considered} list, we may lose the monotonic ordering of points in {\sf Accepted}. That is, a mesh point may have its value of $U$ lowered during the fail-safe procedure such that it is now smaller than some points already in the {\sf Accepted} list. This would have some significant adverse consequences if the fail-safe is called too frequently. However, provided that the fail-safe is called only for a small fraction of mesh points, 
these potential deviations from causality have a negligible impact. 

We note also that the presence of the fail-safe allows some flexibility in the creation of the stencils. As long as the stencils are \textit{reasonably} well refined in the direction of anisotropy, the fail-safe will correct for the few cases where the MAP might be able to sneak by the stencils.

%
%

\section{Performance tests}
\label{sec:quasi:results}

In this section, we set up a family of test problems with nonlinear drift fields with varying rotational components. We compare the performance of the following mesh-based quasi-potential solvers:
\begin{description}
\item[EJM,] the algorithm we have introduced in the preceding sections;
\item[FM-ASR,] the original algorithm of Mirebeau \cite{mirebeau14} which uses linear interpolation for $U$, linear approximation for MAP segments, the right-hand quadrature rule for approximating the geometric action integrals along the MAP segments, 
and  anisotropic stencils constructed from the asymmetric norms $\mathcal{F}^\alpha$ by Algorithm \ref{alg:quasi:stencil} 
(see Section \ref{sec:quasi:stencil});
\item[OLIM-midpoint,] the ordered line integral method \cite{cameron18} featuring an exhaustive search for MAPs with a ``hierarchical update strategy", 
linear interpolation for $U$, linear approximation for MAP segments, and midpoint quadrature rule;
\item[JM,] a jet marching algorithm which we have devised to pinpoint the gain in efficiency due to the targeted search for MAP segments in EJM; it has the same MAP search algorithm as the OLIMs  \cite{cameron18} and the same update procedure as EJM.
\end{description}

\subsection{Nonlinear drift with varying rotational components}
\begin{figure}
(a)\includegraphics{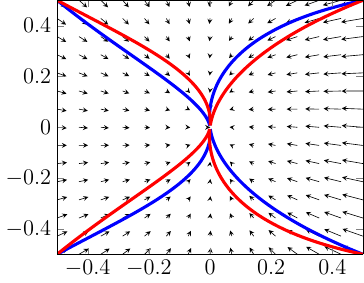}
(b)\includegraphics{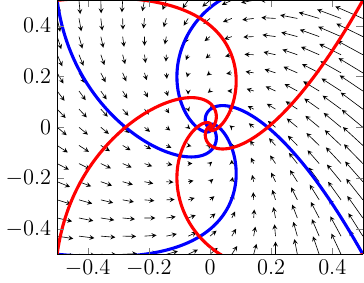}\\
(c)\includegraphics{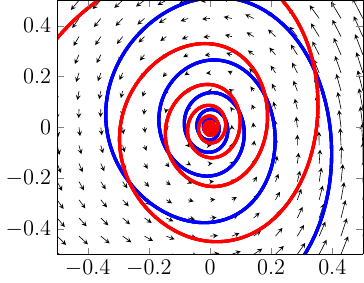}
\caption{Flow lines (blue) and MAPs (red) for drift field $\eqref{eq_drift_nonlinear}$ for three values of $a$. (a): $a=0.1$. (b): $a=1$. (c): $a=10$.}
\label{fig_drifts}
\end{figure}

We have set up a family of test problems with  the nonlinear drift field
\begin{equation}
\label{eq_drift_nonlinear}
	b(x,y) = 
	-\frac{1}{2}
	\begin{bmatrix}
	 4x + 3x^2 \\ 2y
	\end{bmatrix}
	+ \frac{a}{2} \begin{bmatrix}
		-2y
		\\ 4x + 3x^2
	\end{bmatrix}.
\end{equation}
The drift field $b$ admits a stable attracting equilibrium at $\O = (0,0)$ and a saddle at the point $x_0 = (-\sfrac{4}{3},0)$, for any value of $a$. 
The orthogonal decomposition $b = -\frac{1}{2} \nabla U + \ell$ with $\ell \perp \nabla U$ is apparent from \eqref{eq_drift_nonlinear}, so that the quasipotential 
in the basin of  $\O = (0,0)$ is $U(x,y) = 2x^2 + x^3 + y^2$.  
The parameter $a$ controls the rotational component of the drift field. The flow lines of $\eqref{eq_drift_nonlinear}$ 
(sample trajectories of the deterministic system $\dot{X} = b(X)$) are shown in Figure \ref{fig_drifts} with blue curves. 
For $|a|<2^{-3/2}\approx 0.35355$ the eigenvalues of the linearization of $b$ around $\O$ are both real, and there is no spiraling of the flow lines, while for $|a|> 2^{-3/2}$
the origin is a stable spiral point.  MAPs, the red curves in Figure \ref{fig_drifts} behave in the same way as the flow lines, but with a \textit{flipped} potential component.

We compare the performance of EJM with the three other solvers listed above for three different values of $a$: $a=0.1$, $a=1$, and  $a=10$. 
We compute the quasipotential $U$ and its gradient $\nabla U$ on the box $D = [-1,1] \times [-1,1]$ discretized into an 
 $N \times N$  square mesh with common horizontal and vertical mesh spacing $h = \sfrac{2}{(N-1)}$. 
The values of $N$ used in our experiments are: $ N = 2^k+1$ for $ k = 7,...,12$. 
Each solver is terminated once the first \textit{boundary} mesh point is added to {\sf Accepted}, corresponding to the first time the ``wave'' hits the edge of the box. Since this domain $D$ is within the well of attraction of $\O$, the solution $U(x) = 2x^2+x^3 + y^2$ is the quasipotential. We remark that outside the well of attraction, this is not the case.

\begin{figure}[!ht]
    \centering
    \begin{subfigure}[b]{\textwidth}
\centering
\includegraphics{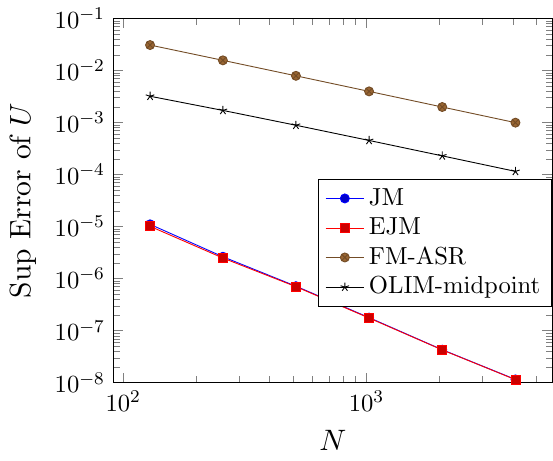}
\hfill
\includegraphics{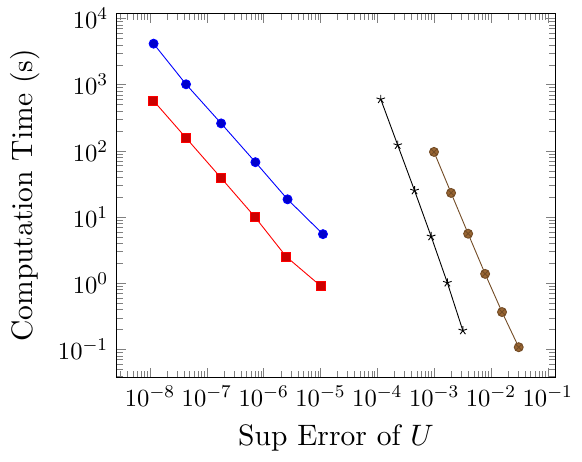}
\vspace{-5mm}
\caption{$a = 0.1$}
\label{fig:quasi:errs0}
\vspace{5mm}
\end{subfigure}
   \begin{subfigure}[b]{\textwidth}
\centering
\includegraphics{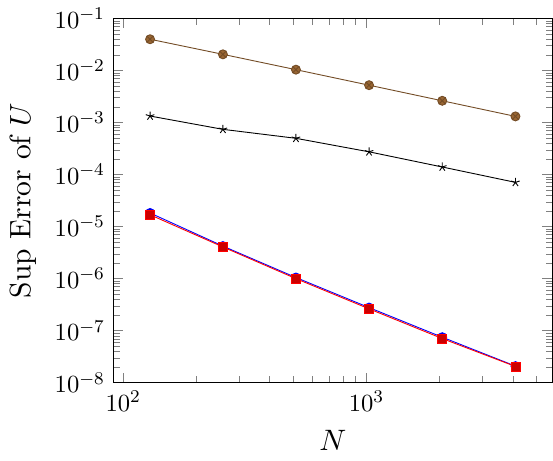}
\hfill
\includegraphics{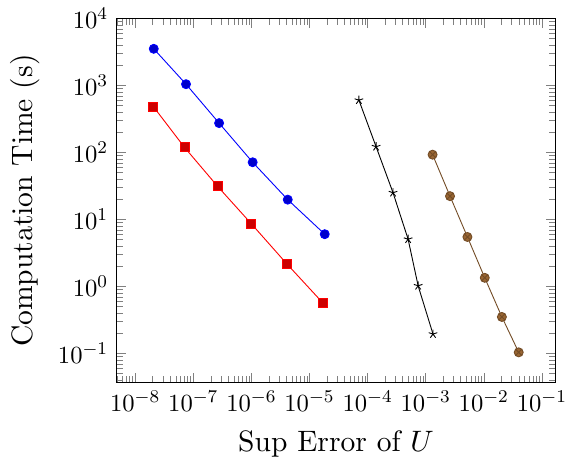}
\vspace{-5mm}
\caption{$a = 1$}
\label{fig:quasi:errs1}
\vspace{5mm}
\end{subfigure}
   \begin{subfigure}[b]{\textwidth}
\centering
\includegraphics{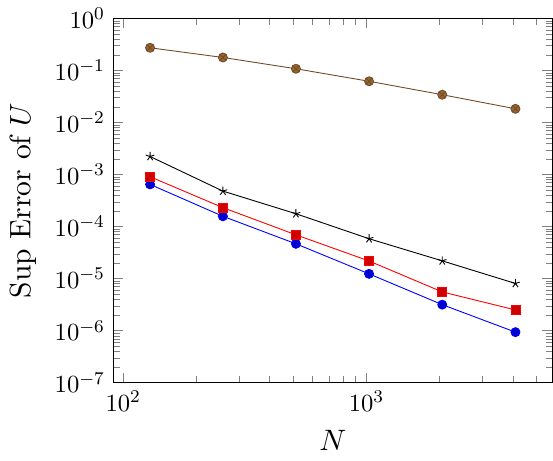}
\hfill
\includegraphics{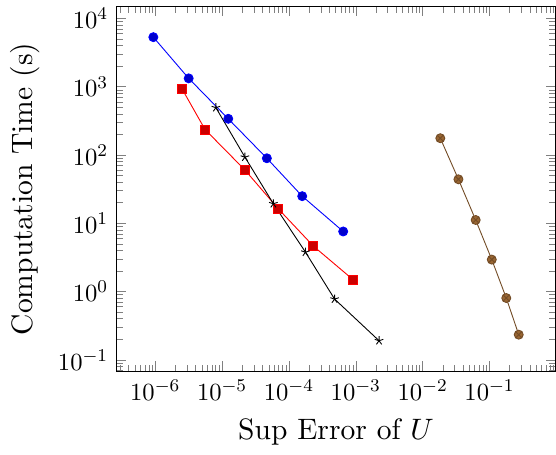}
\vspace{-5mm}
\caption{$a = 10$}
\label{fig:quasi:errs10}
\vspace{5mm}
\end{subfigure}
\caption[Error plots for $U$]{Error plots of $U$ for drift field $\eqref{eq_drift_nonlinear}$ for the three values of $a$. Sup error of $U$ over all Accepted mesh points is plotted against $N$ and computation time.}\label{fig:quasi:errs}
\end{figure}

\textbf{Errors in $U$.} In the left column of Figure \ref{fig:quasi:errs}, we plot the maximum error of $U$ over all {\sf Accepted} points versus $N$ for each of the 4 methods and each of the three values of $a$. Best polynomial fits were computed and are shown in Table \ref{table:quasi:u}. 

The most important takeaway is the clear 2nd order convergence rate in the errors for EJM and JM for the cases of \textit{mild} ($a=0.1)$ and \textit{moderate} ($a=1$) rotational components of the drift (left panel of Figures \ref{fig:quasi:errs0} and \ref{fig:quasi:errs1}). A convergence rate slightly below 2 is observed for the large rotational component at $a=10$.  This slight rate degradation is mostly due to the larger number of long-distance triangle updates associated with the high-rotational components, which cause relatively larger errors in our higher-order approximation scheme.  

 On the other hand, FM-ASR and OLIM-midpoint display the expected first-order convergence rates for $a = 0.1$ and $a = 1$. For the largest grid size $N =4096$ that we sampled, the errors in $U$ by the EJM and the JM are by four and five orders of magnitude smaller than those by OLIM-midpoint and FM-ASR, respectively. 
 
 As expected, the convergence rate for $a=10$ also degrades for FM-ASR. On the contrary, the convergence rate of  OLIM-midpoint increases to almost $1.5$. This curious phenomenon is investigated in \cite{yang19}. Nonetheless, the errors by EJM and JM are a few times smaller than those for OLIM-midpoint for all considered mesh sizes. The convergence rates for JM and EJM for $a=10$ are $1.89$ and $1.82$, respectively.

\begin{table}[h]
\centering
\begin{tabular}{|c|c|c|c|}
\hline
 & $a = 0.1$ & $a = 1$ & $a = 10$\\ \hline \hline 
EJM&$0.14\cdot N^{-1.96}$&$0.20\cdot N^{-1.94}$&$5.91\cdot N^{-1.82}$\\ 
JM&$0.17\cdot N^{-1.98}$&$0.22\cdot N^{-1.95}$&$6.03\cdot N^{-1.89}$\\ 
FM-ASR&$0.72\cdot N^{-0.99}$&$0.70\cdot N^{-0.99}$&$0.02\cdot N^{-0.61}$\\ 
OLIM-midpoint&$0.36\cdot N^{-0.96}$&$0.18\cdot N^{-0.94}$&$1.79\cdot N^{-1.48}$\\   \hline
\end{tabular}
\caption[Sup error of $U$ best fits]{Best fits for the supremum error of $U$ as a function of $N$.}\label{table:quasi:u}
\end{table}

\textbf{Runtime.} To analyze the utility of the anisotropic stencils, we compare the EJM with the JM method. As we see in the left panel of Figure \ref{fig:quasi:errs}, in the $a =0.1$ and $a=1$ case the errors of these two methods are indistinguishable, while in the $a =10$ case they are very close. This is expected since these two methods differ only in the choice of neighborhoods and not in the update procedures.

The important piece of information is the difference in runtime between these two methods, shown in the right panel of Figure \ref{fig:quasi:errs}. For the $a= 0.1$ and $a = 1$ cases, EJM is consistently between 7 and 8 times faster. In the $a = 10$ case, the improvement factor is only about $5$, likely due to a larger number of fail-safe calls due to stencil sparsity.

\textbf{Errors in $\nabla U$.} 
\begin{figure}[!ht]
    \centering
    \begin{subfigure}[b]{\textwidth}
\centering
\includegraphics{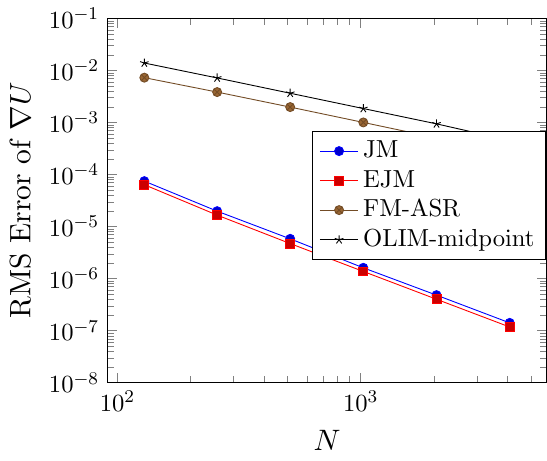}
\hfill
\includegraphics{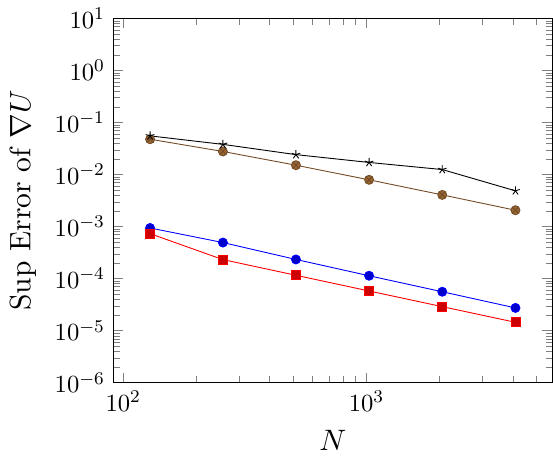}
\vspace{-5mm}
\caption{$a = 0.1$}
\label{fig:quasi:errs0:g}
\vspace{5mm}
\end{subfigure}
   \begin{subfigure}[b]{\textwidth}
\centering
\includegraphics{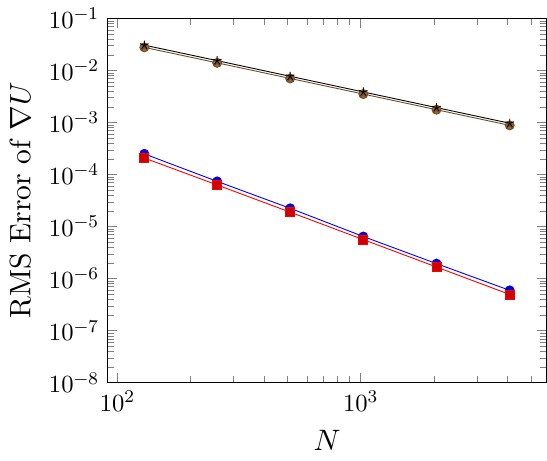}
\hfill
\includegraphics{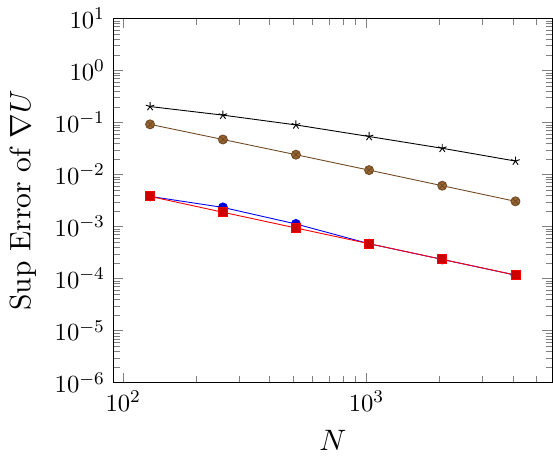}
\vspace{-5mm}
\caption{$a = 1$}
\label{fig:quasi:errs1:g}
\vspace{5mm}
\end{subfigure}
   \begin{subfigure}[b]{\textwidth}
\centering
\includegraphics{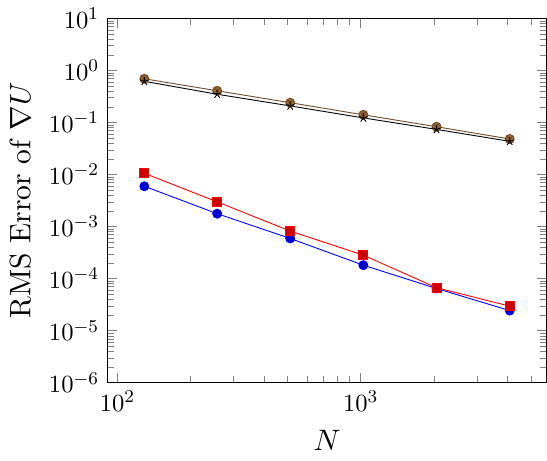}
\hfill
\includegraphics{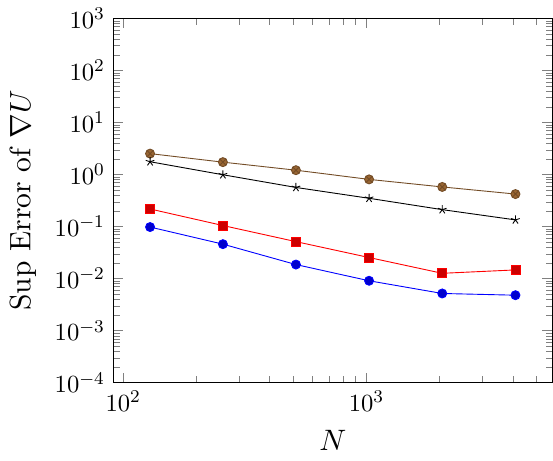}
\vspace{-5mm}
\caption{$a = 10$}
\label{fig:quasi:errs10:g}
\vspace{5mm}
\end{subfigure}
\caption[Error plots of $\nabla U$]{Error plots of $\nabla U$ for drift field $\eqref{eq_drift_nonlinear}$ for the three values of $a$. Sup error and RMS error of $\nabla U$ are plotted against $N$.}\label{fig:quasi:graderrs}
\end{figure}
Next we consider the error of $\nabla U$. In Figure \ref{fig:quasi:graderrs} we plot both the sup error (left panels) and root mean square (RMS) error (right panels) of $\nabla U$ against $N$. Best polynomial fits for these plots are included in Tables \ref{table:quasi:table_gradRMS} and \ref{table:quasi:table_gradSup}. The behavior of sup error and RMS error are markedly different. For the sup error (right panel), all 4 methods, including EJM, display linear in $h$ convergence rates. For the EJM, this is simply due to a few one-point updates making it into the final cut, and is not overly significant.

\begin{table}[h]
\centering
\begin{tabular}{|c|c|c|c|}
\hline
 & $a = 0.1$ & $a = 1$ & $a = 10$\\  \hline \hline 
EJM&$0.40\cdot N^{-1.81}$&$1.02\cdot N^{-1.75}$&$68.13\cdot N^{-1.81}$\\ 
JM&$0.47\cdot N^{-1.81}$&$1.22\cdot N^{-1.75}$&$12.74\cdot N^{-1.59}$\\ 
ASR-endpoint&$0.90\cdot N^{-0.82}$&$3.38\cdot N^{-0.98}$&$146.61\cdot N^{-1.08}$\\ 
OLIM-midpoint&$1.72\cdot N^{-0.99}$&$3.95\cdot N^{-1.00}$&$22.69\cdot N^{-0.75}$\\    \hline
\end{tabular}
\caption[RMS error of $\nabla U$ best fits]{Best fits for the RMS error of $\nabla U$ as a function of $N$.}
\label{table:quasi:table_gradRMS}
\end{table}\vspace{-.5cm}

\begin{table}[h]
\centering
\begin{tabular}{|c|c|c|c|}
\hline
 & $a = 0.1$ & $a = 1$ & $a = 10$\\ \hline \hline 
EJM&$0.12\cdot N^{-1.10}$&$0.50\cdot N^{-1.01}$&$31.81\cdot N^{-1.03}$\\ 
JM&$0.14\cdot N^{-1.03}$&$0.67\cdot N^{-1.04}$&$7.13\cdot N^{-0.92}$\\ 
ASR-endpoint&$14.48\cdot N^{-0.92}$&$66.24\cdot N^{-1.06}$&$578.00\cdot N^{-1.04}$\\ 
OLIM-midpoint&$1.44\cdot N^{-0.65}$&$10.98\cdot N^{-0.77}$&$51.79\cdot N^{-0.72}$\\    \hline
\end{tabular}
\caption[Sup error of $\nabla U$ best fits]{Best fits for the supremum error of $\nabla U$ as a function of $N$.}
\label{table:quasi:table_gradSup}
\end{table}

The RMS errors for EJM, on the other hand, display higher order convergence. Indeed, the best fit convergence rates (Table \ref{table:quasi:table_gradRMS}) are $1.81$, $1.75$ and $1.81$ for $a = 0.1$, $1$ and $10$, respectively. These are close to 2nd order.  On the other hand, OLIM-midpoint displays the expected 1st order in $h$ convergence of $\nabla U$ associated with linear methods which is as good as one can expect from a method with linear approximation to MAP segments. Careful numerical analysis of EJM aiming at establishing convergence rates for $U$ and $\nabla U$ analytically is complicated and left for future work.


\section{The Maier-Stein SDE: sharp estimates for the invariant measure and for the escape time}
\label{sec:MS}
We consider the SDE
\begin{equation}
\label{eq:MS}
	d\left[\begin{array}{c}X_t\\Y_y\end{array}\right] = \begin{bmatrix}
	x-x^3-\beta x y^2 \\ -(1+x^2)y
\end{bmatrix} dt + \sqrt{\e}dW_t, 
\end{equation}
where $\beta>0$ is a parameter, and $dW_t$ is the two-dimensional Brownian motion.
This SDE was introduced by R. Maier and D. Stein in 1993 \cite{MS1993} as an example of a system with a bistable nongradient drift field. The invariance with respect to $y\leftrightarrow-y$ makes it amenable to analysis by means of asymptotic techniques. The system is also symmetric with respect to $x\leftrightarrow-x$. For any $\beta > 0$, the drift $\eqref{eq:MS}$ admits stable attracting equilibria at $\O_-:=(-1,0)$ and $\O_+:=(1,0)$ and a saddle at $\O_*:=(0,0)$. The $y$-axis is the boundary of the basins of attraction of $O_-$ and $\O_+$.

A detailed study of this system is offered by Maier and Stein in a long paper of 1996 \cite{maier96}. 
They considered the \emph{Lagrangian manifold} consisting of MAPs emanating from $\O_+$ in the four-dimensional \emph{lifted space}.
They showed that as the parameter $\beta$ increases, the Lagrangian manifold undergoes a number of bifurcations. If $\beta=1$, the system is gradient with a potential given by
\begin{equation}\label{MSpot}
V(x) = \frac{x^4}{4}-\frac{x^2}{2}+\frac{x^2y^2}{2}+\frac{y^2}{2}+\frac{1}{4}.
\end{equation}
If $1<\beta<4$, the MAPS emanating from $\O_+$ do not intersect, the Lagrangian manifold has no folds, the quasipotential is smooth,
and the WKB approximation \eqref{eq:WKB} for the invariant probability measure is valid. In particular, the MAP connecting $\O_+$ and the saddle $\O_*$ is simply a segment along the $x$-axis. For $\beta>4$, the Lagrangian manifold acquires folds, corresponding to locations where MAPs intersect along the $x$-axis. In fact, the MAP connecting $\O_+$ with the saddle $\O_*$ bifurcates into two symmetric MAPs arching above and below the $x$-axis, respectively. As a result, the quasipotential is nondifferentiable along a certain segment of the $x$-axis for each $\beta>4$. Indeed, the $j$th bifurcation of the Lagrangian manifold (consisting of adding the $j$th fold) occurs at $\beta = (j+1)^3$, $j = 1,2,\ldots$. These bifurcations lead to a failure of the WKB approximation due to a blow-up of the prefactor $C(x)$ at the points where MAPs intersect.

The Maier-Stein SDE has become a benchmark example for testing various methods for computing transition paths and the quasipotential \cite{heymannCPAM,cameron12,dahiya18,kikuchi20}. We use EJM to compute the quasipotential $U$ with respect to $\O_-$ at two values of $\beta$: $\beta = 3$ and $\beta = 10$. At  $\beta = 3$, the system is nongradient but the quasipotential is smooth and the WKB approximation is expected to work. At $\beta = 10$, $\nabla U$ is discontinuous along the segment $[x^\star,0]$ of the $x$-axis where $x^\star\approx -0.668$, and the WKB prefactor is expected to blow up \cite{maier96}. Besides the quasipotential, we also compute the prefactor for the WKB approximation \eqref{eq:WKB} by solving the transport equation \eqref{eq:transport}. This is the first work offering a numerical solution for the WKB prefactor to the best of our knowledge. In addition, we calculate a sharp estimate for the expected escape time from the basin of $\O_-$ using the Bouchet-Reygner formula \cite{bouchet15}. To validate our results, we compare them with the invariant measure and the expected transition time computed using the machinery of the \emph{transition path theory} \cite{eve2006}.

\subsection{The WKB approximation for the invariant measure}
\label{sec:MS_WKB}
Using EJM, we compute the quasipotential with respect to the equilibrium $\O_-=(-1,0)$ for $\beta=3$ and $\beta=10$ in the box $(x,y)\in[-2,0]\times[-1,1]$ discretized to a $2049\times 2049$ mesh. The computation is stopped as soon as the computation reaches the saddle $\O_*$ at the origin. 

\begin{figure}
(a)\includegraphics[width=0.4\textwidth]{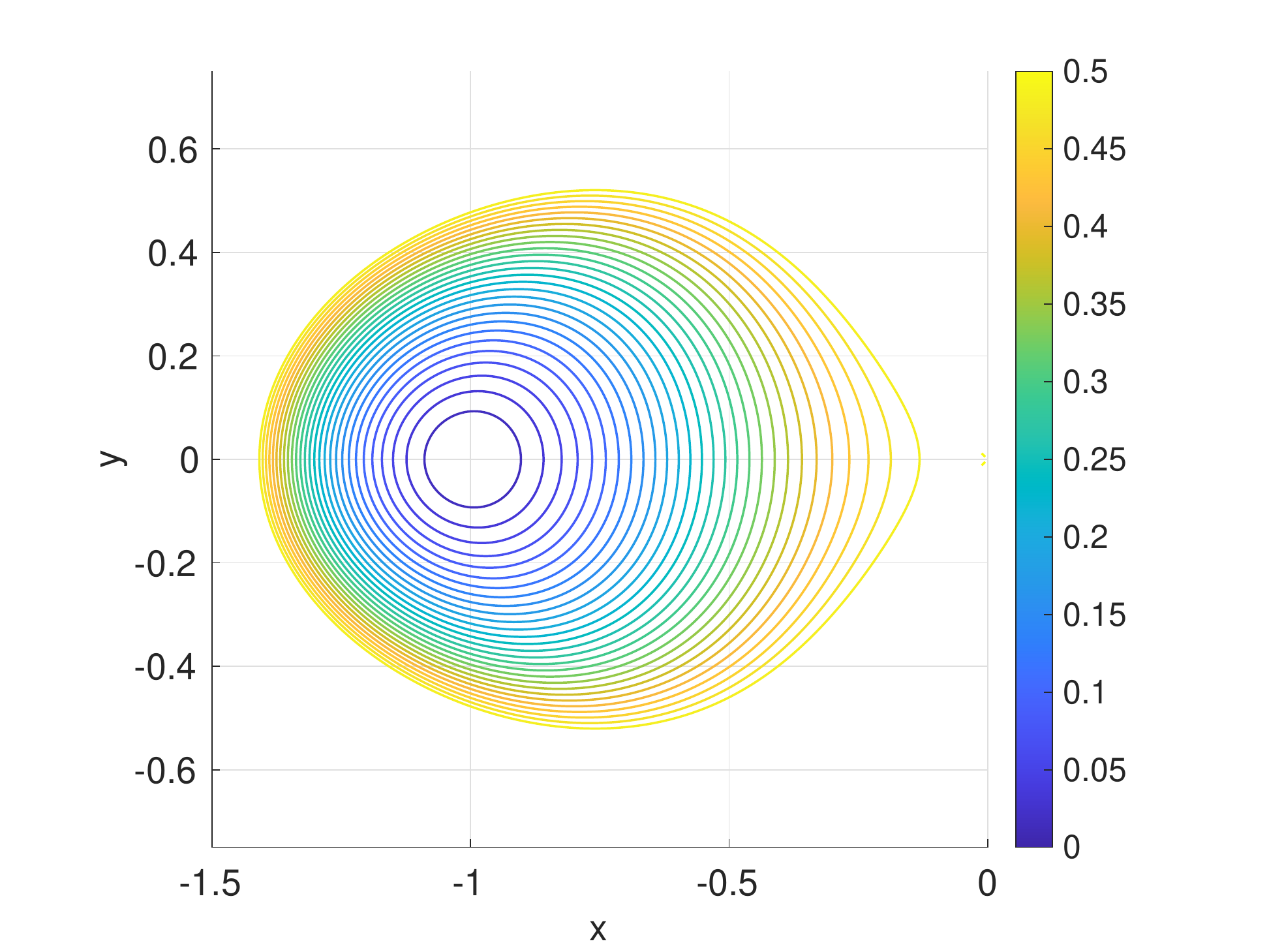}
(b)\includegraphics[width=0.4\textwidth]{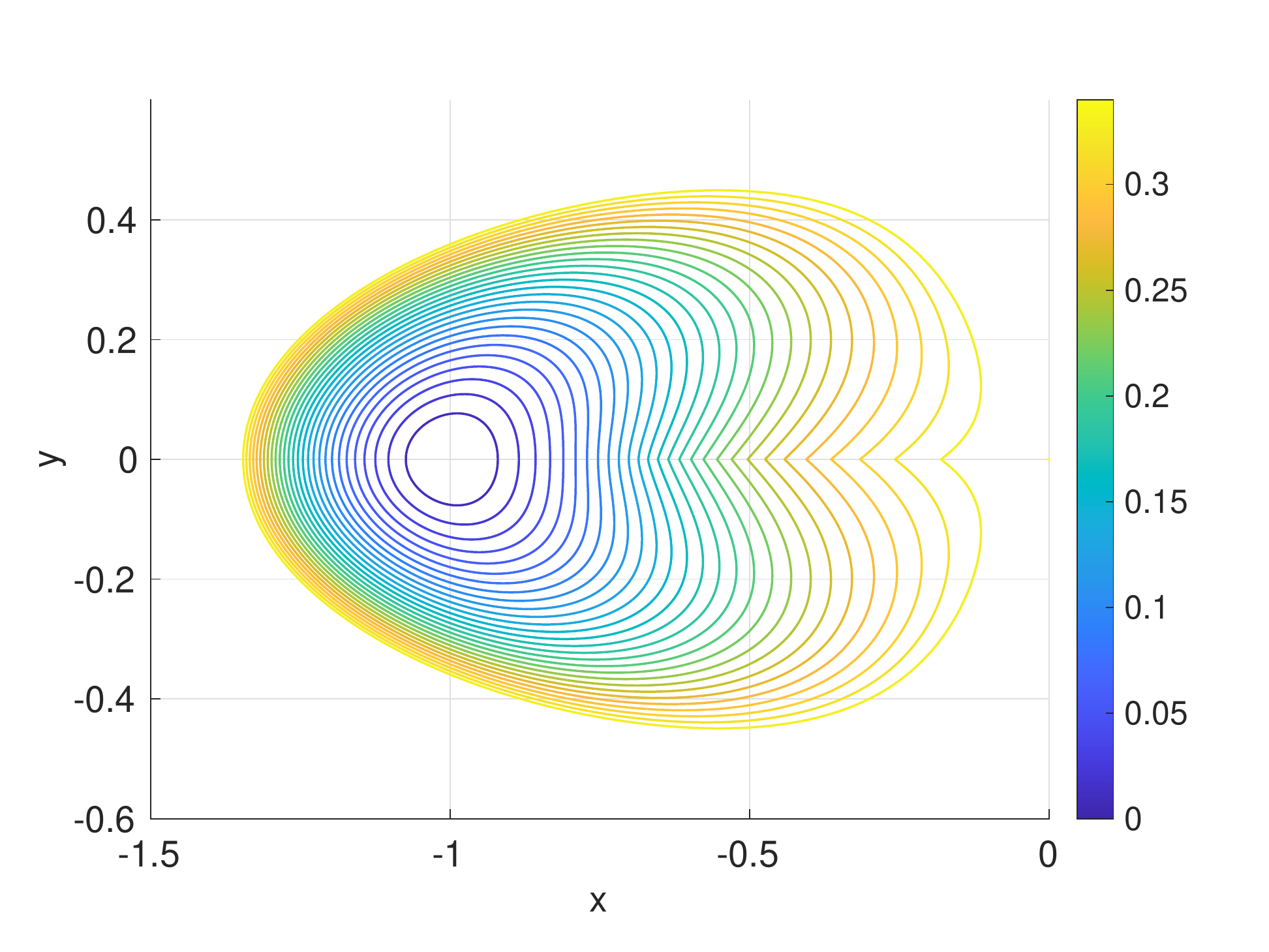}
\caption{The level sets of the quasipotential for the Maier-Stein SDE \eqref{eq:MS} for $\beta=3$ (a) and $\beta = 10$ (b).
}
\label{fig:MSqpot}
\end{figure}

The value of the quasi-potential at the origin $\O_*$ computed by EJM for $\beta=3$ is equal to $0.5$ exactly up to the machine precision. This matches the exact value, which can be computed analytically since the MAP connecting $(-1,0)$ and the origin is a segment of the $x$-axis. Indeed, the computation of the geometric action \eqref{eq_action_geometric} along this MAP \cite{maier96} is given by:
$$
U(\O_*) =\int_{-1}^0\|b(x,0)\|-b(x,0)dx = \int_{-1}^0 2(x-x^3)dx = 0.5.
$$
The computed value of the quasipotential by EJM at the saddle $\O_*$ for $\beta=10$ is $0.340039$. The level sets of the quasipotential for $\beta=3$ and $\beta=10$ are shown in Figure \ref{fig:MSqpot}.
The numerical tests of Section \ref{sec:quasi:results} suggest that the RMS error for the gradient $\nabla U$ decays superlinearly as $h\rightarrow 0$. Therefore, we can afford differentiating the rotational component $l :=b+\tfrac{1}{2}\nabla U$ numerically using central differences and calculate the divergence of the rotational component $\nabla\cdot l$. We observe that the transport equation for the prefactor \eqref{eq:transport} can be written as 
\begin{equation}
    \label{eq:C1}
    (\nabla \cdot l)C + \dot{\phi}\cdot\nabla C = 0,
\end{equation}
where $\dot{ \phi} = b+\tfrac{1}{2}\nabla U$ is the velocity of the MAP with respect to time. Therefore $\dot{\phi}\cdot\nabla C$ is merely the time derivative of $C$ along the MAP.
This allows us to rewrite \eqref{eq:C1} along the MAP as $\tfrac{d}{dt}C+(\nabla \cdot l)C  = 0$ as
\begin{equation}
    \label{eq:C2}
    C(t) = C_0\exp\left\{-\int_0^t(\nabla \cdot l)dt\right\}.
\end{equation}
It is more convenient to integrate along the MAPs with respect to their arclength $dr$ and prescribe endpoints in space rather than in time. Observing that $dr = \|\dot{\phi}\|dt = \|b\|dt$, we obtain the following formula for the WKB prefactor $C$ at the point $(x,y)$ as an integral along the MAP $\varphi$ parametrized by its arclength passing through $(x_0,y_0)$ and arriving at $(x,y)$:
\begin{equation}
    \label{eq:C3}
    C(x,y) = C(x_0,y_0)\exp\left\{-\int_0^L\frac{(\nabla \cdot l(\varphi))}{\|b(\varphi)\|}dr\right\}.
\end{equation}

Thus, we need to propagate the prefactor along the MAPs emanating from $\O_-$. 
To do so numerically, we first we compute the function 
\begin{equation}
    \label{eq:fint}
f(z):=\frac{\nabla \cdot l(z)}{\|b(z)\|},
\end{equation}
on the mesh by differentiating $l = b + \tfrac{1}{2}\nabla U$ using central differences.
Second, we need to initialize this computation in the neighborhood of the attractor.
We observe that the quasipotential can be approximated by the quasipotential for the linearized SDE
\begin{equation}
    \label{eq:lin}
dZ_t = JZ_tdt + \sqrt{\e}dW_t,\quad Z_t=[X_t,Y_t]^\top,
\quad \{J\}_{ij} = \left(\frac{\partial b_i}{\partial Z_j}(\O_-)\right),
\end{equation}
with accuracy $\O(\rho^3)$ in a ball of radius $\rho$ centered at $\O_-$. 
For the linear SDE \eqref{eq:lin}, the quasipotential is of the form $U(z) = z^\top Qz$, and the potential and rotational components are, respectively $Qz$ and $Rz$. The matrices $Q$ and $R$ satisfy: $Q$ is symmetric positive definite, $-Q+R = J$, and  $QR$ is antisymmetric. One can show that ${\sf tr}R = 0$ (see Corollary 5.1 in \cite{cameron17}) which implies that the rotational component is divergence-free. Indeed, 
$$
\nabla\cdot l(z) = \nabla\cdot( Rz) = {\sf tr} R = 0.
$$
Therefore, $C(z)$ is constant for the linearized system in the neighborhood of the attractor. Since the constant is arbitrary, we set it to $1$.
This suggests a simple initialization for the computation of $C$: we pick a level set of the quasipotential $U$ corresponding to some small value and initialize $C=1$ inside this level set. We used the level set $U=0.001$.

\begin{figure}
(a)\includegraphics[width=0.4\textwidth]{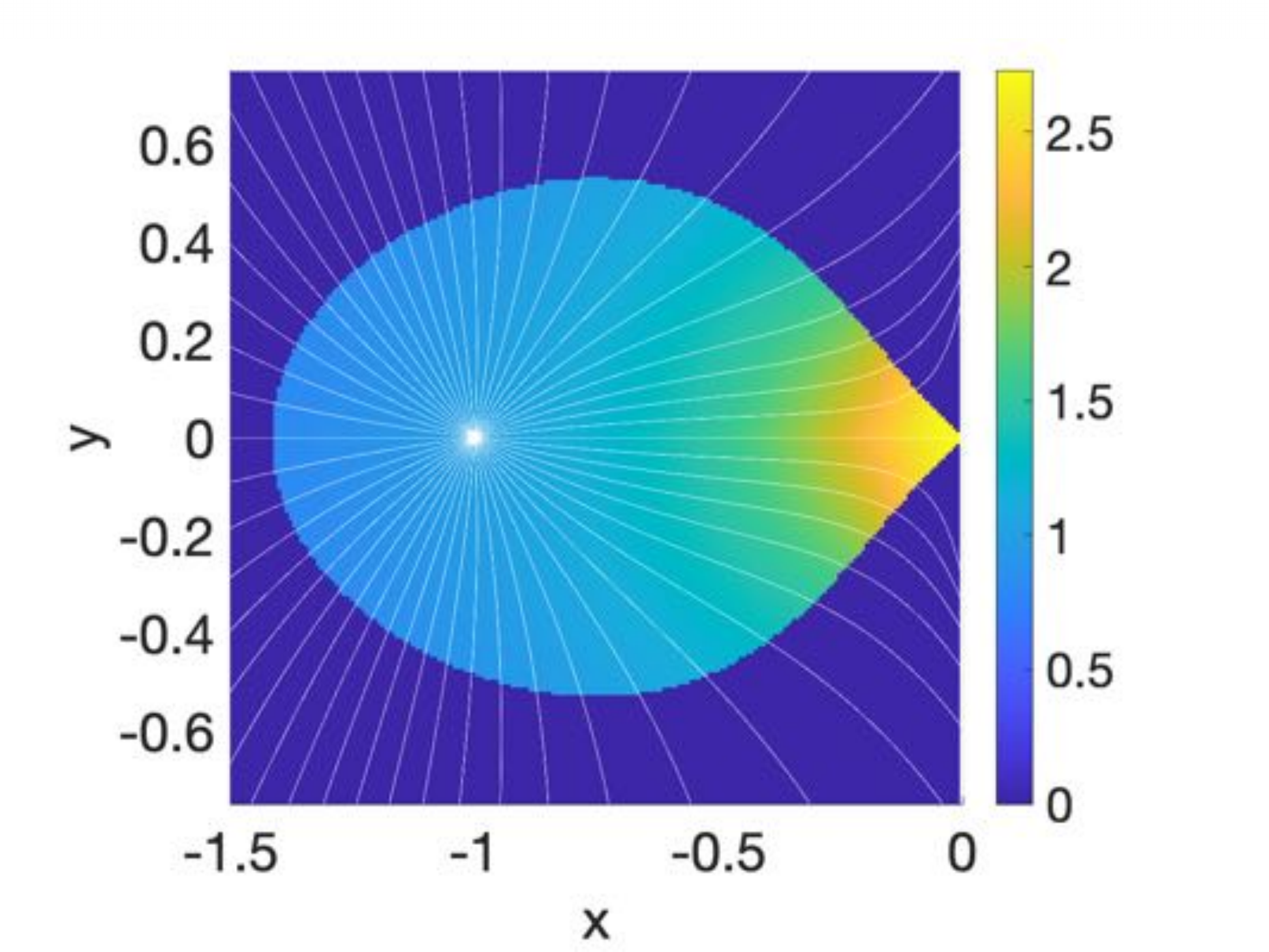}
(b)\includegraphics[width=0.4\textwidth]{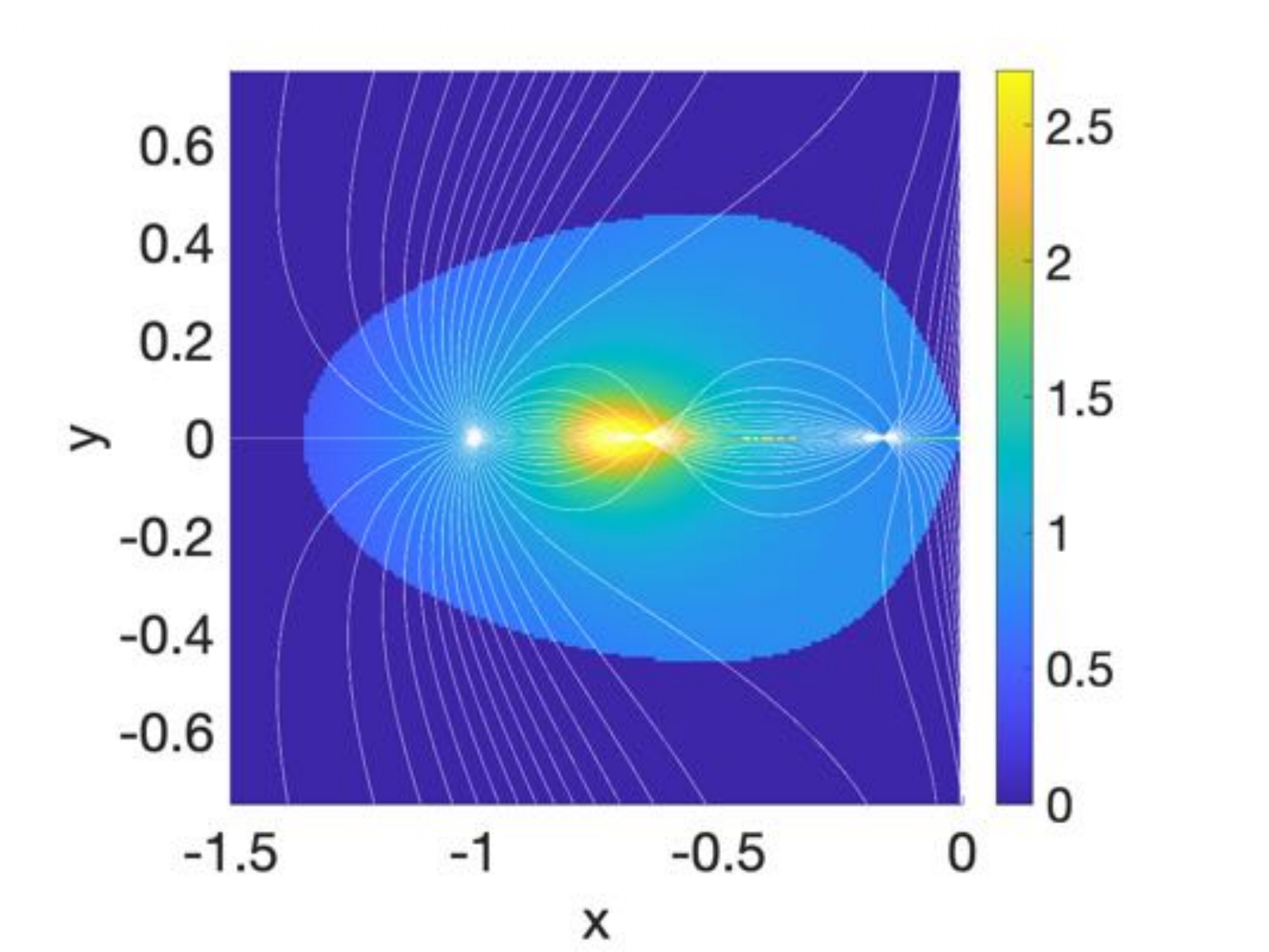}
\caption{The prefactors for the WKB approximation for the invariant measure for the Maier-Stein SDE \eqref{eq:MS} with $\beta=3$ (a) and $\beta = 10$ (b). Thin white curves are the MAPs emanating from $(-1,0)$ computed as described in \cite{aplot2005}.
}
\label{fig:WKBpref}
\end{figure}

Finally, we use the parent-child relationship from the computation of the quasipotential $U$ to compute values of $C$ for the rest of the mesh points exactly in the same order as they were added to {\sf Accepted} in the original EJM procedure.
Each mesh point $z = (x,y)$ whose value was finalized by a triangle update, has two parents, $p_1$ and $p_2$, and a value of the parameter $\lambda$ determining the onset of the MAP segment arriving at $z$ from the interval $[p_1,p_2]$: $p_{\lambda}: = p_1+\lambda(p_2-p_1)$. We use linear interpolation to approximate $C$ at $p_{\lambda}$: $C(p_{\lambda}) = C(p_1)+\lambda (C(p_2)-C(p_1))$ and the backward Euler scheme to approximate $C(z)$ from $\eqref{eq:C3}$
\begin{equation}
    \label{eq:Cscheme}
    C(z) = C(p_\lambda)- \|z-z_{\lambda}\|f(z)C(z),
~~{\rm i.e.,}~~
     C(z) = \frac{C(p_{\lambda})}{f(z)\|z-z_{\lambda}\|}.
\end{equation}
We have chosen this scheme due to its stability. Of course, we could calculate the integral in \eqref{eq:C3} using Simpson's rule and take the exponent of it, but combined with the linear interpolation to obtain $C_0\equiv C(p_{\lambda})$ this results in a numerical instability.
The computed prefactors for $\beta=3$ and $\beta = 10$ are shown in Figure \ref{fig:WKBpref}.
It is evident that the computed prefactor $C$ for $\beta = 10$ blows up\footnotemark[3] exactly at the point where the MAPs cross, as predicted in \cite{maier96}.\footnotetext[3]{The maximum of the numerical solution for $C$ is about $8\cdot10^3$.}

\begin{figure}
(a)\includegraphics[width=0.4\textwidth]{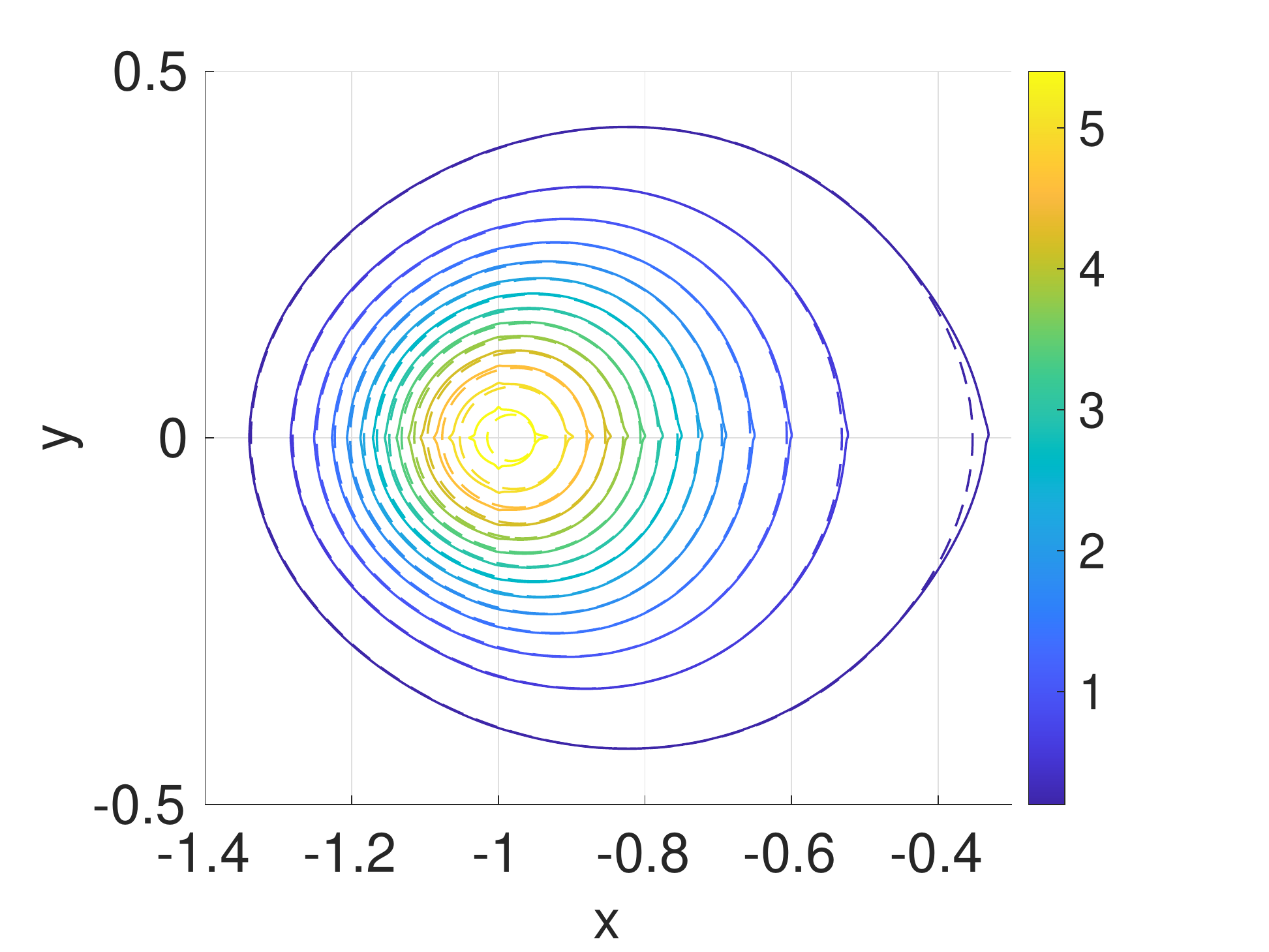}
(b)\includegraphics[width=0.4\textwidth]{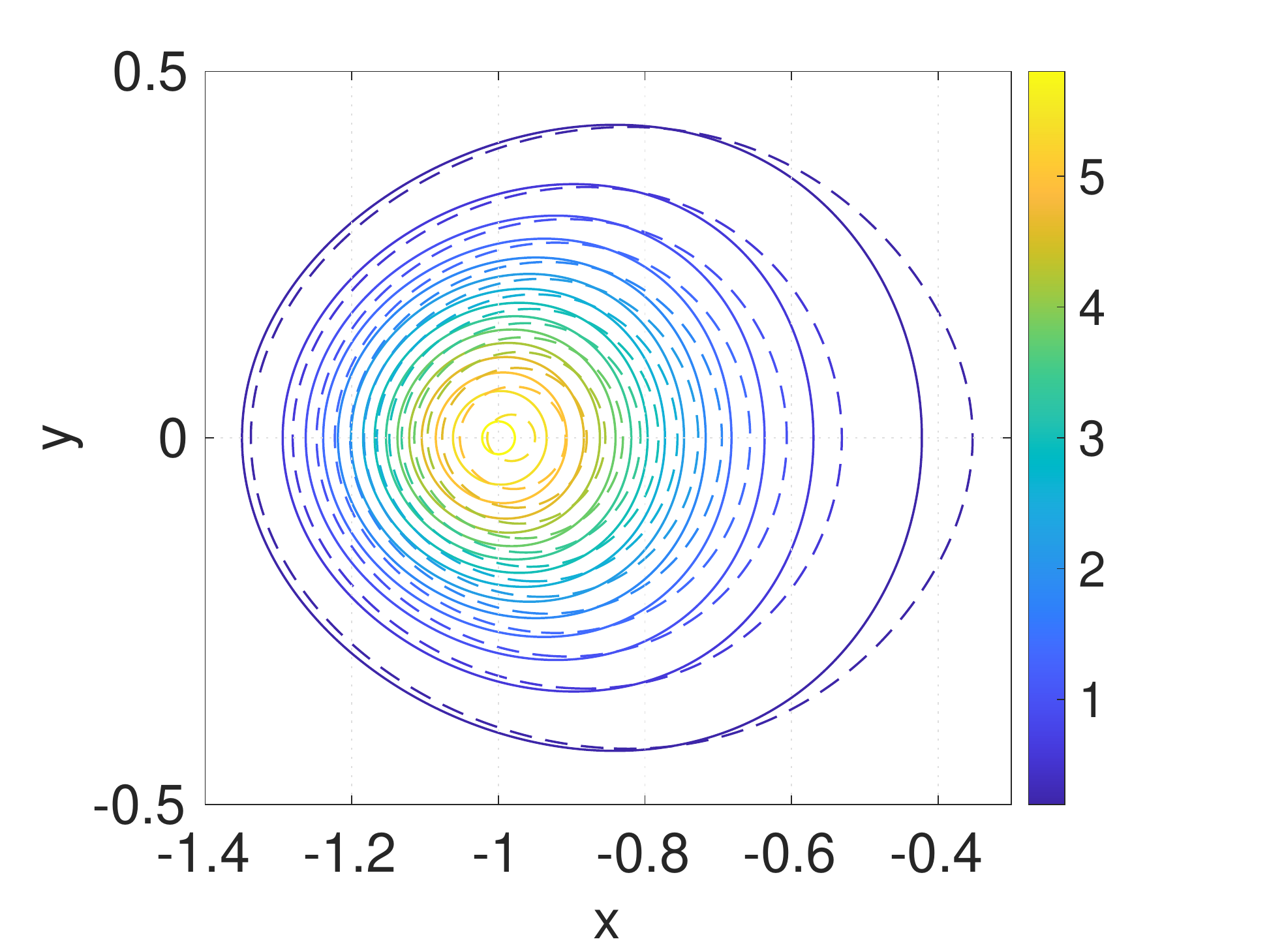}
\caption{(a): The level sets of the WKB approximation for the invariant measure for the Maier-Stein SDE \eqref{eq:MS} with $\beta=3$ (solid curves) and the level sets of the TPT invariant measure (dashed curves).   (b): The level sets of $Z^{-1}e^{-U/\epsilon}$ (solid curves) and the level sets of the TPT invariant measure (dashed curves). In all cases, $\epsilon = 0.1$.
}
\label{fig:WKBmu}
\end{figure}

We pick $\epsilon = 0.1$ and normalize the WKB approximation to the invariant measure 
$$
\mu_{WKB}(z) = C(z)e^{-U(z)/\epsilon},
$$
so that its integral over the mesh is 1. For comparison, we also compute the invariant measure using the transition path theory (TPT) \cite{eve2006} (see Appendix B for details).
The results for $\beta = 3$ are displayed in Figure \ref{fig:WKBmu}(a). To emphasize the effect of the prefactor, we also plot the level sets of $Z^{-1}e^{-U(z)/\epsilon}$ in Figure \ref{fig:WKBmu}(b). 
The invariant measure computed using the WKB approximation  is in a good agreement with the TPT invariant measure. The discrepancy between $Z^{-1}e^{-U(z)/\epsilon}$  and the TPT invariant measure is notably larger.

\subsection{The Bouchet-Reygner sharp estimate for the expected escape time}
Accurate computation of the quasipotential and its gradient allows us to obtain a sharp estimate for the expected escape time using the Bouchet-Reygner formula.
If the quasipotential is twice continuously differentiable at the saddle point $\O_*$ between two stable equilibria $\O_-$ and $\O_+$ in a bistable system, Bouchet and Reygner \cite{bouchet15} showed that the expected escape time from the basin of $\O_-$ is estimated by
\begin{equation}
    \label{eq:quasi:bouchet}
	\E \tau^{\O_- \to \O_*} \approx \frac{2\pi}{\lambda_+^*} \sqrt{\frac{|\mathrm{det} H(\O_*)|}{\mathrm{det} H(\O_-)}} \exp \Big(\int_0^L f(\varphi^*_r)dr \Big) \exp \Big( \frac{U(\O_*)}{\e}\Big),
\end{equation}
where $H$ is the Hessian of the quasipotential, $\lambda_+^*$ is the unique positive eigenvalue of the Jacobian matrix of the drift field at the saddle $\O_*$, $\varphi^*$ is the MAP from $\O_-$ to $\O_*$, and $f$ is given by \eqref{eq:fint}. We remark that we have replaced the integral in \eqref{eq:quasi:bouchet} with respect to time $t$ with the one with respect to the arclength $r$\footnotemark[4]. \footnotetext[4]{There is an error in equation (10) in \cite{dahiya18}: the integrand of the integral should be divided by $\|b + A^{-1} \nabla U\|$.}

\begin{figure}
(a)\includegraphics[width=0.4\textwidth]{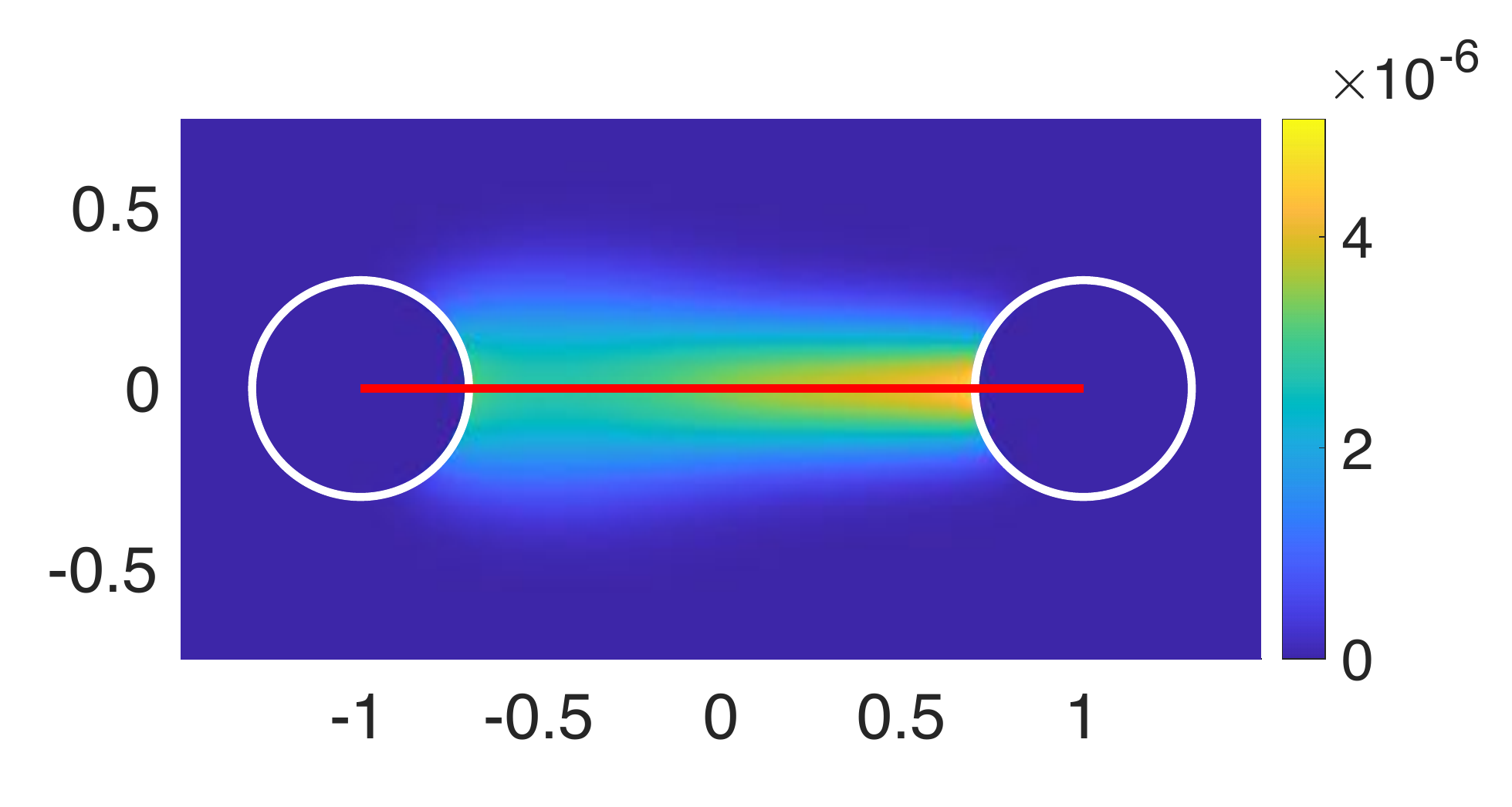}
(b)\includegraphics[width=0.4\textwidth]{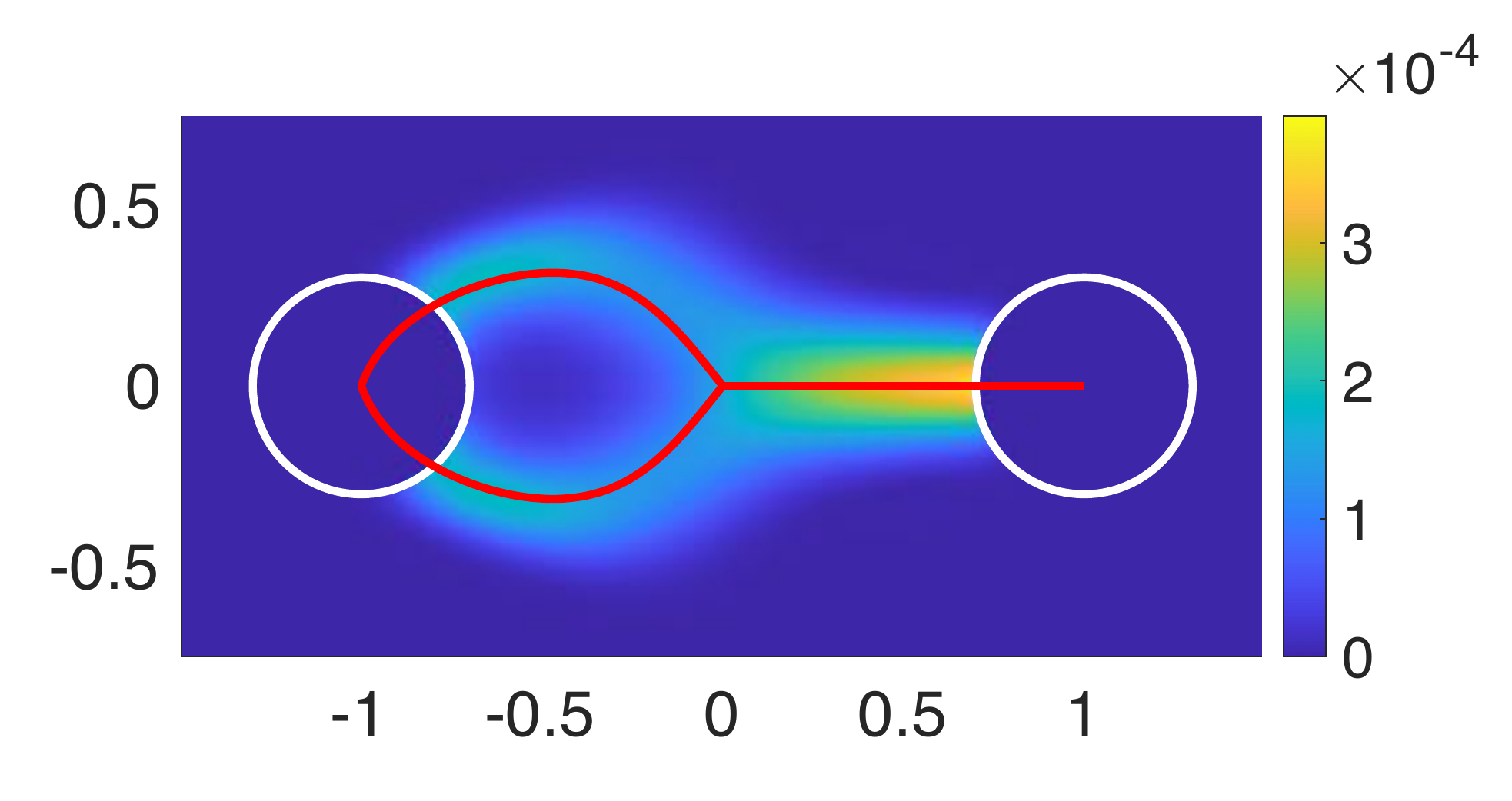}
\caption{The magnitude of the TPT reactive current \cite{eve2006} at $\epsilon=0.04$ and the MAPs connecting $\O_-$ and $\O_+$ for $\beta=3$ (a) and $\beta = 10$ (b).
}
\label{fig:Rcurrent}
\end{figure}

\begin{figure}
(a)\includegraphics[width=0.4\textwidth]{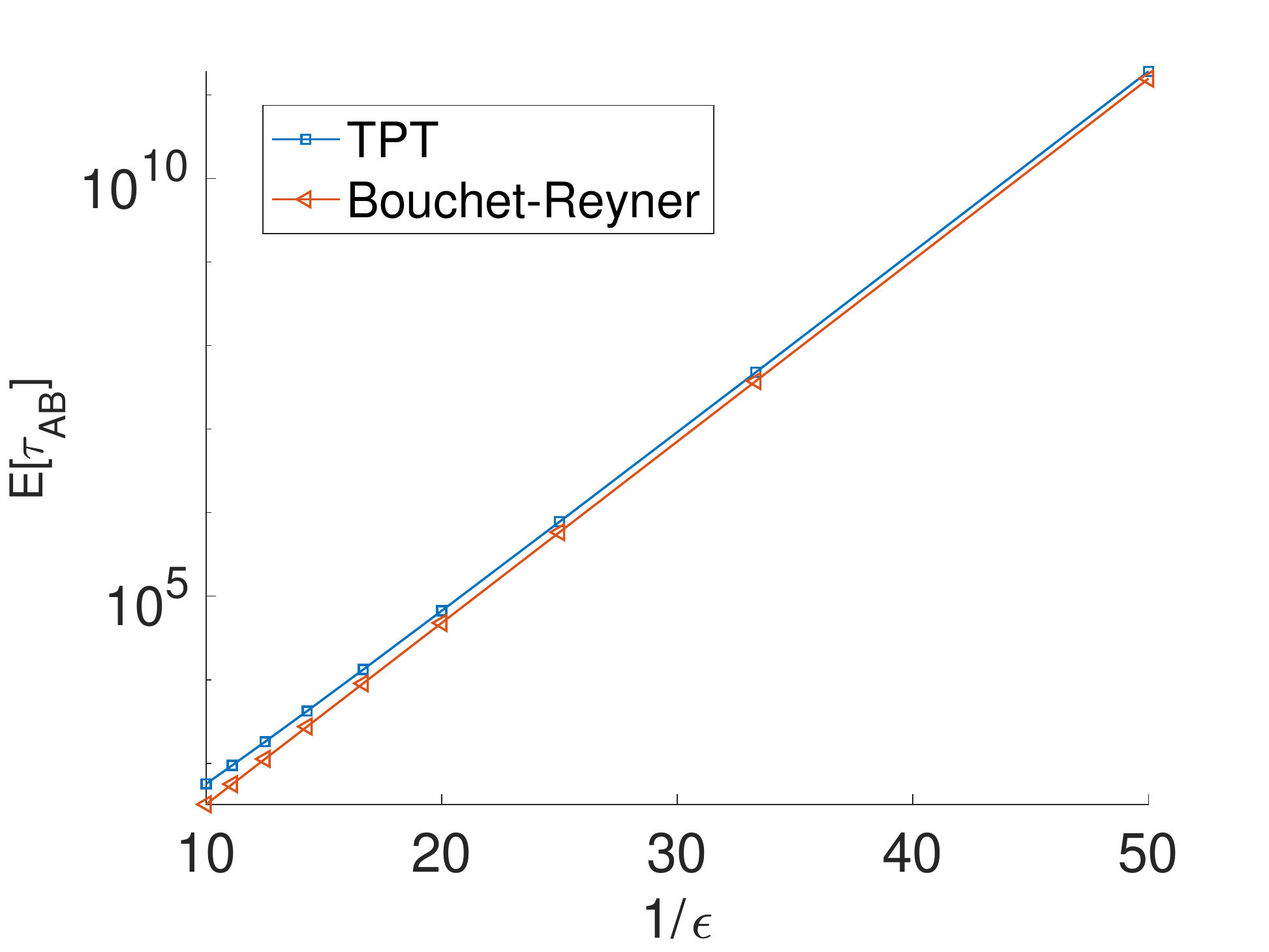}
(b)\includegraphics[width=0.4\textwidth]{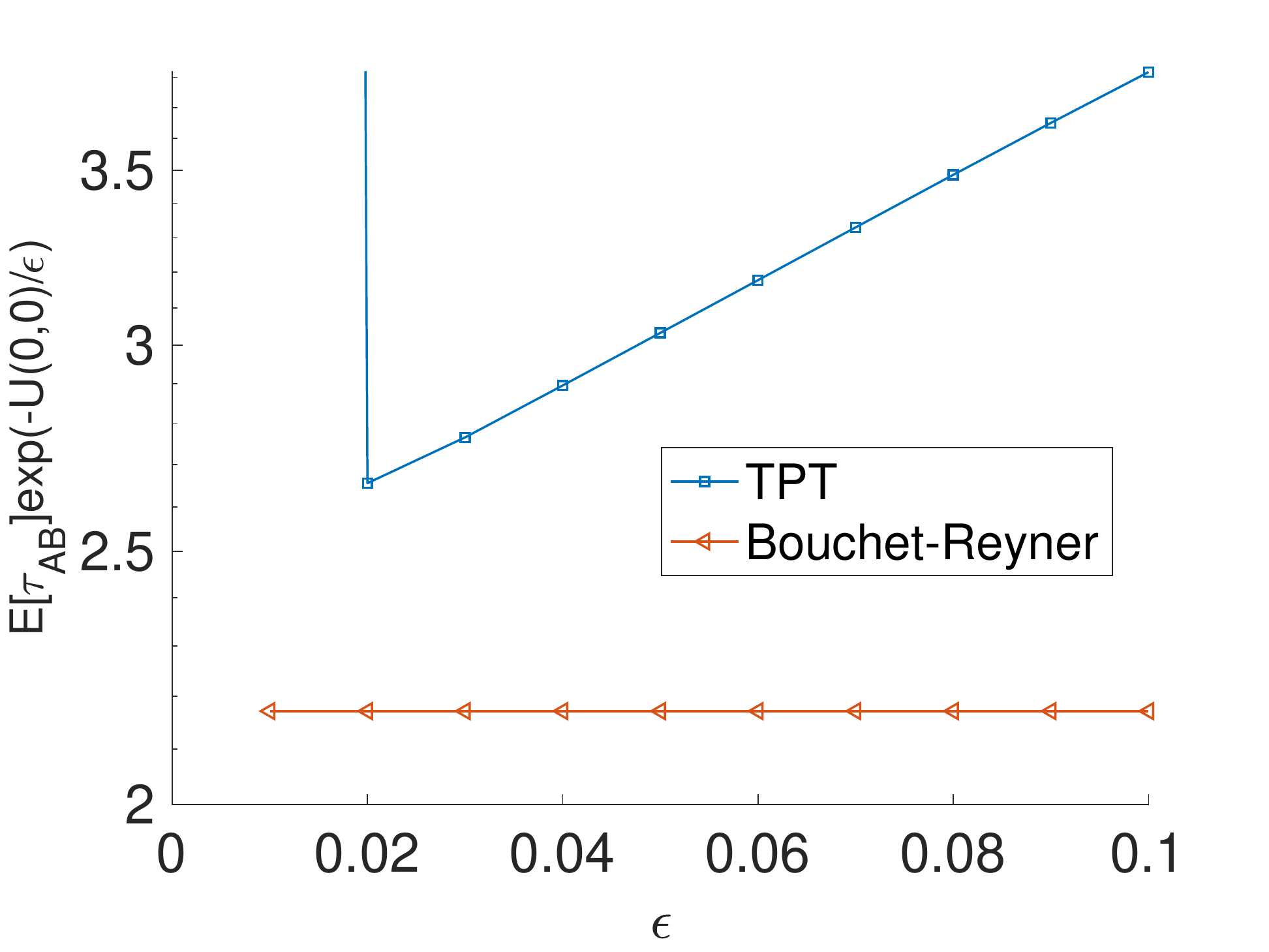}
\caption{(a): The expected escape time $\E \tau^{\O_- \to \O_*} $ from  the basin of $\O_-$ for $\beta=3$ obtained by the Bouchet-Reygner formula and by inverting the TPT rate for $\epsilon = 0.02,0.03,\ldots,0.10$.  (b): The same expected times  $\E \tau^{\O_- \to \O_*} $ at $\beta=3 $ multiplied by the exponent  $\exp( -U(\O_*)/\e)$ in order to expose the prefactor.
}
\label{fig:Etau_a3}
\end{figure}

\begin{figure}
(a)\includegraphics[width=0.4\textwidth]{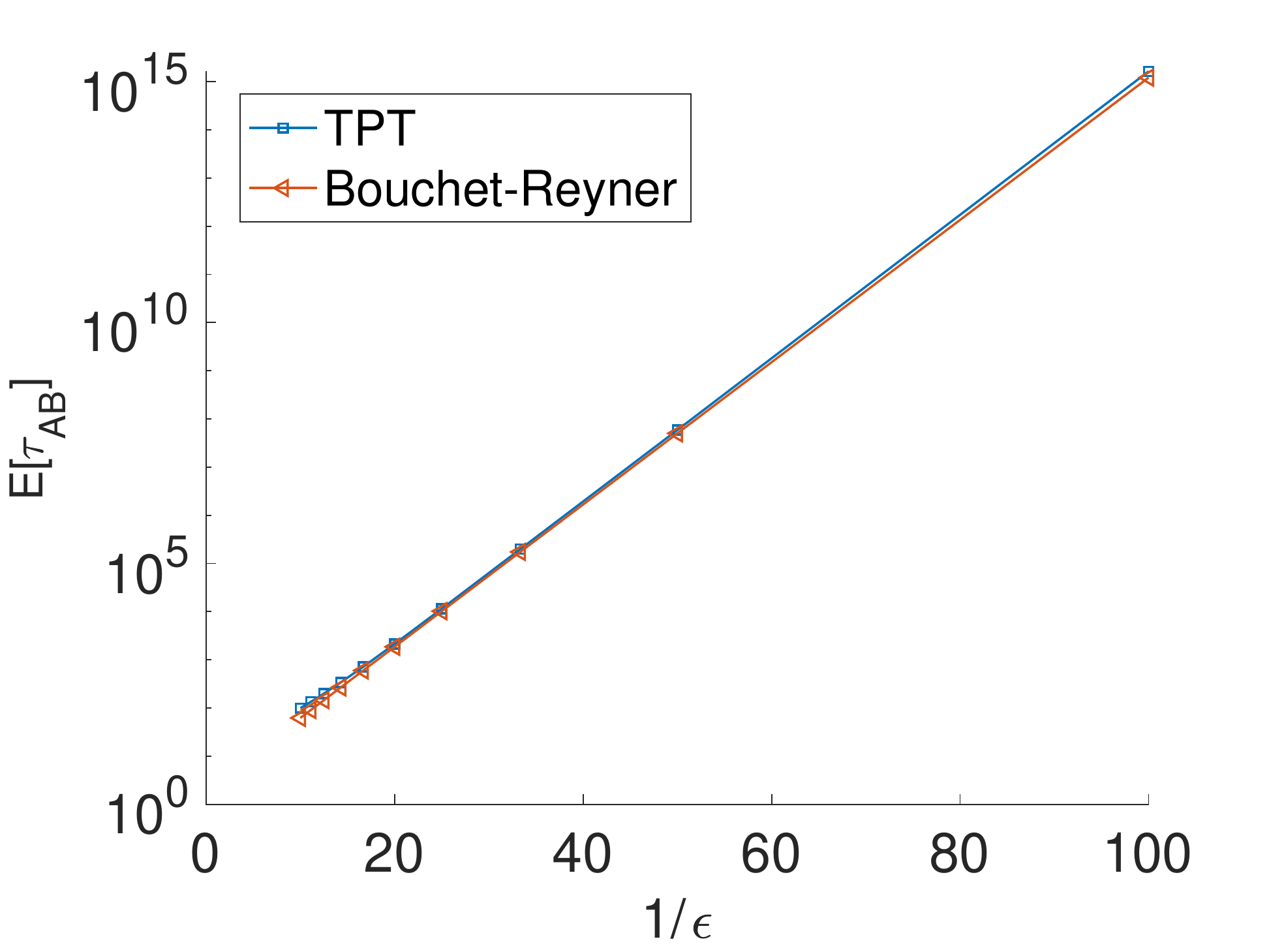}
(b)\includegraphics[width=0.4\textwidth]{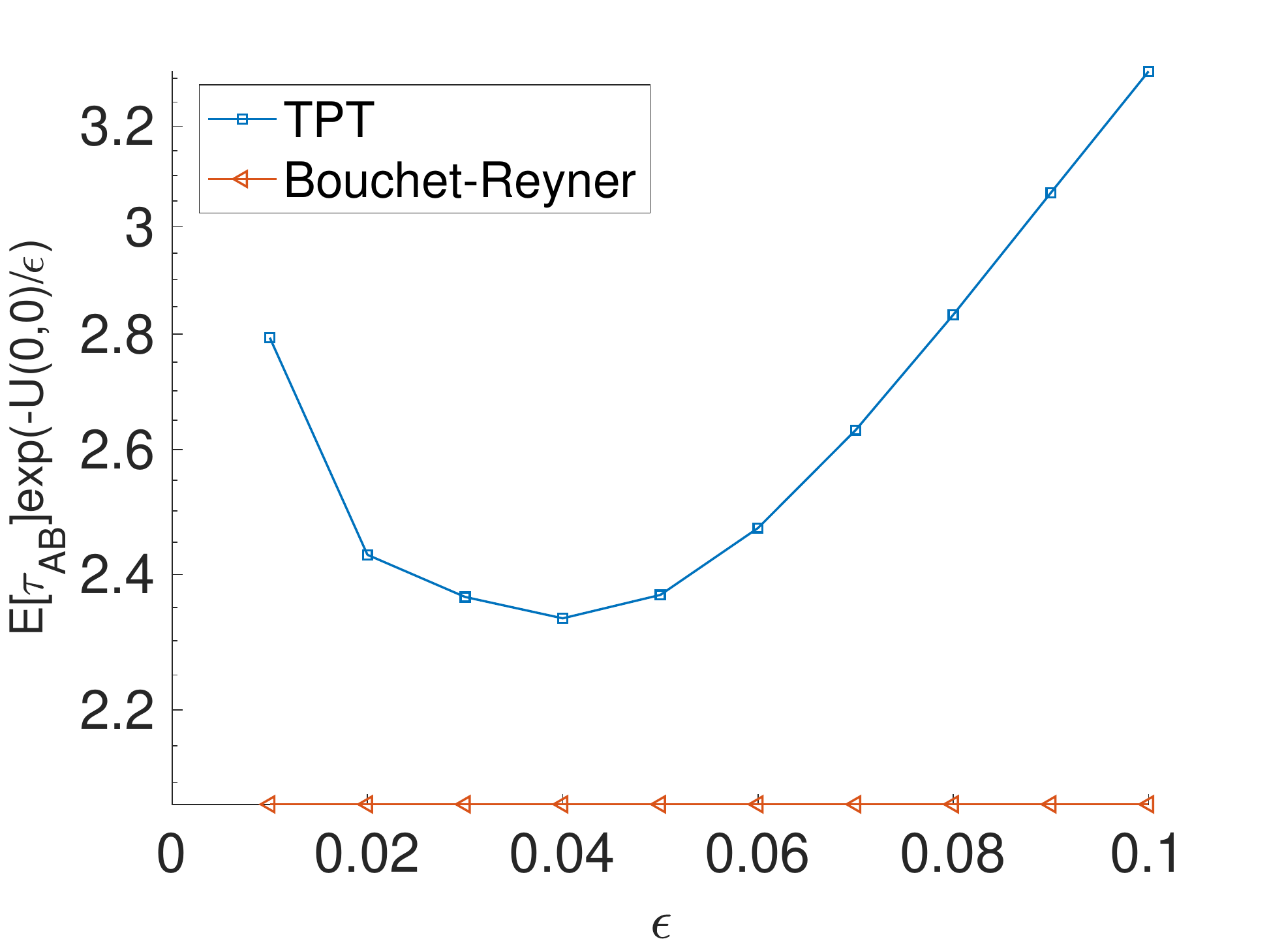}
\caption{(a): The expected escape time $\E \tau^{x_- \to x^*}$ from  the basin of $\O_-$ for $\beta=10$ obtained by the Bouchet-Reygner formula and by inverting the TPT rate for $\epsilon = 0.02,0.03, \ldots,0.10$.  (b): The same expected times  $\E \tau^{x_- \to x^*}$ at $\beta=10$ multiplied by the exponent  $\exp( -U(x^*)/\e)$ in order to expose the prefactor.
}
\label{fig:Etau_a10}
\end{figure}

The TPT reactive current and the MAPs at $\epsilon = 0.04$ for $\beta=3$ and $\beta=10$  are displayed in Figure \ref{fig:Rcurrent}.
The quasipotential with respect to $\O_-$ is smooth for $\beta = 3$ in the left half-plane, hence \eqref{eq:quasi:bouchet} is applicable.
While this is not the case for $\beta=10$, we still can evaluate \eqref{eq:quasi:bouchet} estimating $|\mathrm{det} H(\O_*)|$ by replacing $\O_{*}$ with a nearby mesh point lying in the upper left quadrant. We compare the estimates for  $\E \tau^{\O_- \to \O_*}$ obtained using \eqref{eq:quasi:bouchet} to the inverse TPT rate for $\epsilon=0.01, 0.02,\ldots,0.1$. Details for computing the reactive current and the TPT rate are provided in Appendix B. The results for $\beta=3$ and $\beta =10$ are shown, respectively,  in Figures \ref{fig:Etau_a3} and \ref{fig:Etau_a10}. Note that the implementation of the TPT tools becomes difficult if $\epsilon$ is very small. The plots exhibit an agreement between these expected escape times at best up to $\sim15\%$ for $\beta = 3$ and $\sim12\%$ for $\beta = 10$, which is very good given that both techniques bear certain errors. { Nonetheless, Figures \ref{fig:Etau_a3}(b) and \ref{fig:Etau_a10}(b) reveal an $\epsilon$-dependence of the TPT prefactor which is not captured by the Bouchet-Reygner formula and can be a subject of a future investigation.}


\section{Discussion}
\label{sec:conclusions}
We have developed EJM, the first semi-Lagrangian quasipotential solver exhibiting second order of convergence. Besides the quasipotential, EJM also computes its gradient with nearly second order of convergence. Due to an adaptation of Mirebeau's idea of designing causal stencils, the runtimes of this algorithm are similar to those of OLIM-midpoint. Accurately computed gradients can be differentiated numerically, allowing us to compute the prefactor for the WKB approximation for the invariant measure. We have done it for the Maier-Stein SDE. This is the first reported numerical solution for the prefactor to the best of our knowledge.  

Thus far, we have tested the efficacy of EJM on nonlinear drift fields that have known, \textit{smooth} solutions. In such cases, it has proven \textit{far} more accurate than previous techniques when the drift fields contain at most moderate rotational components. In particular, it displays 2nd order convergence in the error of $U$ and close to 2nd order convergence in the gradient $\nabla U$. { Note that if the gradient were evaluated using finite differences, it would be only first-order accurate.}  Nonetheless, it remains to test our algorithm in more diverse settings.

Admittedly, our solver suffers from an issue of over-complexity in a couple of areas. First, the higher-order updates require minimizations $\eqref{eq:quasi:update_1pt_cubic}$ and $\eqref{eq:quasi:update_triangle_cubic}$ that are performed via Newton's method with first and second-order derivatives computed and implemented by hand. These formulas are extremely lengthy and complicated, which makes EJM a very difficult algorithm to implement from scratch. {Our codes are available on GitHub \cite{EJMpaskal}.} We have yet to explore the efficacy of implementing automatic differentiation techniques.
Second, the Accept/Reject procedure described in Algorithm \ref{alg:quasi:update} and Condition \ref{con:update} are more complicated than in traditional Fast Marching methods, which accept a new proposed update $U_{\mathrm{new}}$ value if it is smaller than any previously proposed values. Instead, EJM requires a more complicated rule. The two key issues that cause this are the (a) unacceptably large errors coming from triangle updates computed over large distances and (b) the presence of \textit{undesired} local minimizers of $\eqref{eq:quasi:update_triangle_cubic}$. The procedure presented in this article is our best attempt to simultaneously circumvent both of these issues, but it is possible there exist more concise solutions.

Numerical analysis of EJM is left for future work. Establishing theoretical convergence guarantees for the quasipotential and its gradient is an interesting and challenging problem in its own right. 

{ The promotion of the EJM to 3D will involve two components:  a targeted search for the fastest characteristic and at least cubically accurate local approximations. The targeted search has been implemented in higher dimensions for the eikonal equation in Riemannian metrics \cite{mirebeau143D,mirebeau19}. An upgrade for 3D of our version of the inverse stencil generator (Algorithm 4) is straightforward. The extension of local approximations presents a bigger challenge. In the 3D version of the OLIM \cite{potter20} simplex updates with linear interpolation in their triangular bases were used. While one can build bicubic interpolants in square faces of cubic mesh cells following the recipe in \cite{nave10}, this will make the solver less flexible as it will require four {\sf Accepted} points to attempt a 3D update. On the other hand, an interpolation based on values of the quasi-potential and its gradient at the vertices of a triangle is only second-order accurate \cite{farin86}, which may lead to deterioration of quadratic accuracy of the quasipotential solver extended to 3D in a manner similar to the OLIM. Another concern is that the extension of the EJM to 3D may exacerbate the issues originating from the imperfection of stencils whose construction is based on the linear approximation for MAP segments. In 2D, we address these issues by introducing a sophisticated Accept/Reject rule and the Fail-Safe procedure.
In sum, the extension of the EJM to 3D is possible but requires resolving some technical issues.}  

{ Extending the EJM for SDEs with nondegenerate geometric noise, i.e., for $dX_t = b(X)t)dt +\sqrt{\epsilon}\sigma(X_t)dW_t$ where $\det(\sigma(X_t)\sigma^\top(X_t)) > \delta > 0$ is straightforward. We chose not to do it in this work as it would somewhat complicate formulas to be evaluated at each step of the algorithm and require significant additional testing. The extension of the OLIM to this case \cite{dahiya18} has shown that the numerical errors are sensitive to the ratio of the largest and the smallest eigenvalues of the diffusion matrix $\sigma(X_t)\sigma^\top(X_t)$.}

In summary, we have achieved our goal to design a second-order accurate and efficient quasipotential solver. It did provide us with a tool for computing the full WKB approximation for the invariant probability measure and the sharp Bouchet-Reygner estimates for the expected escape time. However, the resulting method is complicated and involves a number of safeguards such as the fail-safe procedure and a sophisticated Accept/Reject rule.

\section*{Data availability statement}
The datasets generated during and/or analyzed during the current study are available in the GitHub repository, \href{https://github.com/npaskal/EfficientJetMarcher}{https://github.com/npaskal/EfficientJetMarcher}.

\section*{Declarations}
%

\begin{itemize}
\item Funding: this work was partially supported by  NSF CAREER grant DMS-1554907 and AFOSR MURI grant FA9550-20-1-0397.
\item Conflict of interest/Competing interests: no conflict of interest.
\item Ethics approval: Not Applicable.
\item Consent to participate: Not Applicable.
\item Consent for publication: Not Applicable.
\item Availability of data and materials: the manuscript has associated data available on GitHub \cite{EJMpaskal}.
\item Code availability: Our codes are available on GitHub \cite{EJMpaskal}.
\end{itemize}


 \appendix
\setcounter{equation}{0}
\renewcommand{\theequation}{\Alph{section}-\arabic{equation}}
    \setcounter{lemma}{0}
    \renewcommand{\thelemma}{\Alph{section}\arabic{lemma}}

\section{Appendix: Solving minimization problems arising in update procedures}
\begin{figure}
\includegraphics[width=0.75\textwidth]{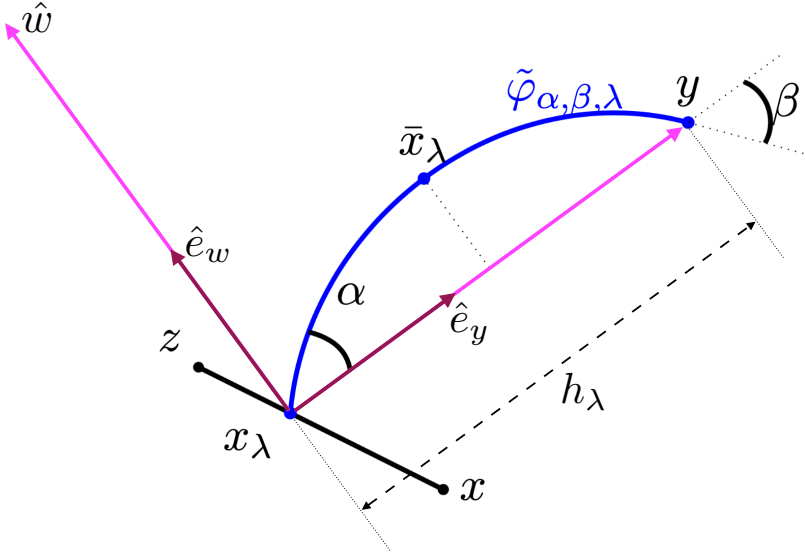}
\caption{An illustration for solving the triangle update optimization problem \eqref{eq:quasi:update_triangle_cubic}.}
\label{fig:A}
\end{figure}
The minimization problems \eqref{eq:quasi:update_1pt_cubic} and \eqref{eq:quasi:update_triangle_cubic} for one-point update and triangle update are solved using Newton's method.
In order to implement it, one needs to calculate the gradient and the Hessian of the objective function.
First, we will redefine the objective function for triangle update in a way more amenable for taking derivatives.
For each update, we introduce a local Euclidean coordinate system with the origin at the point $x_{\lambda}$ and axes directed along 
\begin{equation}
\label{eq:eyw}
\hat{e}_y:=\frac{y-x_{\lambda}}{h_{\lambda}}\quad{\rm and}\quad \hat{e}_w:=R\hat{e}_y,
\end{equation}
where $h_{\lambda}: = \|y-x_{\lambda}\|$ and $R=\left[\begin{array}{rr}0&-1\\1&0\end{array}\right]$ is the 90 degree counter clockwise rotation matrix (Figure \ref{fig:A}).
The cubic curve $\tilde{\varphi}_{\alpha,\beta,\lambda}$ approximating the MAP segment arriving at $y$ from $x_{\lambda}$ 
is represented as a graph of a cubic polynomial in this coordinate system. Let $a_0:=\tan\alpha$ and $a_1:=\tan\beta$.
Then 
\begin{equation}
\label{eq:q}
\tilde{\varphi}_{\alpha,\beta,\lambda}(r) = \left[\begin{array}{c}r\\q(r)\end{array}\right],\quad{\rm where}\quad
q(r) = a_0r-\frac{a_0}{h_{\lambda}}r^2+
\frac{a_0+a_1}{h_{\lambda}^2}
r^2(r-h_{\lambda}),\quad r\in[0,h_{\lambda}].
\end{equation}
Note that the cubic polynomial $q(r)$ satisfies $q(0)=q(h_{\lambda}) = 0$, $q'(0) = a_0$, $q'(h_{\lambda})=a_1$, and the cubic curve $\tilde{\varphi}_{\alpha,\beta,\lambda}$ 
is parametrized by the arclength of its chord.
Then the velocity vector for $\tilde{\varphi}_{\alpha,\beta,\lambda}$ is
\begin{equation}
\label{eq:qprime}
\dot{\tilde{\varphi}}_{\alpha,\beta,\lambda}(r) = \left[\begin{array}{c}1\\q'(r)\end{array}\right],\quad{\rm where}\quad
q'(r) = a_0-\frac{2a_0}{h_{\lambda}}r+
\frac{a_0+a_1}{h_{\lambda}^2}
(3r^2-2h_{\lambda}r).
\end{equation}
Therefore, the geometric action along $\tilde{\varphi}_{\alpha,\beta,\lambda}$ is
\begin{equation}
\label{eq:ga1}
\int_0^{h_{\lambda}}\|b(\tilde{\varphi}_{\alpha,\beta,\lambda})\|\|\dot{\tilde{\varphi}}_{\alpha,\beta,\lambda}\| - 
\langle b(\tilde{\varphi}_{\alpha,\beta,\lambda}),\dot{\tilde{\varphi}}_{\alpha,\beta,\lambda}\rangle dr,
\end{equation}
where $\langle\cdot,\cdot\rangle$ denotes the dot product, and the integrand can be detailed as
\begin{equation}
\label{eq:ga2}
\|b(x_{\lambda} + r\hat{e}_y+q(r) \hat{e}_w)\|\sqrt{1 +[q'(r)]^2} - \langle b(x_{\lambda} + r\hat{e}_y+q(r) \hat{e}_w),\hat{e}_y + q'(r) \hat{e}_w\rangle.
\end{equation}
To apply Simpson's rule, we need to evaluate \eqref{eq:ga2} at the endpoints $r=0$ and $r=h_{\lambda}$, and at the midpoint
$$
\bar{x}_{\lambda} = x_{\lambda} + \frac{h_{\lambda}}{2}\hat{e}_y+q\left(\frac{h_{\lambda}}{2}\right) \hat{e}_w 
= \frac{y+x_{\lambda}}{2} + \frac{h_{\lambda}}{8}(a_0-a_1).
$$
We show how one can write the dot product in \eqref{eq:ga2} in a compact manner:
\begin{align*}
\langle b(\tilde{\varphi}_{\alpha,\beta,\lambda}(r) ),\hat{e}_y + q'(r) \hat{e}_w\rangle &= 
\langle b(\tilde{\varphi}_{\alpha,\beta,\lambda}(r) ),\hat{e}_y + q'(r) R\hat{e}_y\rangle\\
=\langle b(\tilde{\varphi}_{\alpha,\beta,\lambda}(r) ),(I + Rq'(r))\hat{e}_y\rangle &= 
\langle\hat{e}_y,(I + Rq'(r))^\top b(\tilde{\varphi}_{\alpha,\beta,\lambda}(r) )\rangle.
\end{align*}
We denote the matrix $I + Rq'(r)$ by $R_{q'(r)}$:
$R_{q'(r)} = \left[\begin{array}{cc}1&-q'(r) \\ q'(r) & 1\end{array}\right]$.
Observing that $R_{q'(r)}^\top\equiv R_{-q'(r)}$ and $q'(\sfrac{h_\lambda}{2}) = -\tfrac{1}{4}(a_0+a_1)$, we write out the objective function for triangle update:
\begin{align}
G(a_0,a_1,\lambda) &=p(\lambda) +\frac{h_{\lambda}}{6}\left(\|b(x_{\lambda})\| \sqrt{1+a_0^2} 
+ 4\|b(\bar{x}_{\lambda})\|\sqrt{1+\frac{(a_0+a_1)^2}{16}}
+\|b(y)\| \sqrt{1+a_1^2}\right)\notag\\
&-\frac{1}{6}\left\langle y-x_{\lambda},\left(R_{-a_0}b(x_{\lambda}) + 4R_{\frac{a_0+a_1}{4}}b(\bar{x}_{\lambda}) 
+R_{-a_1}b(y)\right)\right\rangle. \label{eq:G}
\end{align}
Here, the polynomial $p(\lambda)$ is the Hermite interpolant for $U$ along the segment $[x,z]$ given by
\begin{equation}
\label{eq:plam}
p(\lambda) = u_0p_0(\lambda) + u_1p_0(1-\lambda)+u_0'p_1(\lambda) - u_1'p_1(1-\lambda),
\end{equation}
where $u_0=U(x)$, $u_1=U(z)$, $u_0' = (z-x)\cdot\nabla U(x)$, $u_0' = (z-x)\cdot\nabla U(z)$, $p_0(\lambda) = 1-3\lambda^2+2\lambda^3$, and $p_1(\lambda) = \lambda(1-\lambda)^2$.

The objective function for one-point update is obtained in a similar manner:
\begin{align}
G_1(a_0,a_1) &=U(x) +\frac{h}{6}\left(\|b(x)\| \sqrt{1+a_0^2} + 4\|b(\bar{x})\|\sqrt{1+\frac{(a_0+a_1)^2}{16}}
+\|b(y)\| \sqrt{1+a_1^2}\right)\notag\\
&-\frac{1}{6}\left\langle y-x,\left(R_{-a_0}b(x) + 4R_{\frac{a_0+a_1}{4}}b(\bar{x}) +R_{-a_1}b(y)\right)\right\rangle. \label{eq:G1}
\end{align}
We conduct minimization with respect to the slopes $a_0$ and $a_1$ rather than with respect to the angles $\alpha$ and $\beta$.
Now, as we have the objective functions \eqref{eq:G} and \eqref{eq:G1}, it is straightforward to calculate their first and second derivatives.

Below we will write out the derivatives for $G(a_0,a_1,\lambda)$. The derivatives for $G_1(a_0,a_1)$ are readily obtained from them.
We use the following shorthands:
\begin{align*}
& \delta  := z-x,
\\ & q_{\lambda} := \frac{d\bar{x}_{\lambda}}{d\lambda} = \frac{1}{2}( \delta  - \frac{a_0-a_1}{4} R \delta ),
\\ & h_\lambda' = \frac{d h_{\lambda}}{d\lambda} = - \frac{\langle{y}-x_{\lambda},\delta  \rangle}{h_{\lambda}},
\\ & h_\lambda'' = \frac{d^2 h_{\lambda}}{d\lambda^2} = \frac{|\delta |^2}{h_{\lambda}} - \frac{|\langle {y} - x_{\lambda}, \delta  \rangle|^2}{h_{\lambda}^3},
\\ & \hat{w}_{\lambda} := R (y - x_\lambda)
\end{align*}
Also, we will use the following notation for the derivatives of $b$:
$Db$ is the Jacobian matrix for $b$: $(Db)_{ij} = \tfrac{\partial b_i}{\partial x_j}$, $i,j = 1,2$, and 
$D^2b$ is the second derivative tensor: $(D^2b)_{ijk} = \tfrac{\partial^2 b_i}{\partial x_j\partial x_k}$, and $(u D^2 b v)_i = u_j(D^2b)_{ijk}v_k$.

{\bf Gradient:}
\begin{align*}
	\frac{\partial G}{\partial a_0} & = \frac{h_\lambda}{6} \Big[  \frac{a_0}{\sqrt{1+a_0^2}}|b(x_{\lambda})| +
	 \frac{1}{4} \frac{a_0+a_1}{\sqrt{1+\frac{(a_0+a_1)^2}{16}}}|b(\bar{x}_{\lambda})| \\
	& + \frac{h_{\lambda}}{2}  \sqrt{1+\frac{(a_0+a_1)^2}{16}} \frac{\langle b(\bar{x}_{\lambda}),Db(\bar{x}_{\lambda}) \hat{w}_{\lambda} \rangle}{|b(\bar{x}_{\lambda})|}   
	\\ & +  \langle \hat{w}_{\lambda},b(\bar{x}_{\lambda})-b(x_{\lambda}) \rangle
	- \frac{1}{2} \langle R_{(a_0+a_1)/4} Db(\bar{x}_{\lambda})\hat{w}_{\lambda},{y}-x_{\lambda} \rangle \Big];
\end{align*}
\begin{align*}
	\frac{\partial G}{\partial a_1} & = \frac{h_\lambda}{6} \Big[  \frac{a_1}{\sqrt{1+a_1^2}}|b({y})|   + \frac{1}{4} \frac{a_0+a_1}{\sqrt{1+\frac{(a_0+a_1)^2}{16}}}|b(\bar{x}_{\lambda})| \\
	& - \frac{h_{\lambda}}{2}  \sqrt{1+\frac{(a_0+a_1)^2}{16}} \frac{\langle b(\bar{x}_{\lambda}),Db(\bar{x}_{\lambda}) \hat{w}_{\lambda} \rangle}{|b(\bar{x}_{\lambda})|} 
	\\ & + \langle \hat{w}_{\lambda},b(\bar{x}_{\lambda})-b({y}) \rangle + \frac{1}{2} \langle R_{(a_0+a_1)/4} Db(\bar{x}_{\lambda})\hat{w}_{\lambda},{y}-x_{\lambda} \rangle \Big];
\end{align*}
\begin{align*}
	\frac{\partial G}{\partial \lambda} & = p'(\lambda) + \frac{h_\lambda'}{6} \Big[ \sqrt{1+a_0^2} |b(x_{\lambda})| 
	+4 \sqrt{1+\frac{(a_0+a_1)^2}{16}} |b(\bar{x}_{\lambda})|  + \sqrt{1+a_1^2} |b({y})|\Big] 
	\\ & + \frac{h_\lambda}{6} \Big[\sqrt{1+a_0^2} \frac{\langle b(x_{\lambda}),Db(x_{\lambda})\delta  \rangle}{|b(x_{\lambda})|} + 4 \sqrt{1+\frac{(a_0+a_1)^2}{16}} \frac{\langle b(\bar{x}_{\lambda}),Db(\bar{x}_{\lambda}) q_{\lambda} \rangle}{|b(\bar{x}_{\lambda})|}\Big]
	 \\ & + \frac{1}{6} \langle \delta , R_{-a_0}b(x_{\lambda}) + 4R_{(a_0+a_1)/4} b(\bar{x}_{\lambda}) + R_{-a_1}b({y}) \rangle
	 \\ & - \frac{1}{6} \langle {y} - x_{\lambda}, R_{-a_0} Db(x_{\lambda})\delta  + 4 R_{(a_0+a_1)/4} Db(\bar{x}_{\lambda}) q_{\lambda} \rangle,
\end{align*}

{\bf Hessian:}
\begin{align*}
	\frac{\partial ^2G}{\partial a_0 \partial a_1} & = - \frac{h_\lambda^2}{16}\sqrt{1+\frac{(a_0+a_1)^2}{16}} \frac{1}{|b(\bar{x}_{\lambda})|} \Big[ |Db(\bar{x}_{\lambda}) \hat{w}_{\lambda}|^2 + \langle b(\bar{x}_{\lambda}),\hat{w}_{\lambda} D^2b(\bar{x}_{\lambda})\hat{w}_{\lambda} \rangle  \\
	&- \frac{|\langle b(\bar{x}_{\lambda}),Db(\bar{x}_{\lambda})\hat{w}_{\lambda} \rangle|^2}{|b(\bar{x}_{\lambda})|^2} \Big]
	\\ & + \frac{1}{4}|b(\bar{x}_{\lambda})|  \Big(1 + \frac{(a_0+a_1)^2}{16}\Big)^{-3/2} + \frac{h_\lambda}{16}\langle {y} - x_{\lambda}, R_{(a_0+a_1)/4} \hat{w}_{\lambda} D^2b(\bar{x}_{\lambda})\hat{w}_{\lambda} \rangle,
\end{align*}
\begin{align*}
	\frac{\partial^2 G}{\partial a_0 \partial \lambda} & = \frac{h_{\lambda}'}{6} \Big[ \frac{a_0}{\sqrt{1+a_0^2}} |b(x_{\lambda})| + 
	\frac{1}{4}\frac{a_0+a_1}{\sqrt{1+\frac{(a_0+a_1)^2}{16}}}|b(\bar{x}_{\lambda})| \\ 
	&+ \frac{h_{\lambda}}{2} \sqrt{1+\frac{(a_0+a_1)^2}{16}}  \frac{\langle b(\bar{x}_{\lambda}),Db(\bar{x}_{\lambda})\hat{w}_{\lambda} \rangle}{|b(\bar{x}_{\lambda})|} \Big]
	\\ & + \frac{h_{\lambda}}{6} \Big[ \frac{a_0}{\sqrt{1+a_0^2}} \frac{\langle b(x_{\lambda}),Db(x_{\lambda})\delta  \rangle}{|b(x_{\lambda})|} + \frac{1}{4}\frac{a_0+a_1}{\sqrt{1+\frac{(a_0+a_1)^2}{16}}} \frac{\langle b(\bar{x}_{\lambda}),Db(\bar{x}_{\lambda}) q_{\lambda} \rangle}{|b(\bar{x}_{\lambda})|}  \Big]
	\\ & + \frac{h_{\lambda}^2}{12}\sqrt{1+\frac{(a_0+a_1)^2}{16}} \Big[ \frac{\langle Db(\bar{x}_{\lambda})\hat{w}_{\lambda},Db(\bar{x}_{\lambda})q_{\lambda}  + \langle b(\bar{x}_{\lambda}), \hat{w}_{\lambda} D^2b(\bar{x}_{\lambda}) q_{\lambda} \rangle}{|b(\bar{x}_{\lambda})|}  
	\\ & - \frac{\langle b(\bar{x}_{\lambda}),Db(\bar{x}_{\lambda}) R \delta  \rangle}{h_{\lambda} |b(\bar{x}_{\lambda})|} 
  - \frac{\langle b(\bar{x}_{\lambda}),Db(\bar{x}_{\lambda})q_{\lambda} \rangle \langle b(\bar{x}_{\lambda},Db(\bar{x}_{\lambda})\hat{w}_{\lambda} \rangle}{|b(\bar{x}_{\lambda})|^3}  \Big]
	\\ &+ \frac{1}{6} \langle \delta , -R b(x_{\lambda}) +R b(\bar{x}_{\lambda}) + \frac{h_{\lambda}}{2} R_{(a_0+a_1)/4} Db(\bar{x}_{\lambda})\hat{w}_{\lambda} \rangle 
	\\ & -\frac{1}{6} \langle {y} - x_{\lambda}, -RDb(x_{\lambda}) \delta  + R Db(\bar{x}_{\lambda})q_{\lambda} + \frac{h_{\lambda}}{2} R_{(a_0+a_1)/4}\hat{w}_{\lambda} D^2b(\bar{x}_{\lambda})q_{\lambda}\\
	& -\frac{1}{2}R_{(a_0+a_1)/4} Db(\bar{x}_{\lambda})R\delta  \rangle,
\end{align*}
\begin{align*}
	\frac{\partial^2 G}{\partial a_1 \partial \lambda} & = \frac{h_{\lambda}'}{6} \Big[  \frac{1}{4}\frac{a_0+a_1}{\sqrt{1+\frac{(a_0+a_1)^2}{16}}}|b(\bar{x}_{\lambda})| - \frac{h_{\lambda}}{2} \sqrt{1+\frac{(a_0+a_1)^2}{16}}  \frac{\langle b(\bar{x}_{\lambda}),Db(\bar{x}_{\lambda})\hat{w}_{\lambda} \rangle}{|b(\bar{x}_{\lambda})|} \\
	&+ \frac{a_1}{\sqrt{1+a_1^2}} |b({y})|\Big]
	 + \frac{h_{\lambda}}{6} \Big[  \frac{1}{4}\frac{a_0+a_1}{\sqrt{1+\frac{(a_0+a_1)^2}{16}}} \frac{\langle b(\bar{x}_{\lambda}),Db(\bar{x}_{\lambda}) q_{\lambda} \rangle}{|b(\bar{x}_{\lambda})|}  \Big]
	\\ & + \frac{h_{\lambda}^2}{12}\sqrt{1+\frac{(a_0+a_1)^2}{16}} \Big[ -\frac{\langle Db(\bar{x}_{\lambda})\hat{w}_{\lambda},Db(\bar{x}_{\lambda})q_{\lambda}  + \langle b(\bar{x}_{\lambda}), \hat{w}_{\lambda} D^2b(\bar{x}_{\lambda}) q_{\lambda} \rangle}{|b(\bar{x}_{\lambda})|}  
	\\ & + \frac{\langle b(\bar{x}_{\lambda}),Db(\bar{x}_{\lambda}) R \delta  \rangle}{h_{\lambda} |b(\bar{x}_{\lambda})|} 
  + \frac{\langle b(\bar{x}_{\lambda}),Db(\bar{x}_{\lambda})q_{\lambda} \rangle \langle b(\bar{x}_{\lambda},Db(\bar{x}_{\lambda})\hat{w}_{\lambda} \rangle}{|b(\bar{x}_{\lambda})|^3}  \Big]
	\\ &+ \frac{1}{6} \langle \delta , -R b({y}) +R b(\bar{x}_{\lambda}) - \frac{h_{\lambda}}{2} R_{(a_0+a_1)/4} Db(\bar{x}_{\lambda})\hat{w}_{\lambda} \rangle 
	\\ & -\frac{1}{6} \langle {y} - x_{\lambda}, + R Db(\bar{x}_{\lambda})q_{\lambda} - \frac{h_{\lambda}}{2} R_{(a_0+a_1)/4}\hat{w}_{\lambda} D^2b(\bar{x}_{\lambda})q_{\lambda} +\frac{1}{2}R_{(a_0+a_1)/4} Db(\bar{x}_{\lambda})R\delta  \rangle,
\end{align*}
\begin{align*}
\frac{\partial^2 G }{\partial a_0^2} 
& = \frac{h_{\lambda}}{6} \Big[ (1+a_0^2)^{-3/2}|b(x_{\lambda})| + \frac{1}{4} \Big(1+\frac{(a_0+a_1)^2}{16} \Big)^{-3/2}|b(\bar{x}_{\lambda})| \Big]
\\ & + \frac{h_{\lambda}^2}{96}  \frac{a_0+a_1}{\sqrt{1+\frac{(a_0+a_1)^2}{16}}} \frac{\langle b(\bar{x}_{\lambda}),Db(\bar{x}_{\lambda}) \hat{w}_{\lambda} \rangle}{|b(\bar{x}_{\lambda})|}
\\ & + \frac{h_{\lambda}^3}{96}\sqrt{1+\frac{(a_0+a_1)^2}{16}} \Big[ \frac{ |Db(\bar{x}_{\lambda})\hat{w}_{\lambda}|^2 + \langle b(\bar{x}_{\lambda}),\hat{w}_{\lambda} D^2b(\bar{x}_{\lambda})\hat{w}_{\lambda} \rangle}{|b(\bar{x}_{\lambda})|} - \frac{\langle b(\bar{x}_{\lambda}),Db(\bar{x}_{\lambda})\hat{w}_{\lambda}\rangle^2}{|b(\bar{x}_{\lambda})|^3}
\Big]
\\ & +\frac{h_{\lambda}^2}{24} \langle \hat{w}_{\lambda}, Db(\bar{x}_{\lambda})\hat{w}_{\lambda}\rangle - \frac{h_{\lambda}^2}{96}\langle R_{(a_0+a_1)/4}\hat{w}_{\lambda} Db(\bar{x}_{\lambda})\hat{w}_{\lambda},{y} - x_{\lambda} \rangle,
\end{align*}
\begin{align*}
\frac{\partial^2 G }{\partial a_1^2} 
& = \frac{h_{\lambda}}{6} \Big[ (1+a_1^2)^{-3/2}|b({y})| + \frac{1}{4} \Big(1+\frac{(a_0+a_1)^2}{16} \Big)^{-3/2}|b(\bar{x}_{\lambda})| \Big]
\\ & - \frac{h_{\lambda}^2}{96}  \frac{a_0+a_1}{\sqrt{1+\frac{(a_0+a_1)^2}{16}}} \frac{\langle b(\bar{x}_{\lambda}),Db(\bar{x}_{\lambda}) \hat{w}_{\lambda} \rangle}{|b(\bar{x}_{\lambda})|}
\\ & + \frac{h_{\lambda}^3}{96}\sqrt{1+\frac{(a_0+a_1)^2}{16}} \Big[ \frac{ |Db(\bar{x}_{\lambda})\hat{w}_{\lambda}|^2 + \langle b(\bar{x}_{\lambda}),\hat{w}_{\lambda} D^2b(\bar{x}_{\lambda})\hat{w}_{\lambda} \rangle}{|b(\bar{x}_{\lambda})|} - \frac{\langle b(\bar{x}_{\lambda}),Db(\bar{x}_{\lambda})\hat{w}_{\lambda}\rangle^2}{|b(\bar{x}_{\lambda})|^3}
\Big]
\\ & -\frac{h_{\lambda}^2}{24} \langle \hat{w}_{\lambda}, Db(\bar{x}_{\lambda})\hat{w}_{\lambda}\rangle - \frac{h_{\lambda}^2}{96}\langle R_{(a_0+a_1)/4}\hat{w}_{\lambda} Db(\bar{x}_{\lambda})\hat{w}_{\lambda},{y} - x_{\lambda} \rangle,
\end{align*}
\begin{align*}
\frac{\partial^2G }{\partial \lambda^2} 
& =p''(\lambda) + \frac{h_\lambda''}{6} \Big[ \sqrt{1+a_0^2} |b(x_{\lambda})| 
	+4 \sqrt{1+\frac{(a_0+a_1)^2}{16}} |b(\bar{x}_{\lambda})|  + \sqrt{1+a_1^2} |b({y})|\Big] 
	\\ & + \frac{h_{\lambda}'}{3} \Big[\sqrt{1+a_0^2} \frac{\langle b(x_{\lambda}),Db(x_{\lambda})\delta  \rangle}{|b(x_{\lambda})|} + 4 \sqrt{1+\frac{(a_0+a_1)^2}{16}} \frac{\langle b(\bar{x}_{\lambda}),Db(\bar{x}_{\lambda}) q_{\lambda} \rangle}{|b(\bar{x}_{\lambda})|} \Big]
	\\ & + \frac{h_{\lambda}}{6}  \sqrt{1+a_0^2} \Big[ \frac{ |Db(x_{\lambda})\delta |^2 + \langle b(x_{\lambda}), \delta  D^2b(x_{\lambda})\delta  \rangle}{|b(x_{\lambda})|}  - \frac{ \langle b(x_{\lambda}),Db(x_{\lambda})\delta  \rangle^2}{|b(x_{\lambda})|^3}
	\Big]
	\\ & + \frac{2h_{\lambda}}{3} \sqrt{1+\frac{(a_0+a_1)^2}{16}} \Big[ \frac{ |Db(\bar{x}_{\lambda})q_{\lambda}|^2 + \langle b(\bar{x}_{\lambda}),q_{\lambda} Db(\bar{x}_{\lambda}) q_{\lambda} \rangle}{|b(\bar{x}_{\lambda})|} - \frac{\langle b(\bar{x}_{\lambda}),Db(\bar{x}_{\lambda})q_{\lambda}\rangle^2}{|b(\bar{x}_{\lambda})|^3}    \Big]
	\\ & + \frac{1}{3} \langle \delta , R_{-a_0} Db(x_{\lambda})\delta  + 4R_{(a_0+a_1)/4} Db(\bar{x}_{\lambda})q_{\lambda}  \rangle
	\\ & -\frac{1}{6} \langle {y} - x_{\lambda},R_{-a_0} \delta  D^2b(x_{\lambda}) \delta  + 4 R_{(a_0+a_1)/4} q_{\lambda} D^2b(\bar{x}_{\lambda})q_{\lambda} \rangle.
\end{align*}


\section{Appendix: estimating transition rate using the transition path theory}
The transition path theory \cite{eve2006} is a mathematical framework for quantifying transitions between metastable states in stochastic systems. The generator of the stochastic process governed by SDE \eqref{eq_SDE} is given by
\begin{equation}
    \label{eq:gen}
    \mathcal{L}: = b(x)\cdot\nabla + \frac{\e}{2} \Delta.
\end{equation}
After discretizing the generator to an $N\times N$ mesh, it can be viewed as a generator matrix of a Markov jump process on this mesh, where jumps are enabled only between nearest neighbors. Denoting the discretized generator by $L$, we easily obtain its adjoint $L^*$ by transposing it. The eigenvector of $L^*$ corresponding to its zero eigenvalue is normalized so that the sum of its entries is one and it provides the invariant probability distribution $\mu$ on the mesh. The computational domain for finding $\mu$ was set to $D = [-2,0]\times[-1,1]$, and $N=1025$. The homogeneous Neumann boundary conditions were used at $\partial D$ which correspond to reflecting boundary condition. The symmetry of the Maier-Stein drift field with respect to the $y$-axis and the extremely low probability of reaching the other parts of $\partial D$ justify this choice.

The transition process between given regions $A$ and $B$ is described by the vector field called the \emph{reactive current}. The sets $A$ and $B$ are usually picked as neighborhoods of attractors of $\dot{x}=b(x)$. For the Maier-Stein SDE, We chose $A$ and $B$ to be balls of radius $0.3$ surrounding $\O_-$ and $\O_+$.  
The reactive current for SDE \eqref{eq_SDE} is given by \cite{eve2010}
\begin{equation}
    \label{eq:Rcurrent}
    J_R = q_+q_-J + \frac{\mu\epsilon}{2}\left(q_-\nabla q_+ - q_+\nabla q_-\right),
\end{equation}
where $J = \mu b-\tfrac{\e}{2}\nabla \mu$ is the probability current, and $q_+$ and $q_-$ are the forward and backward committor functions. The forward committor $q_+(x)$ is the probability that the process starting at $x$ first hits $B$ rather than $A$. The backward committor is the probability that the process arriving at $x$ last exited from $A$ rather than from $B$. The forward and backward committors are the solutions to the following boundary-value problems:
\begin{equation}
    \label{eq:comm}
    \begin{cases}\mathcal{L}q_+=0,&x\in\R^n\backslash(A\cup B),\\
    q_+ = 0, & x\in\partial A,\\
    q_+ = 1,& x\in\partial B;\end{cases}\qquad
   \begin{cases}\mathcal{L^\dagger}q_-=0,&x\in\R^n\backslash(A\cup B),\\
    q_- = 1, & x\in\partial A,\\
    q_- = 0,& x\in\partial B.\end{cases}
\end{equation}
The operator $\mathcal{L}^\dagger$ is the generator for the time-reversed process. It is easier to find its discrete counterpart using the formula: $L^\dagger = M^{-1}L^\top M$, where $M = {\sf diag}\{\mu\}$ is the diagonal matrix with the invariant distribution $\mu$ along its diagonal \cite{metzner2009}. To compute the reactive current, we used the computational domain $[-1.5,1.5]\times[-0.75,0.75]\backslash (A\cup B)$ with homogeneous Neumann boundary conditions at the outer boundary. The box $[-1.5,1.5]\times[-0.75,0.75]$ was discretized to $N\times N$ mesh with $N=1025$.

The reaction rate, i.e., the mean number of transitions from $A$ to $B$ per unit time is equal to the average flux of the reactive current through any surface $S$ dividing $A$ and $B$:
\begin{equation}
 \label{eq:TPTrate}
 \nu_{AB} = \int_S J_R\cdot\hat{n}d\sigma,
\end{equation}
where $\hat{n}$ is the unit normal to the surface $S$ pointing towards $B$ and $d\sigma$ is the surface element. For the Maier-Stein SDE, we chose the $y$-axis as the dividing surface.

\section*{Acknowledgements}
This work was partially supported by NSF CAREER Grant DMS-1554907 (MC) and AFOSR MURI grant FA9550-20-1-0397 (MC). We are grateful to the anonymous reviewer for helping us to improve this manuscript.

%
%

\bibliographystyle{spmpsci}      

\begin{thebibliography}{}
%
%
\bibitem{aplot2005}
Beri, S., Mannella, R., Luchinsky, D.G., Silchenko, A.N., McClintock, P.V.E.:
  Solution of the boundary value problem for optimal escape in continuous
  stochastic systems and maps.
Physical Review E \textbf{72}, 036131 (2005)

\bibitem{bouchet15}
Bouchet, F., Reygner, J.: Generalisation of the Eyring-Kramers transition rate
  formula to irreversible diffusion processes.
Annales Henri Poincar{\'e} \textbf{17}, 3499--3532 (2015)

\bibitem{cameron12}
Cameron, M.K.: Finding the quasipotential for nongradient SDEs.
Physica D: Nonlinear Phenomena \textbf{241}, 1532--1550 (2012)

\bibitem{cameron17}
Cameron, M.K.: Construction of the quasi-potential for linear SDEs using false
  quasi-potentials and a geometric recursion (2017).
ArXiv:1801.00327

\bibitem{lorenz63}
Cameron, M.K., Yang, S.: Computing the quasipotential for highly dissipative
  and chaotic sdes. an application to stochastic Lorenz'63.
Communications in Applied Mathematics and Computational Science
  (CAMCoS) \textbf{14}(2), 207--246 (2019)

\bibitem{dahiya18}
Dahiya, D., Cameron, M.K.: Finding the quasipotential for nongradient SDEs.
Physica D: Nonlinear Phenomena \textbf{382--383}, 33--45 (2018)

\bibitem{cameron18}
Dahiya, D., Cameron, M.K.: Ordered line integral methods for computing the
  quasipotential.
 Journal of Scientific Computing \textbf{75}(3), 1351--1384 (2018)

\bibitem{dijkstra59}
Dijkstra, E.W.: A note on two problems in connexion with graphs.
Numerische Mathematik \textbf{1}, 269--271 (1959)

\bibitem{e04}
E, W., Ren, W., Vanden-Eijnden, E.: Minimum action method for the study of rare
  events.
Communications on Pure and Applied Mathematics \textbf{57}(5),
  637--656 (2004)

\bibitem{eve2006}
E, W., Vanden-Eijnden, E.: Towards a theory of transition paths.
Journal of Statistical Physics \textbf{123}(3), 503--523 (2006)

\bibitem{eve2010}
E, W., Vanden-Eijnden, E.: Transition-path theory and path-finding algorithms
  for the study of rare events.
Annual review of physical chemistry \textbf{61}, 391--420 (2010)

\bibitem{farin86}
Farin, G.E.: Triangular bernstein-b{\'e}zier patches.
Computer Aided Geometric Design \textbf{3}, 83--127 (1986)

\bibitem{freidlin12}
Freidlin, M.I., Wentzell, A.D.: Random Perturbations of Dynamical Systems,
  third edition.
 Springer-Verlag (2012)

\bibitem{heymannCPAM}
Heymann, M., Vanden-Eijnden, E.: The geometric minimum action method: A least
  action principle on the space of curves.
 Communications on Pure and Applied Mathematics \textbf{61},
  1052--1117 (2008)

\bibitem{heymann08}
Heymann, M., Vanden-Eijnden, E.: Pathways of maximum likelihood for rare events
  in nonequilibrium systems: Application to nucleation in the presence of
  shear.
 Physical Review Letters \textbf{100}, 140601 (2008)

\bibitem{kikuchi20}
Kikuchi, L., Singh, R., Cates, M.E., Adhikari, R.: Ritz method for transition
  paths and quasipotentials of rare diffusive events.
Phys. Rev. Research \textbf{2}, 033208 (2020)

\bibitem{lin19}
Lin, L., Yu, H., Zhou, X.: Quasi-potential calculation and minimum action
  method for limit cycle.
Journal of Nonlinear Science \textbf{29}(3), 961--991 (2019)

\bibitem{MS1993}
Maier, R.S., Stein, D.L.: Escape problem for irreversible systems.
Physical review. E, Statistical physics, plasmas, fluids, and related
  interdisciplinary topics \textbf{48}(2), 931--938 (1993)

\bibitem{maier96}
Maier, R.S., Stein, D.L.: A scaling theory of bifurcations in the symmetric
  weak-noise escape problem.
Journal of Statistical Physics \textbf{83}, 291--357 (1996)

\bibitem{metzner2009}
Metzner, P., Sch{\"u}tte, C., Vanden-Eijnden, E.: Transition path theory for
  Markov jump processes.
Multiscale Modeling and Simulation \textbf{7}, 1192--1219 (2009)

\bibitem{mirebeau143D}
Mirebeau, J.M.: Anisotropic fast-marching on cartesian grids using lattice
  basis reduction.
 SIAM J. Numer. Anal. \textbf{52}, 1573--1599 (2014)

\bibitem{mirebeau14}
Mirebeau, J.M.: Efficient fast marching with Finsler metrics.
Numerische Mathematik \textbf{126}, 515--557 (2014)

\bibitem{mirebeau19}
Mirebeau, J.M., Portegies, J.M.: Hamiltonian fast marching: A numerical solver
  for anisotropic and non-holonomic eikonal PDEs.
 Image Process. Line \textbf{9}, 47--93 (2019)

\bibitem{nave10}
Nave, J.C., Rosales, R.R., Seibold, B.: A gradient-augmented level set method
  with an optimally local, coherent advection scheme.
Journal of Computational Physics \textbf{229}, 3802--3827 (2010)

\bibitem{Nolting16}
Nolting, B.C., Abbott, K.C.: Balls, cups, and quasi-potentials: quantifying
  stability in stochastic systems.
 Ecology \textbf{97}(4), 850--864 (2016)

\bibitem{EJMpaskal}
Paskal, N.: Github repository: Efficient jet marcher.
\url{https://github.com/npaskal/EfficientJetMarcher} (2021)

\bibitem{poppe18}
Poppe, G., Schaefer, T.: Computation of minimum action paths of the stochastic
  nonlinear schr{\"o}dinger equation with dissipation.
 Journal of Physics A: Mathematical and Theoretical \textbf{51}(33),
  335102 (2018)

\bibitem{potter20}
Potter, S.F., Cameron, M.K.: Jet marching methods for solving the eikonal
  equation.
 SIAM Journal on Scientific Computing \textbf{43}(6), A4121--A4146
  (2021)

\bibitem{sethian96}
Sethian, J.A.: A fast marching level set method for monotonically advancing
  fronts.
Proceedings of the National Academy of Sciences of the United States
  of America \textbf{93}(4), 1591--1595 (1996)

\bibitem{sethian99}
Sethian, J.A.: Level Set Methods and Fast Marching Methods Evolving Interfaces
  in Computational Geometry, Fluid Mechanics, Computer Vision, and Materials
  Science.
Cambridge University Press (1999)

\bibitem{sethian01}
Sethian, J.A., Vladimirsky, A.: Ordered upwind methods for static
  hamilton-jacobi equations.
Proceedings of the National Academy of Sciences of the United States
  of America \textbf{98}(20), 11069--11074 (2001)

\bibitem{sethian03}
Sethian, J.A., Vladimirsky, A.: Ordered upwind methods for static
  Hamilton-Jacobi-Bellman equations: theory and algorithms.
 SIAM Journal on Numerical Analysis \textbf{41}(1), 325--363 (2003)

\bibitem{talkner87}
Talkner, P.: Mean first passage time and the lifetime of a metastable state.
Zeitschrift f\"{u}r Physik B Condensed Matter \textbf{68}, 201--207
  (1987)

\bibitem{tao18}
Tao, M.: Hyperbolic periodic orbits in nongradient systems and
  small-noise-induced metastable transitions.
Physica D: Nonlinear Phenomena \textbf{363}, 1--17 (2016)

\bibitem{yang19}
Yang, S., Potter, S.F., Cameron, M.K.: Computing the quasipotential for
  nongradient SDEs in 3d.
Journal of Computational Physics \textbf{379}, 325--350 (2019)

\bibitem{zhou08}
Zhou, X., Ren, W., E, W.: Adaptive minimum action method for the study of rare
  events.
The Journal of Chemical Physics \textbf{128}(10), 104111 (2008)


\end{thebibliography}


\end{document}